\documentclass[11pt]{article}
\usepackage{amsmath, amsfonts, amsthm, amssymb,verbatim,graphicx,color,a4wide}
\usepackage{epsfig}

\newtheorem{claim}{Claim}
\numberwithin{equation}{section}
\usepackage{caption}
\usepackage{subcaption}
\newcommand \barv {\overline v}
\newcommand \bse {\begin{subequations}}
\newcommand \ese {\end{subequations}}

\newcommand \del  \partial

\newcommand{\auth}{\textsc}
\newcommand \bel {\begin{equation}\label}
\newcommand \be    {\begin{equation}}
\newcommand \ee    {\end{equation}}

\newcommand \eps     \epsilon
\newcommand \bei {\begin{itemize}}
\newcommand \eei {\end{itemize}}
\newcommand \sgn {\text{sgn}} 
\newcommand \map \mapsto 

\usepackage[most]{tcolorbox} 
\newtcolorbox{done}{%
     enhanced, breakable, size=minimal, colframe=white, parbox=false, after={\par}, 
     before upper={\indent}, colback=white, 
     overlay = {\draw[line width=2pt] (frame.north east) -|
                       ([xshift=3mm]frame.east)|-(frame.south east);},
     overlay first={\draw[line width=2pt] (frame.north east) -|
                           ([xshift=3mm]frame.south east);},
     overlay middle={\draw[line width=2pt] ([xshift=3mm]frame.north east) -- 
                              ([xshift=3mm]frame.south east);},
     overlay last={\draw[line width=2pt] ([xshift=3mm]frame.north east)|-
                          (frame.south east);},
}

\begin{document}

\title{Asymptotic structure of cosmological Burgers flows 
in one and 
\\
two space dimensions: a numerical study} 


\author{Yangyang Cao$^1$, Mohammad A. Ghazizadeh$^2$, and Philippe G. LeFloch$^1$}

\date{}  

\footnotetext[1]{Laboratoire Jacques-Louis Lions \& Centre National de la Recherche Scientifique,
Sorbonne Universit\'e, 4 Place Jussieu, 75252 Paris, France.
Email: {\sl caoy@ljll.math.upmc.fr, contact@philippelefloch.org.}
\newline
$^2$ Department of Civil Engineering, University of Ottawa, Ottawa, ON K1N 6N5, Canada. 
Email: {\sl sghaz023@uottawa.ca} 
\hfill Completed in July 2019.
}

\maketitle

\begin{abstract} 
We study the cosmological Burgers model, as we call it, which is a nonlinear hyperbolic balance law (in one and two spatial variables) posed on an expanding or contracting background. We design a finite volume scheme that is fourth-order in time and second-order in space, and allows us to compute weak solutions containing shock waves. Our main contribution is the study of the asymptotic structure of the solutions as the time variable approaches infinity (in the expanding case) or zero (in the contracting case). We discover that a saddle competition is taking place which involves, on one hand, the geometrical effects of expanding or contracting nature and, on the other hand, the nonlinear interactions between shock waves. 
\end{abstract}

\setcounter{tocdepth}{1}
\tableofcontents


\section{Introduction}

\paragraph*{The balance law of interest.}

We investigate numerically the global dynamics of a compressible fluid containing shock waves and evolving on a curved background spacetime of a contracting or expanding type. Motivated by the (inviscid) Burgers equation that has played such a central role in standard fluid dynamics, we consider here its relativistic version 
\bel{bueq30-2D}
a \, v_t + f(v)_x + g(v)_y  = a_t h(v),  \qquad (x, y) \in [0,L]^2, 
\ee
which we refer to as the {\bf cosmological Burgers model}. This equation provides a simple setup for designing and testing shock-capturing schemes in a curved spacetime background and investigating the asymptotic behavior of weak solutions; see \cite{PLF-review} for a derivation and review of such models. 

In \eqref{bueq30-2D}, the unknown is a function $v = v(t,x,y) \in (- 1 /\eps, 1/\eps)$ representing the main velocity component of a fluid vector field, and  $1/\eps$ represents the speed of light.
The fluxes $f=f(v)$ and $g=g(v)$ and the source function $h=h(v)$ are given smooth functions. We formulate the evolution on the domain $[0,L]^2$ with vanishing boundary conditions. 
A typical choice of flux and source functions is 
\bel{1Dflux}
f(v) = g(v) = {1 \over 2} v^2, \qquad h(v) = - v (1- \eps^2 v^2),
\ee
which allows us to recover the standard Burgers equation by taking the limit
 $a \to 1$ and $\eps \to 0$.  


\paragraph*{The geometric background of interest.} 

The function $a = a(t)>0$ describes a geometric background of contracting or expanding type.
Shock wave solutions to
nonlinear hyperbolic equations such as \eqref{bueq30-2D} are only defined in the {\sl forward time direction}
and, since the equation is singular at $t = 0$, it is natural to distinguish between two initial value problems corresponding to the following range of the time variable:  

\bei

\item In the range $t \in [t_0, + \infty)$, the background is assumed to be expanding toward the future in the sense that 
$a(t)$ increases monotonically to $+ \infty$ and data are prescribed at some $t_0>0$. 

\item In the range $t \in [t_0,0)$, the background is  assumed to be contracting toward the future in the sense that
$a(t)$ decreases monotonically to $0$ and data are prescribed at some $t_0<0$. 

\eei
\noindent A  typical choice is the function 
$a(t) = a_0 (t / t_0)^\alpha$, which we can 
normalize by taking $a_0 = 1$ and $t_0 = \pm 1$, in which $\alpha \in(0, 1)$ represents the rate of contraction or expansion of the background:
\bel{equa-talpha} 
a(t) = |t|^\alpha.
\ee 

Our model is motivated from the full Euler system posed on 
the so-called FLRW background (after Friedmann--Lema\^{i}tre--Robertson--Walker) describing 
a homogeneous and isotropic cosmology, for which a typical exponent is $\alpha = 2/3$.


\paragraph*{The strategy of this paper.}
  
We introduce a shock-capturing, high-order finite volume method for computing the weak solutions to \eqref{bueq30-2D}. Our numerical algorithm is sufficiently robust and accurate in order to investigate the propagation and nonlinear interaction of shock waves in presence of the curved geometry of interest. Our main challenge is then to determine the {\sl asymptotic behavior} of the flow which will turn out to be highly complex, both in the expanding and the contracting regimes. 

We recall that the inviscid Burgers equation has played a central role in the development of shock-capturing schemes in non-relativistic fluid dynamics. More recently, a generalization of the standard Burgers equation has been introduced and investigated on curved spacetimes by LeFloch and collaborators \cite{ALO,PLF1,PLFM,LMO} who took into account various geometrical effects.


We are going to discretize \eqref{bueq30-2D} via the finite volume methodology by keeping the structure of the equation at the discrete level. The numerical algorithm proposed below enjoys the following features:

\bei

\item {\sl Consistency with the divergence part.} Our scheme is consistent with the divergence part of the balance law and, therefore according to the Lax-Wendroff theorem, correctly computes weak solutions containing shock waves.

\item {\sl Second-order accuracy in space.} This is achieved by introducing a piecewise linear approximation and a min-mod limiter in order to prevent oscillations (similar to the Gibbs phenomena) in the vicinity of discontinuities of the solutions. This is an essential property for an accurate computation of shock waves in fluid flows. 

\item {\sl Fourth-order accuracy in time.}  A very high accuracy in time turned out to be important in the present context, since the background geometry may become singular as time evolves and we are interested in accurately computing the long-time asymptotics of the solutions. We rely here on a fourth-order Runge-Kutta discretization in order to achieve the desired accuracy. 
 
\end {itemize}


\paragraph*{Outline of this paper and main results.}

Our main contribution is a study of the asymptotic behavior of the solution as the time variable approaches infinity (in the expanding case) or approaches zero (in the contracting case). We discover that a competition is taking place which involves, on one hand, the geometrical effects of expanding or contracting nature and, on the other hand, the nonlinear interactions between shock waves. 

This paper is organized as follows. In Section~\ref{sec:BurgersEq}, we describe some properties of the cosmological Burgers model and describe the class of spatially homogeneous solutions. 
In Section~\ref{sec:alg1D},  working in the so-called cosmological time (denoted by $\tau$ below), we design a finite volume scheme for the $(1+1)$--cosmological Burgers equation which has the desired accuracy in space and in time. 
Next, in Section~\ref{section--4} we investigate the global dynamics of $(1+1)$--cosmological Burgers flows: 
in the expanding case, the fluid is coming to a rest in the late-time limit $\tau \to +\infty$ and, interesting, our scheme is sufficiently robust in order to capture a rescaled version of the solution which describe the small-scale features in this flow: we discover that the solution approaches an N-wave profile containing {\sl finitely many jumps} that no longer interact together in this late-time limit. This is reminiscent of phase transition phenomena. 
Analogous conclusions are then reached for the $(1+1)$--equation in the contracting case, and next in Section~\ref{FVforBurgers} for the same problems but now posed two spatial dimensions. A generalization of our method and numerical experiments to the full Euler systems of compressible fluids is presented in the companion paper \cite{CGL2}.


\section{Cosmological Burgers flows}
\label{sec:BurgersEq}

\paragraph*{A rescaled time variable.}

In this section, we describe various properties of the $(2+1)$ cosmological Burgers model given \eqref{bueq30-2D}.
It is interesting to introduce a new time variable, denoted by $\tau$ so that, after setting, 
\bel{attau}
a_t(t) = m(\tau) \text{ with } a(t) d\tau = dt, 
\ee
the balance law \eqref{bueq30-2D} in terms of the unknown $v = v(\tau, x, y) \in (-1,1)$ reads 
\bel{RBEat1}
v_\tau + f(v)_x + g(v)_y  = m(\tau) h(v), 
\qquad 
\tau \neq 0, \quad x,y \in [0,L], 
\ee
and in the following it will be convenient to formulate our numerical scheme in this time variable. Recall that we distinguish between two cases:  

\bei

\bse
\item In the expanding case, we have $t \in [1, + \infty)$, and our typical function is 
$a(t) = t^\alpha$ (with $\alpha \in (0, 1)$). With 
\bel{taut}
\tau = {t^{1-\alpha} \over 1-\alpha}, 
\qquad
m(\tau) = {\kappa \over \tau},
\qquad
\kappa = {\alpha \over 1- \alpha} \in (0, +\infty), 
\ee
the equation \eqref{RBEat1} with flux \eqref{1Dflux} and geometry function \eqref{equa-talpha} reads  
\bel{RstBEat}
v_\tau + {1 \over 2}(v^2)_x + {1 \over 2}(v^2)_y = - {\kappa \over  \tau} v (1-v^2),
\qquad 
\tau \in [\kappa+1, + \infty). 
\ee 
\ese

\bse
\item In the contracting case, we have $t \in [-1,0)$ and our typical function is 
$a(t) = (- t)^\alpha$ ($\alpha \in (0, 1)$). With 
\bel{taut1}
\tau = - {(- t)^{1-\alpha} \over 1-\alpha}, 
\qquad
m(\tau) =  {\kappa \over \tau}, 
\qquad
\kappa = {\alpha \over 1- \alpha} \in (0, +\infty), 
\ee
the equation \eqref{RBEat1} with flux \eqref{1Dflux} and geometry function \eqref{equa-talpha} reads 
\bel{RstBEat1}
v_\tau + {1 \over 2}(v^2)_x + {1 \over 2}(v^2)_y = - {\kappa \over  \tau} v (1-v^2),
\qquad 
\tau \in [-\kappa-1, 0). 
\ee 
\ese

\eei


\paragraph*{The non-relativistic limit.}

In the limit  $\eps \to  0$ when \eqref{1Dflux} is assumed, the balance law becomes 
$a \, v_t + {1 \over 2} (v^2 )_x + {1 \over 2} (v^2 )_y + v a_t = 0$ and, therefore,  
\bel{bueq4}
\big( a \, v \big)_t + {1 \over 2}(v^2)_x + {1 \over 2}(v^2)_y = 0. 
\ee
This is a {\sl conservation law} and, in fact, a weighted version of the standard Burgers equation.


\paragraph*{Conservation form for regular solutions.}

For {\sl sufficiently regular} solutions, our balance law can be transformed to a {\sl conservation law}, namely: 
\bel{1cl}
\bigg({\tau^\kappa v\over (1-v^2)^{1/2}}\bigg)_\tau + \bigg( {\tau^\kappa \over (1-v^2)^{1/2}}\bigg)_x + \bigg( {\tau^\kappa \over (1-v^2)^{1/2}}\bigg)_y= 0.
\ee 
However, this transformation is not valid for {\sl weak} solutions and, therefore, we will not use it in the following. 


\paragraph*{Spatially homogeneous solutions.}

Spatially homogeneous solutions are solutions $v=v(\tau)$ depending on the time variable only. Such solutions are relevant in describing the long time behavior of solutions, and are characterized by the ordinary differential equation 
\bel{BODE10}
v_\tau =- m(\tau) \, v (1-v^2),
\ee
equivalent to
$\big({1 \over v} + {v \over 1- v^{2}} \big) v_\tau = - m(\tau)$. 
Given any $v_0 \in (-1, 1)$, the solution $v = v(\tau)$ satisfying the initial condition $v(\tau_0) = v_0$ is given explicitly by 
$$
{e^{M(\tau)} v (\tau) \over (1 - v(\tau)^2)^{1/2}} 
= {e^{M(\tau_0)} v_0 \over (1-v_0^2)^{1/2}}, 
\qquad
M(\tau) = \int^\tau m(s)ds
$$
or, equivalently, 
\bel{eq:sol45} 
v(\tau) = {v_0 \over \sqrt{v_0^2 + (1 - v_0^2) e^{ 2 \int_{\tau_0}^\tau m(s) ds}}}. 
\ee

\bse

\begin{figure}[htbp]
\centering
  \subcaptionbox{}{\includegraphics[width=2.8in, trim={0 0 0 0 1cm},clip]{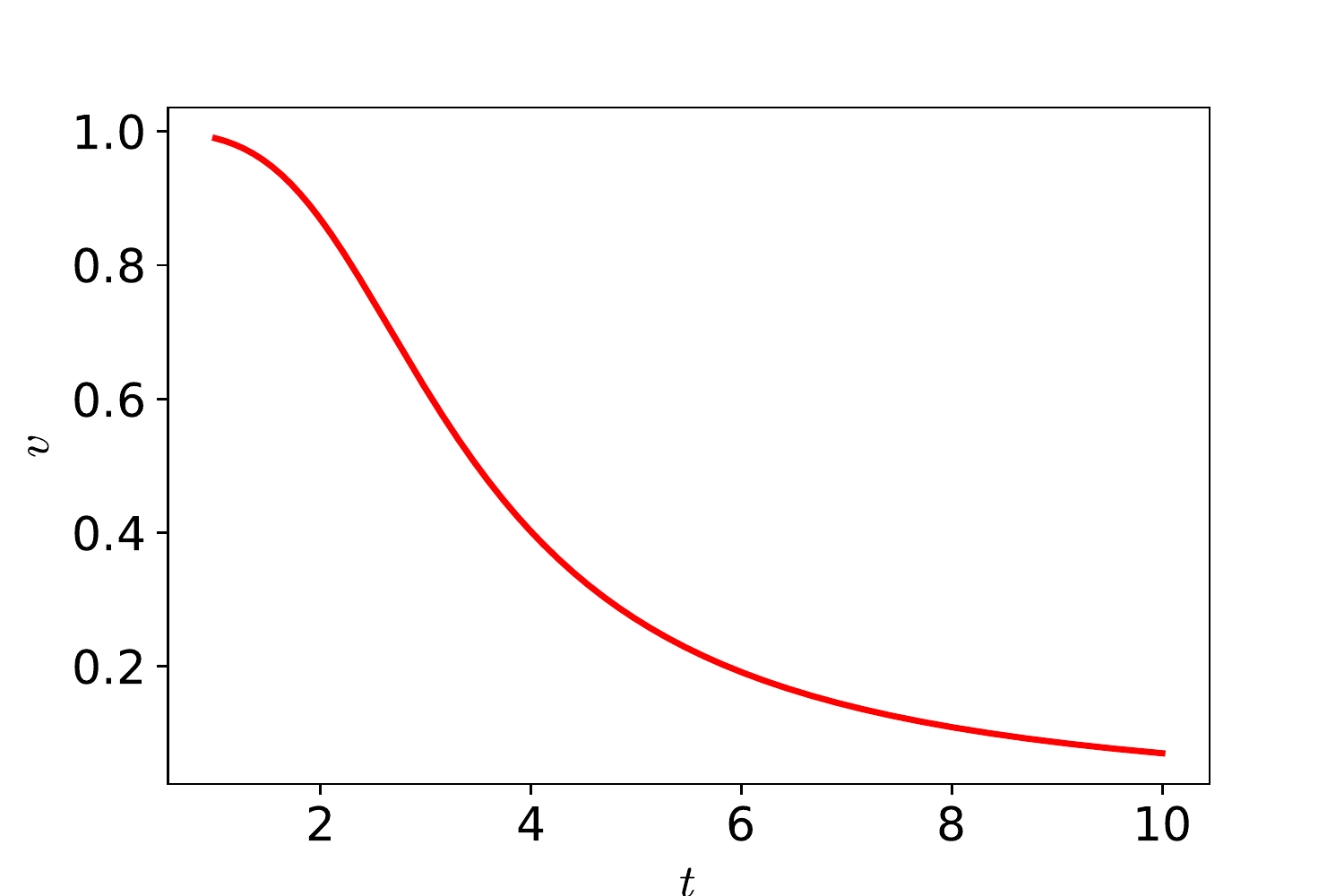}}
  \subcaptionbox{}{\includegraphics[width=2.8in, trim={0 0 0 0 1cm},clip]{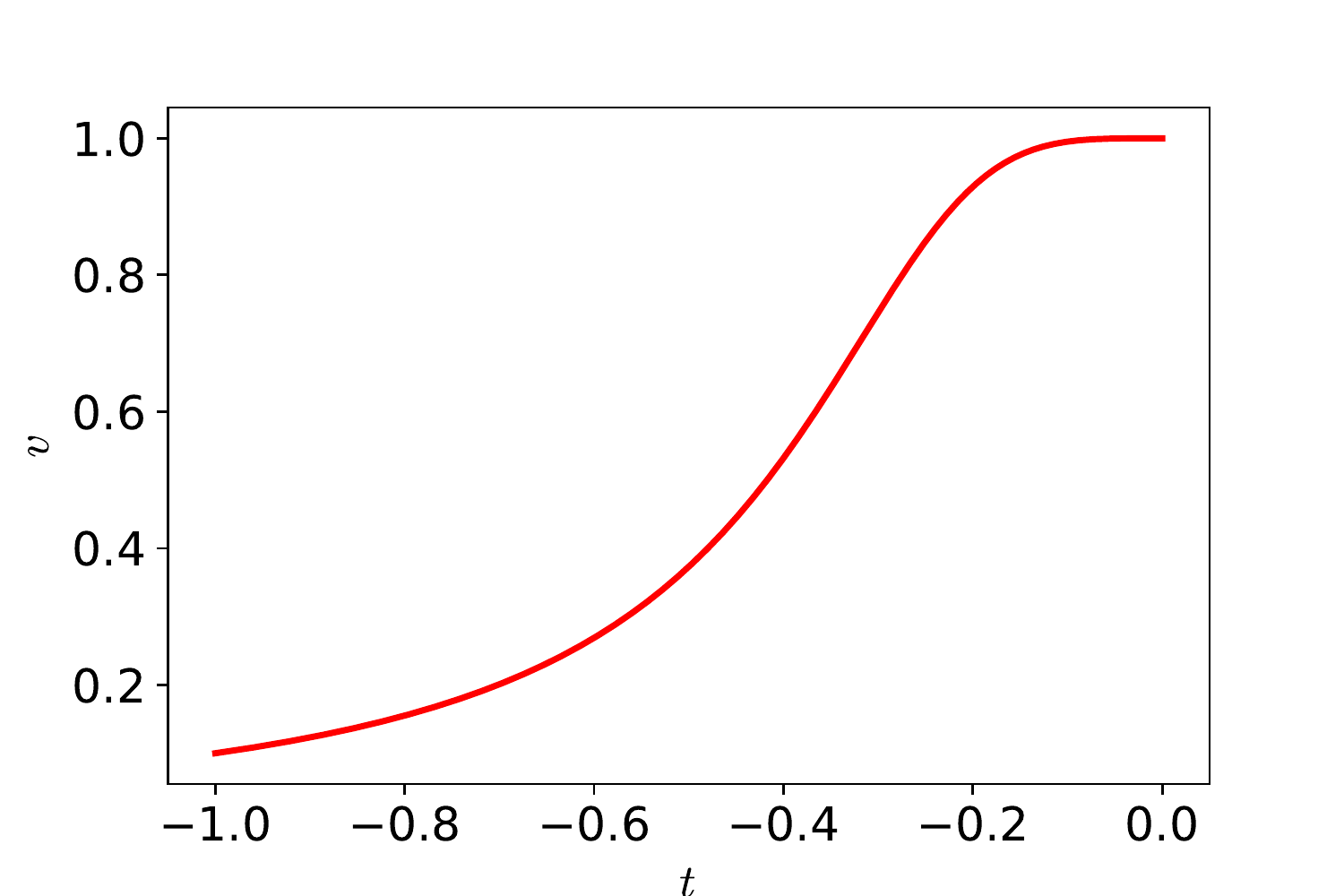}}
\caption{(a) Homogenous solution on an expanding background. (b) Homogenous solution on a contracting background.}
\label{FIG-Homo}
\end{figure}

We can specialize our conclusion above to the case $m(\tau) = {\kappa \over \tau}$ on an expanding background  and we find 
$M(\tau) = \log \tau^\kappa$, so that 
\bel{vtau11}
v(\tau) = {v_0 \over \sqrt{v_0^2 + \tau^{2 \kappa} (1 - v_0^2)}} \qquad \text{(expanding case).} 
\ee
On a contracting background with the function $m(\tau) = {\kappa \over \tau}$ we find
$M(\tau) = \log (-\tau)^\kappa$, and the spatially homogeneous solutions are  
\bel{vtau12}
v(\tau) = {v_0 \over \sqrt{v_0^2 + (-\tau)^{2 \kappa} (1 - v_0^2)}} \qquad \text{(contracting case).} 
\ee
\ese
Figure~\ref{FIG-Homo} contains a plot of the spatially homogeneous solutions, which clearly enjoy the following properties: 

\bei
\bse

\item All spatially homogeneous solutions satisfy $|v| <1 $, as required.  

\item On an expanding background $\tau \to + \infty$, $v(\tau) \simeq \pm \tau^{- \kappa}$ (up to a positive multiplicative constant); thus, the solution converges to $0$: 
\be
\lim_{\tau \to + \infty} v(\tau) = 0 \quad \text{(homogeneous solutions on an expanding background).}
\ee

\item On a contracting background $\tau \to 0$, $\pm 1 + v(\tau) \simeq \pm (-\tau)^{2 \kappa}$
(up to a positive multiplicative constant). Therefore, 
\be
\lim_{\tau \to 0 \atop \tau < 0} v(\tau) = \pm 1
\quad \text{(homogeneous solutions on a contracting background).}
\ee
\ese

\eei


\section{A finite volume scheme for $(1+1)$--cosmological Burgers flows} 
\label{sec:alg1D}
 
\subsection{The first-order Godunov discretization}

We begin with $(1+1)$-dimensional equations and present a discretization of the cosmological Burgers model 
\bel{bueq300}
v_\tau + f(v)_y = m(\tau)\, h(v), \qquad y \in [0,L], 
\ee
when an initial value is specified at some time $\tau = \tau_0$
\bel{1Dv0}
v_0(y) = v(\tau_0, y), \qquad y \in [0,L]. 
\ee
For definiteness, we write our scheme for the expanding case where $\tau_0>0$. 
We follow the finite volume methodology and a time-length $\Delta \tau$ is introduced together with 
the discrete times $\tau_n = \tau_0 + n \Delta \tau$ for $n=0, 1, \ldots$, as well as 
a space-length $\Delta y$ and discrete spatial points $y_j = j \Delta y \in [0,L]$ 
and $y_{j +1/2} = (j +1/2) \Delta y \in [0,L]$ (for a suitable range of integers $j$).

Using the notation
\bel{averV}
\aligned
& v_j^n \simeq {1 \over \Delta y}  \int ^{y_{j +1/2}}_{y_{j -1/2}} v(\tau_n, y)dy,
\qquad
s_j^n \simeq {1 \over \Delta y \Delta \tau}  \int ^{y_{j +1/2}}_{y_{j -1/2}} \int_{\tau_n}^{\tau_{n+1}} m(\tau)\, h(v) \, dy d\tau, 
\endaligned
\ee
we intergate the balance law \eqref{bueq300} on the slab $[\tau_n, \tau_{n+1}] \times [y_{j-1/2}, y_{j+1/2}]$ and obtain
\bel{firstg}
v^{n+1}_j = v^n_j - {\Delta \tau \over \Delta y} \Big(f^n_{j+1/2} - f^n_{j-1/2}\Big) + \Delta \tau s_j^n. 
\ee


For a first-order approximation (at this stage) we define the source $s_j^n$ to be 
\bel{sdd}
s_j^n =  m(\tau_n)\, h(v_j^n),
\ee
and for the numerical flux we set 
\bel{NFLux} 
f^n_{j+1/2} = f(v^n_j, v^n_{j+1}), 
\ee
in which for the two-point flux $f=f(v,w)$ we can choose, for instance, the Godunov flux $f_G$ is determined by solving the Riemann problem. Specifically, for any convex flux function such as \eqref{1Dflux} and thus satisfying the normalization
\be
f(0) = f'(0) = 0, 
\ee
we have the following explicit expressions ($v_l, v_r \in [-1,1]$):
\bei

\bse
\label{equa388}

\item Case $v_l > v_r$:
\bel{Shockflux}
f_G(v_l^n, v_r^n) 
= \begin{cases} 
f(v_l^n),   
\,
& f(v_r^n) - f(v_l^n) \leq 0,
\\
f(v_r^n),
& f(v_r^n) - f(v_l^n) \geq 0,
\\
0,  
& \text{otherwise}.
\end{cases}
\ee

\item Case $v_l^n \leq v_r^n$:
\bel{Rereflux}
f_G(v_l^n, v_r^n) 
= \begin{cases} 
f(v_l^n),   
\,
& f'(v_l^n) > 0,
\\
f(v_r^n),
&  f'(v_r^n) < 0,
\\
f(0),
&  \text{otherwise}.
\end{cases}
\ee

\ese

\eei


The following restriction on the time step is also imposed (which we express directly for our quadratic flux): 
\bei 

\item As far as the nonlinear propagation is concerned, we require that $\Delta \tau$ satisfies the so-called CFL (Courant-Friedrichs-Lewy) condition at any given $n$: 
\bel{CFL0}
{\Delta \tau  \over \Delta y} \max_j | v^n_j | < 1, 
\ee
together with further conditions in the expanding or contracting cases.

\item {\it Expanding background.} In this case, we expect that $v^n_j \to 0$ as $\tau \to +\infty$, so that 
$\Delta \tau$ is not restricted by \eqref{CFL0} for sufficiently large times. 
A second stability condition is required which is motivated by the following discretization
$v^{n+1} = v^n \Big( 1 - (\Delta \tau) (\kappa/\tau_n) \, (1 - (v^n)^2) \Big)$ 
of the ODE \eqref{BODE10}: 
$$
\max_j 
\Big|1- {\kappa \over \tau_n} \Delta \tau(1- (v^n_j)^2)\Big| <1,
$$
that is, 
$$
\Delta \tau \leq {2 \tau_n \over  \kappa (1- (v^n_j)^2)}.
$$
To summarize, the following stability condition is chosen on an expanding background:
\bel{CFL1}
\Delta \tau \leq \min \Bigg({\Delta y \over \max_j \big| v^n_j \big|}, {2 \tau_n \over  \kappa (1- (v^n_j)^2)}\Bigg).
\ee

\item {\it Contracting background.} In this case, we expect that 
$v^n_j \to \pm1$ as $\tau \to 0$,  so that $\Delta \tau$ would approach a constant if we would impose \eqref{CFL0}, 
whereas the time $\tau \in [-1, 0)$ is bounded above by $0$. 
It is natural to select time-increments that are approaching zero so that $\ldots < \tau^n < \tau^{n+1} < \ldots < 0$ only reach zero asymptotically as $n \to +\infty$. In the following paper, we propose to do so on a linear way with respect to $\tau_n$ and, in addition, to use time-increment that are proportional to $1/\kappa$ since a larger $\kappa$ means a stiffer ODE problem: 
$$
\Delta \tau \leq \min\big(1, {1\over \kappa}\big) |\tau_n|.
$$
Therefore, the following stability condition is required for a contracting background:
\bel{CFL1-con}
\Delta \tau \leq \min \Bigg({\Delta y \over \max_j \big| v^n_j \big|},  \min\big(1, {1\over \kappa}\big) |\tau_n|\Bigg).
\ee

\eei   


\subsection{Temporal discretization}

After integrating the equation \eqref{bueq300} over interval $[y_{j-1/2}, y_{j+1/2}]$ we arrive
 at the semi-discrete finite volume scheme
\bel{FVb0}
{d \over d\tau} v_j
= - {1 \over \Delta y}  \big(f_{j+1/2} - f_{j-1/2} \big) + m(\tau) \, h(v_j),
\ee
where
\[
f_{j-1/2} = f_G(v_{j-1},v_j), \qquad f_{j+1/2} = f_G(v_j, v_{j+1}),
\]
and the Godunov flux $f_G$ is defined by \eqref{equa388}. 
\bse
To shorten the notation, we introduce 
\bel{Gg}
G(v, \tau)_j = - {1 \over \Delta y} \big(f_{j+1/2} - f_{j-1/2}\big) + m(\tau) \, h(v_j),
\ee
and we rewrite \eqref{FVb0} as 
\bel{FVb11}
{d \over d\tau}  v_j = G(v, \tau)_j.
\ee
A fourth-order Runge-Kutta discretization from the initial data $v(\tau_0) = v_0$  is now applied.
We denote by $v_j^n$ the numerical solution given by \eqref{FVb11} at some time $\tau_n$, hence
\bel{FVb12}
{d \over d\tau}  v_j^n = G(v, \tau_n)_j,
\ee
and then we define
\bel{vk1234}
\aligned
& v_j^{0,n } = v_j^n,    
\quad
&& {K_j^{1,n}} 
= G\big(v^{0, n}, \tau_n\big)_j,
\\
& v_j^{1,n} 
= 
v_j^{0,n} + {\Delta \tau \over 2} K_j^{1,n},  
\quad
&&  K_j^{2, n} 
= 
G\big( v^{1, n}, \tau_n + {\Delta \tau \over 2} \big)_j,
\\
& v_j^{2, n} 
= 
v_j^{0, n} + {\Delta \tau \over 2} K_j^{2,n}, 
\quad  
&& K_j^{3, n} 
= 
G\big( v^{2,n}, \tau_n + {\Delta \tau \over 2} \big)_j,
\\
& v_j^{3, n} 
= 
v_j^{0,n} +  \Delta \tau K_j^{3,n},  
\quad
&&  K_j^{4, n} 
= 
G\big( v^{3, n}, \tau_n + \Delta \tau \big)_j,
\\
& v_j^{n+1} 
= 
v_j^{0, n} + {\Delta \tau \over 6}  \big({K_j^{1, n} + 2 K_j^{2, n}+ 2 K_j^{3, n} + K_j^{4, n}}\big).
\endaligned
\ee
\ese
Furthermore, the same stability conditions \eqref{CFL1} and \eqref{CFL1-con} are assumed.


\subsection{Second-order spatial discretization}
\label{Secondspace}

In order to improve the accuracy of the algorithm, we design a second-order version of our scheme, based on a piecewise linear reconstruction. We introduce the piecewise linear reconstruction
\bel{PLA}
v_j^n(y) = v_j^n + (y - y_j) \delta_j^n, 
\ee
where $\delta_j^n$ represents the local slope of the numerical solution in each cell. 
\bse
In order to prevent a Gibbs-type phenomena where oscillations would arise near discontinuities, the following limiter is 
applied:  
\bel{MCL}
\delta_j^n \Delta y 
= 
\begin{cases} 
\sgn(v_{j+1}^n - v_{j-1}^n)\min \Big( 2|v_j^n - v_{j-1}^n|,  2|v_{j+1}^n - v_j^n|, {1 \over 2} \big|v_{j+1}^n - v_{j -1}^n\big|\Big),   
\,
& \eta_j^n >0,
\\
0,
& \text{otherwise},
\end{cases}
\ee
in which we have set 
\bel{eta}
\eta_j^n = (v_{j+1}^n - v_j^n)(v_j^n - v_{j-1}^n).
\ee
The values of the reconstruction at the interfaces are denoted by $v_{j,L}^n$ and $v_{j,R}^n$, that is, 
\bel{Vinterface}
v_{j,L}^n = v_j^n - {\Delta y\over 2} \delta_j^n, 
\qquad
v_{j,R}^n = v_j^n + {\Delta y\over 2} \delta_j^n.
\ee
At each interface $y_{j+1/2}$, we apply our first-order scheme with the states replaced by 
the left- and right-hand values $v_{j, R}^n$ and $v_{j+1, L}^n$. This ensures that the algorithm provides a
 second-order approximation.  

\ese 


\section{Global dynamics of $(1+1)$--cosmological Burgers flows}
\label{section--4}

\subsection{Validation of the numerical algorithm}

\begin{figure}[htbp]
\centering
  \subcaptionbox{}{\includegraphics[width=2.8in, trim={0 0 0.0 1cm},clip]{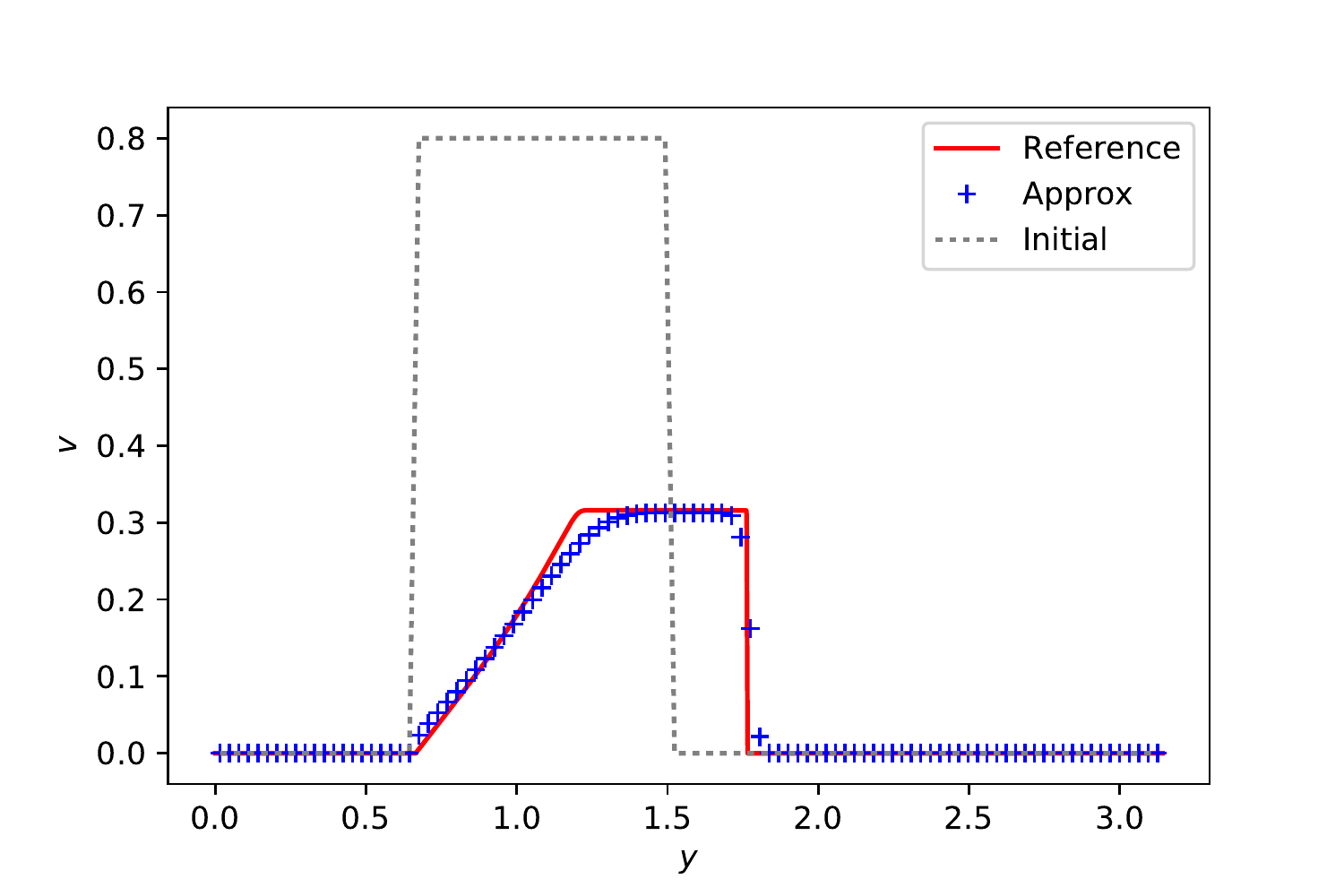}}
  \subcaptionbox{}{\includegraphics[width=2.8in, trim={0 0 0.0 1cm},clip]{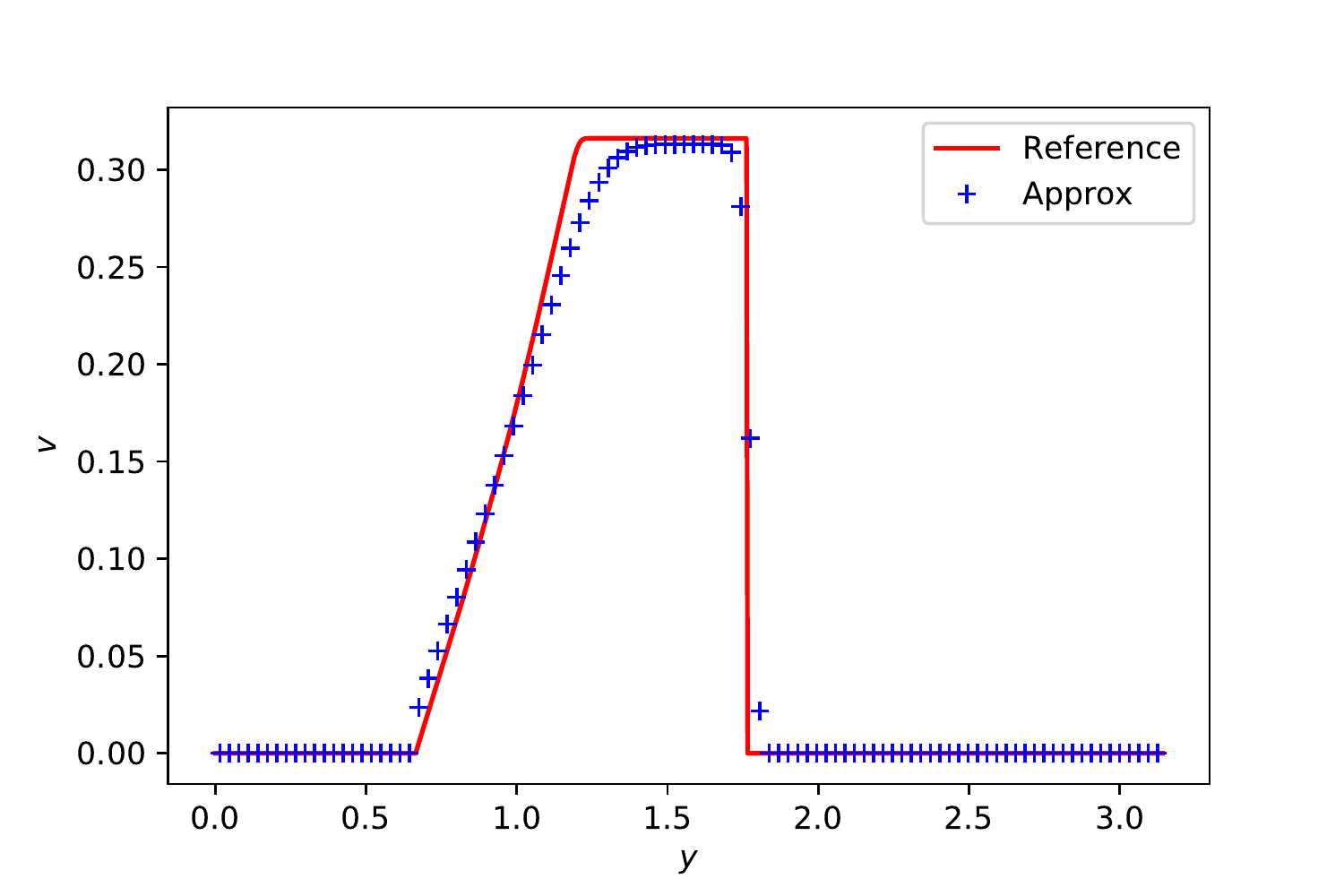}}
\caption{(a) Numerical solution: first-order method. (b) Magnified solution.}
\label{FIG-10}
\end{figure}

\begin{figure}[htbp]
\centering
  \subcaptionbox{}{\includegraphics[width=2.8in, trim={0 0 0.0 1cm},clip]{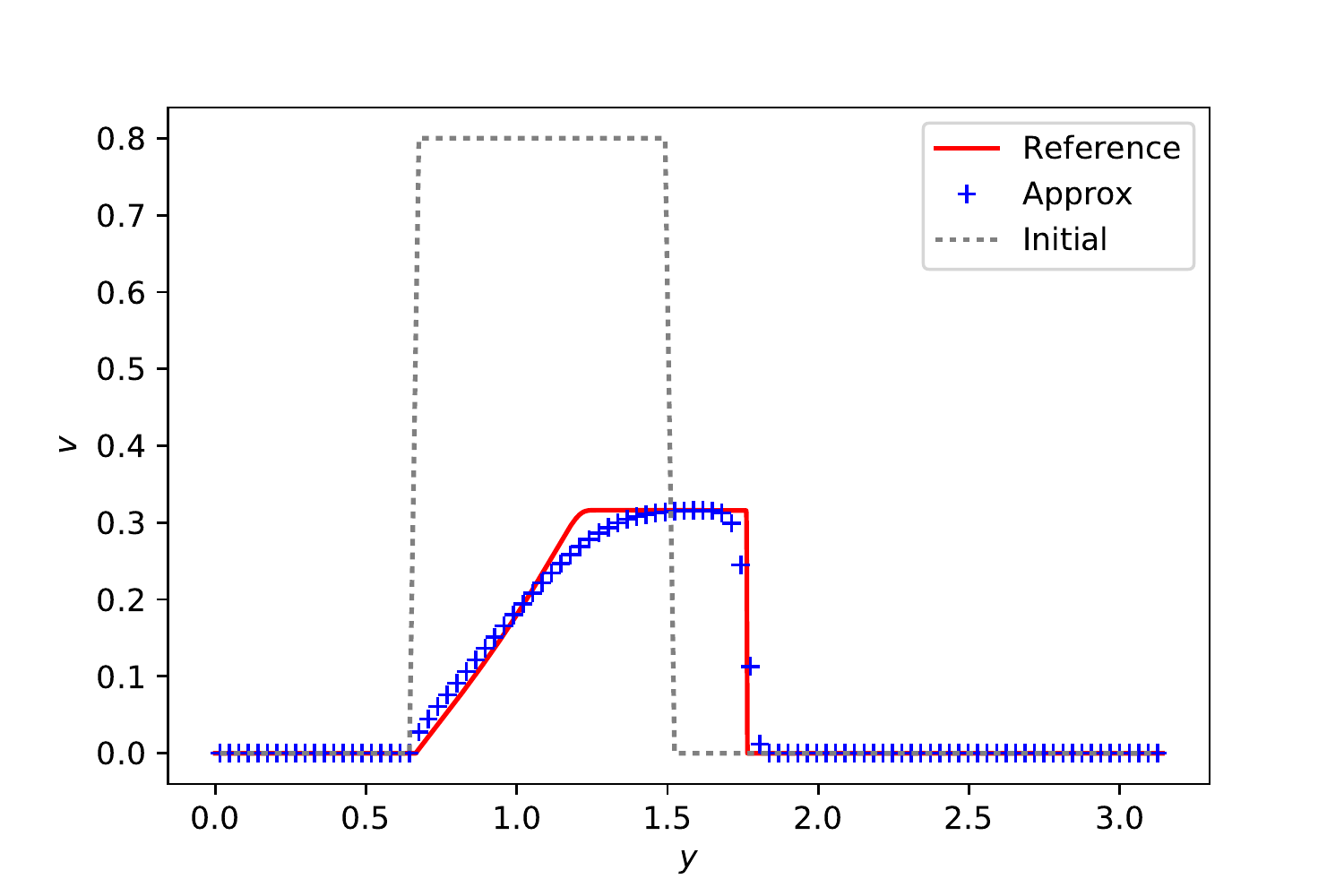}}
  \subcaptionbox{}{\includegraphics[width=2.8in, trim={0 0 0.0 1cm},clip]{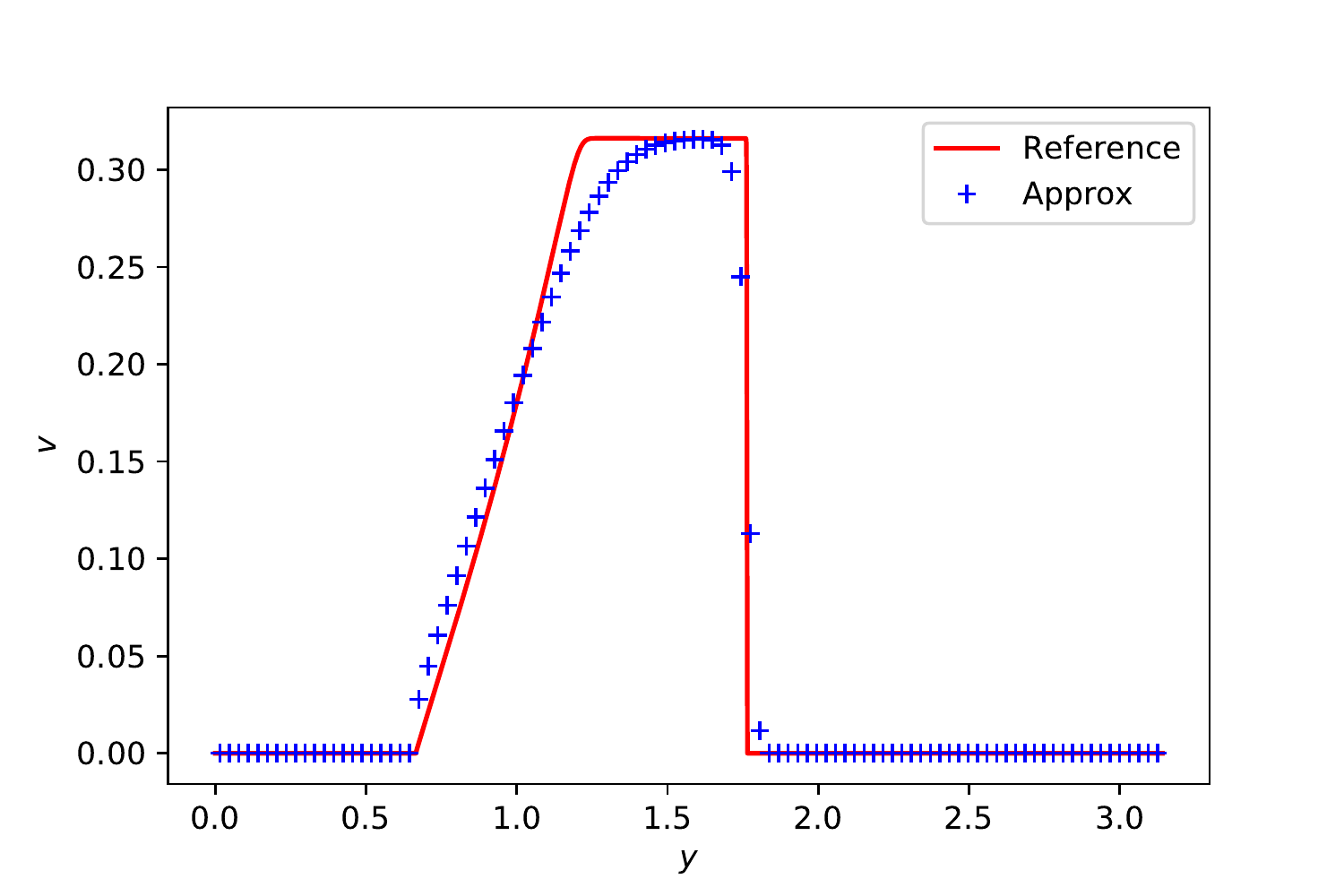}}
\caption{(a) Numerical solution:  first-order in space and Runge-Kutta in time. (b)  Magnified solution.}
\label{FIG-11}
\end{figure}

\begin{figure}[htbp]
\centering
  \subcaptionbox{}{\includegraphics[width=2.8in, trim={0 0 0.0 1cm},clip]{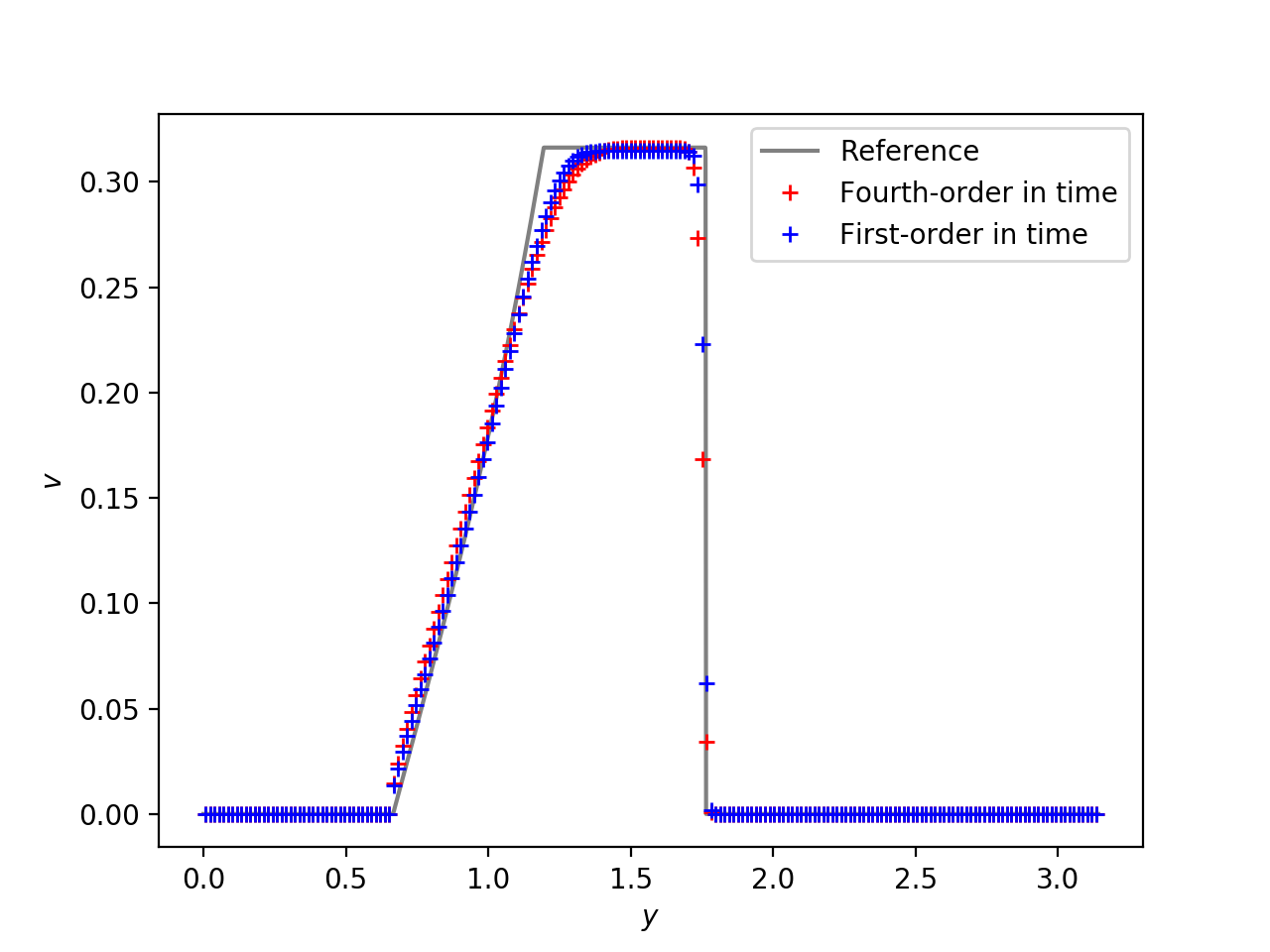}}
  \subcaptionbox{}{\includegraphics[width=2.8in, trim={0 0 0.0 1cm},clip]{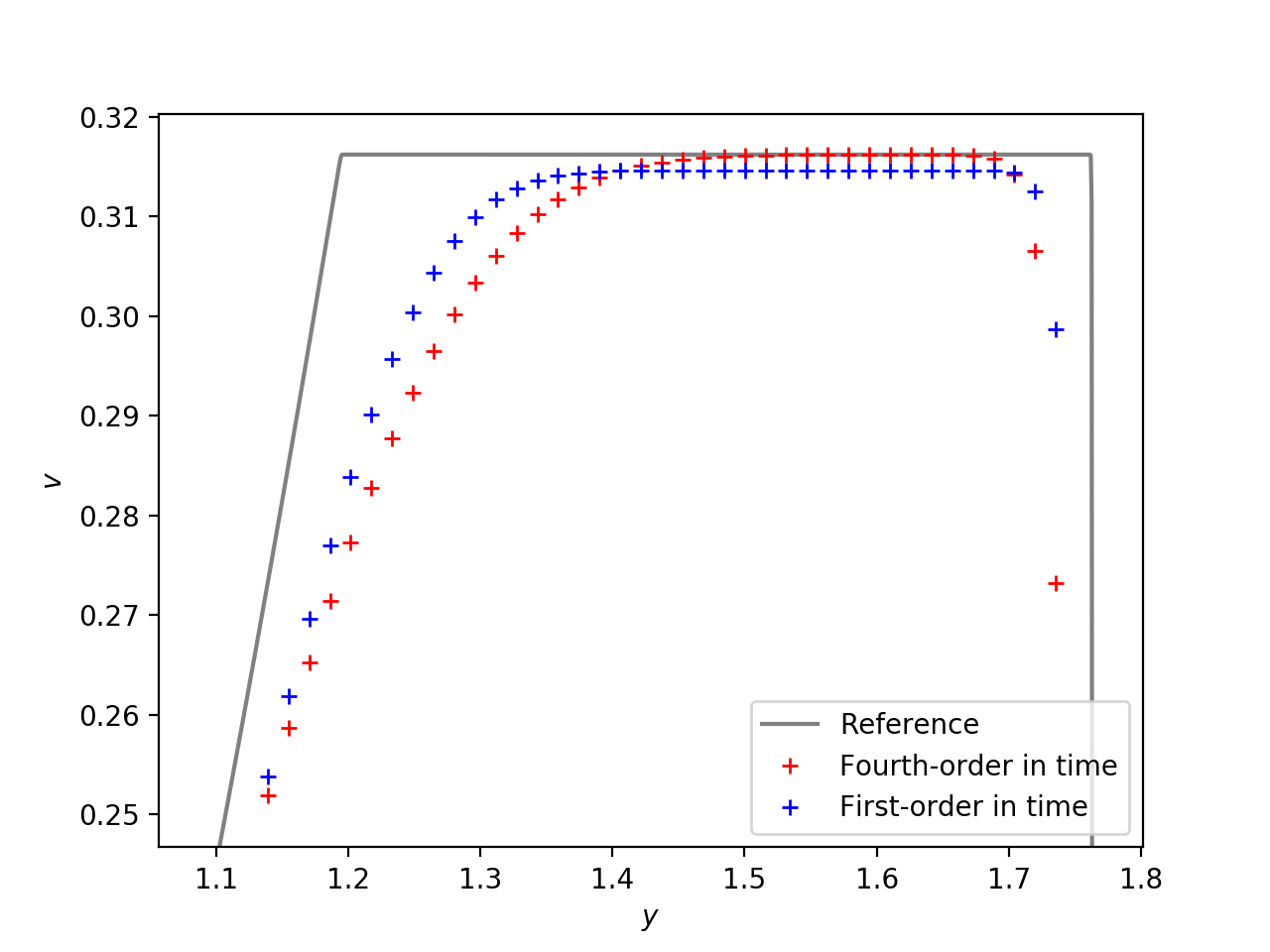}}
\caption{(a) First-order in space and first-order in time compared to first-order in space and fourth-order in time. (b)  Magnified solution.}
\label{FIG-140}
\end{figure}

\begin{figure}[htbp]
\centering
  \subcaptionbox{}{\includegraphics[width=2.8in, trim={0 0 0.0 1cm},clip]{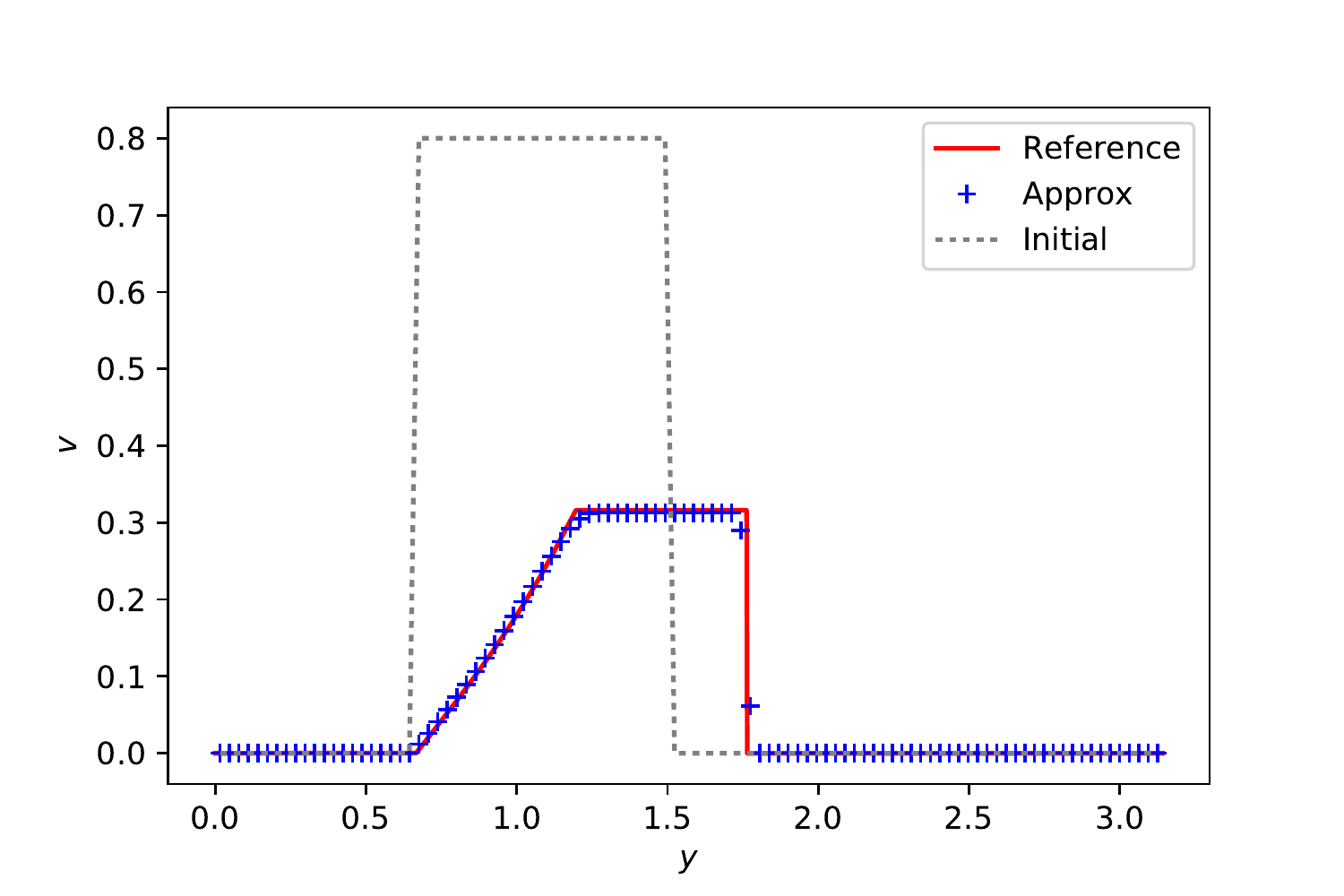}}
  \subcaptionbox{}{\includegraphics[width=2.8in, trim={0 0 0.0 1cm},clip]{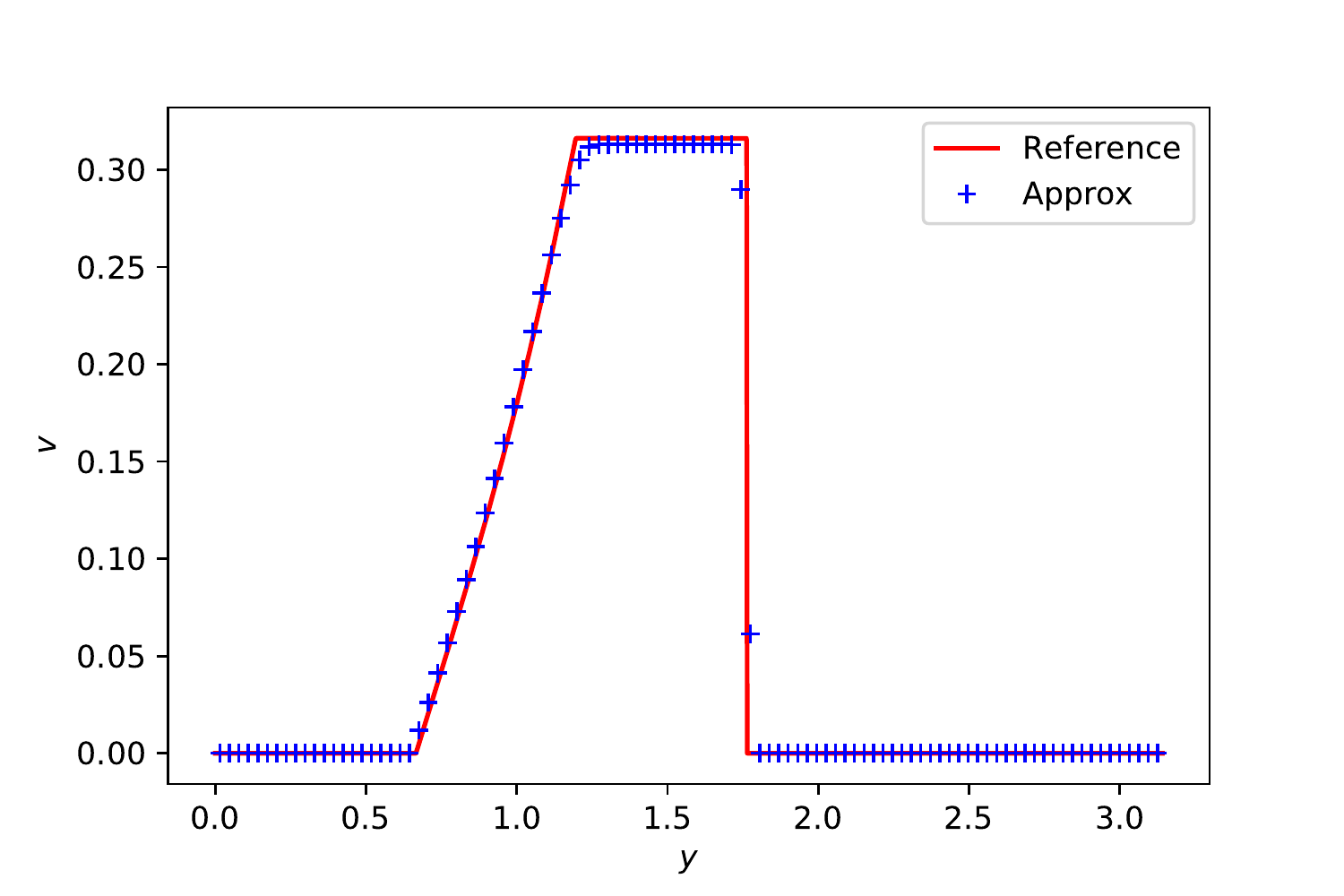}}
\caption{(a) Numerical solution: second-order scheme. (b)  Magnified solution.}
\label{FIG-12}
\end{figure}

\begin{figure}[htbp]
\centering
  \subcaptionbox{}{\includegraphics[width=2.8in, trim={0 0 0.0 1cm},clip]{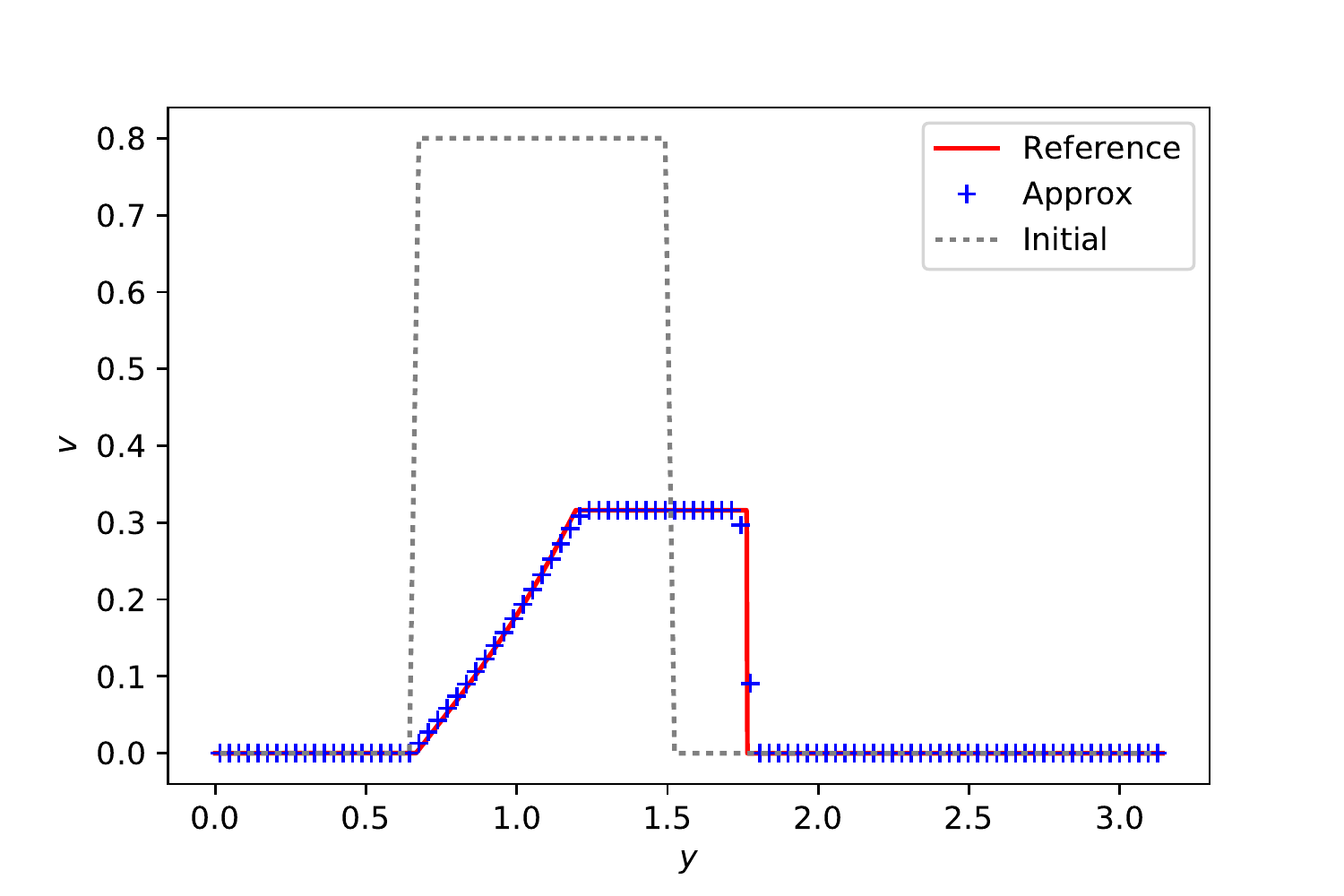}}
  \subcaptionbox{}{\includegraphics[width=2.8in, trim={0 0 0.0 1cm},clip]{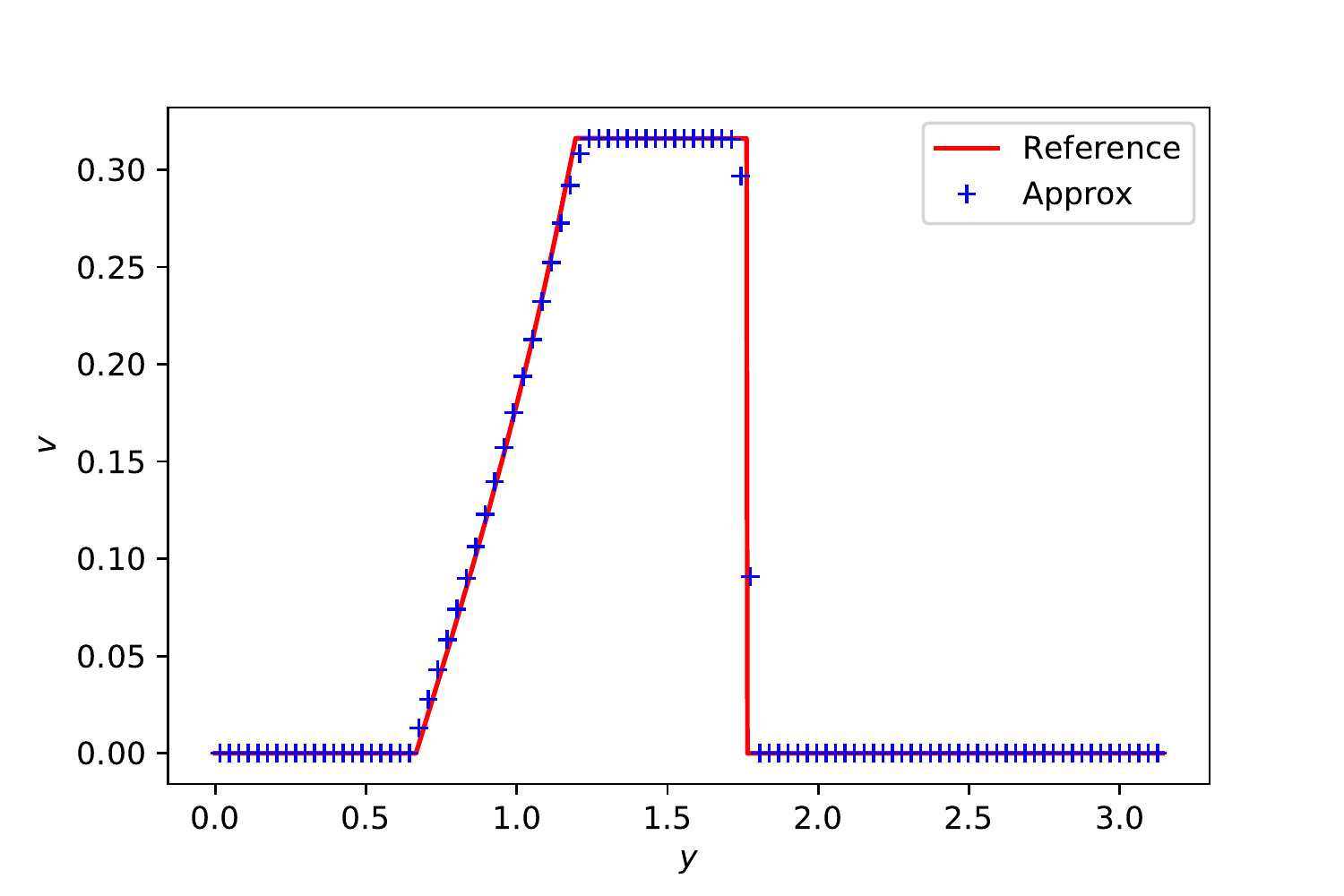}}
\caption{(a) Numerical solution: second-order in space and Runge-Kutta in time. (b)  Magnified solution.}
\label{FIG-13}
\end{figure}

We now perform several numerical experiments with the relativistic Burgers equation, expressed in the $\tau$-variable, that is,   
\bel{bueq40}
v_\tau + f(v)_y =m(\tau) h(v),
\qquad
y\in [0, \pi], 
\ee
with flux $f(v) = {1\over 2} v^2$ and source term
$h(v) = - v(1- v^2)$, while the geometric function is $m(\tau) = {\kappa \over \tau}$ with $\kappa >0$.
We begin with an initial data containing a single jump discontinuity, say, 
\bel{ina111}
v_0(y) = \begin{cases}0.8,   
\quad \,
& 0.666 \leq y < 1.5,
\\
0,
& \text{otherwise},
\end{cases}
\ee
We denote by $J$ the total number of grid cells in space, and $J = 5000$ is chosen in order have a very fine grid. The numerical solution will serve as a``reference solution'', since on such a grid is presumably very close to the exact solution. The numerical results are presented in Figures \ref{FIG-10}, \ref{FIG-11},  \ref{FIG-140}, 
\ref{FIG-12},  and \ref{FIG-13}, 
and are given at several order of accuracy: 
first-order in space and first-order in time;
first-order in space and fourth-order in time;
second-order in space and second-order in time;
second-order in space and fourth-order in time. 
We observe that the second-order in space and fourth-order in time discretization significantly provides the best possible accuracy for the solution. These results fully justify the involved construction we have proposed in the previous section.


\begin{figure}[htbp]
\centering
  \subcaptionbox{}{\includegraphics[width=2.8in, trim={0 0 0.0 1cm},clip]{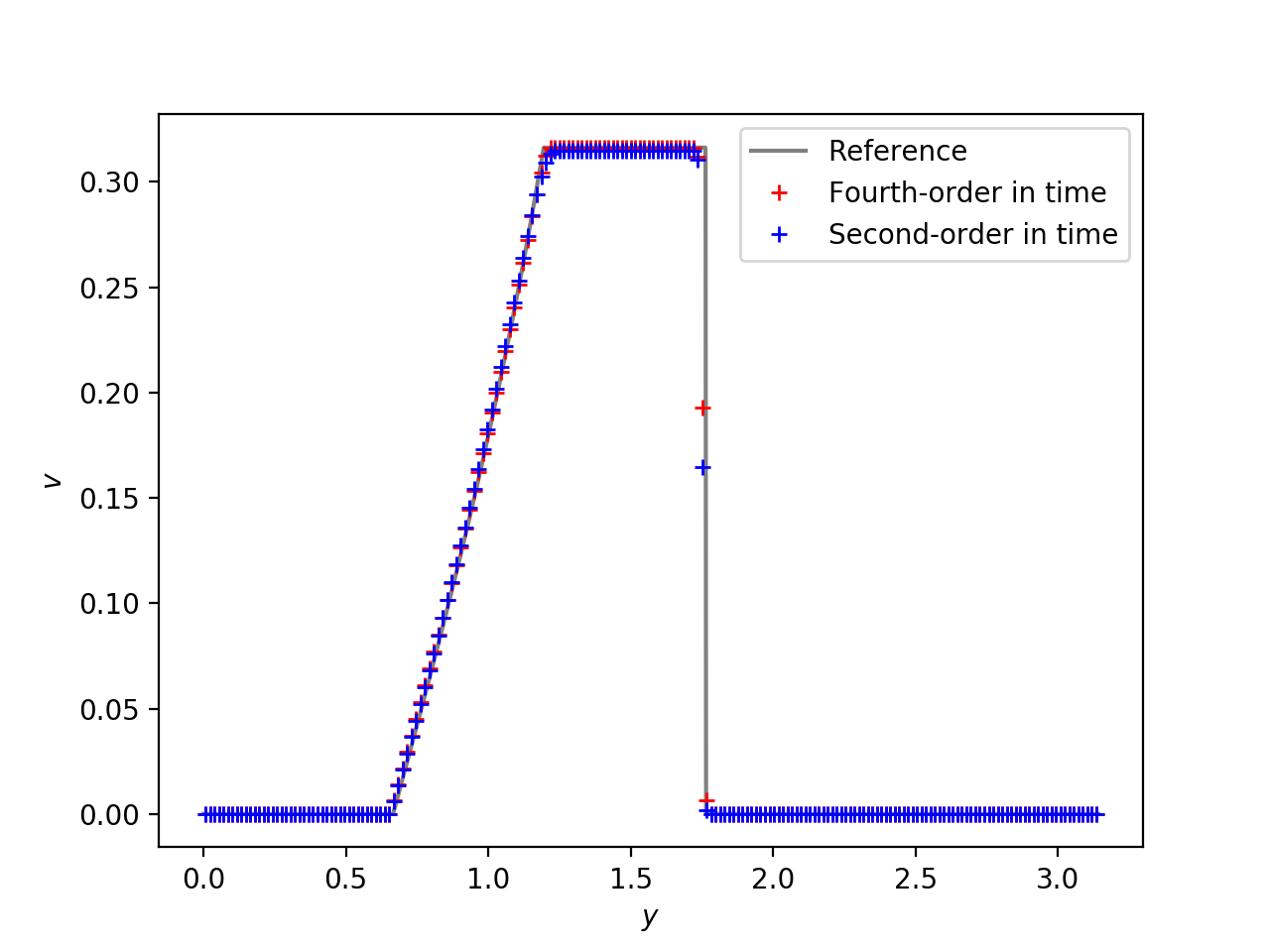}}
  \subcaptionbox{}{\includegraphics[width=2.8in, trim={0 0 0.0 1cm},clip]{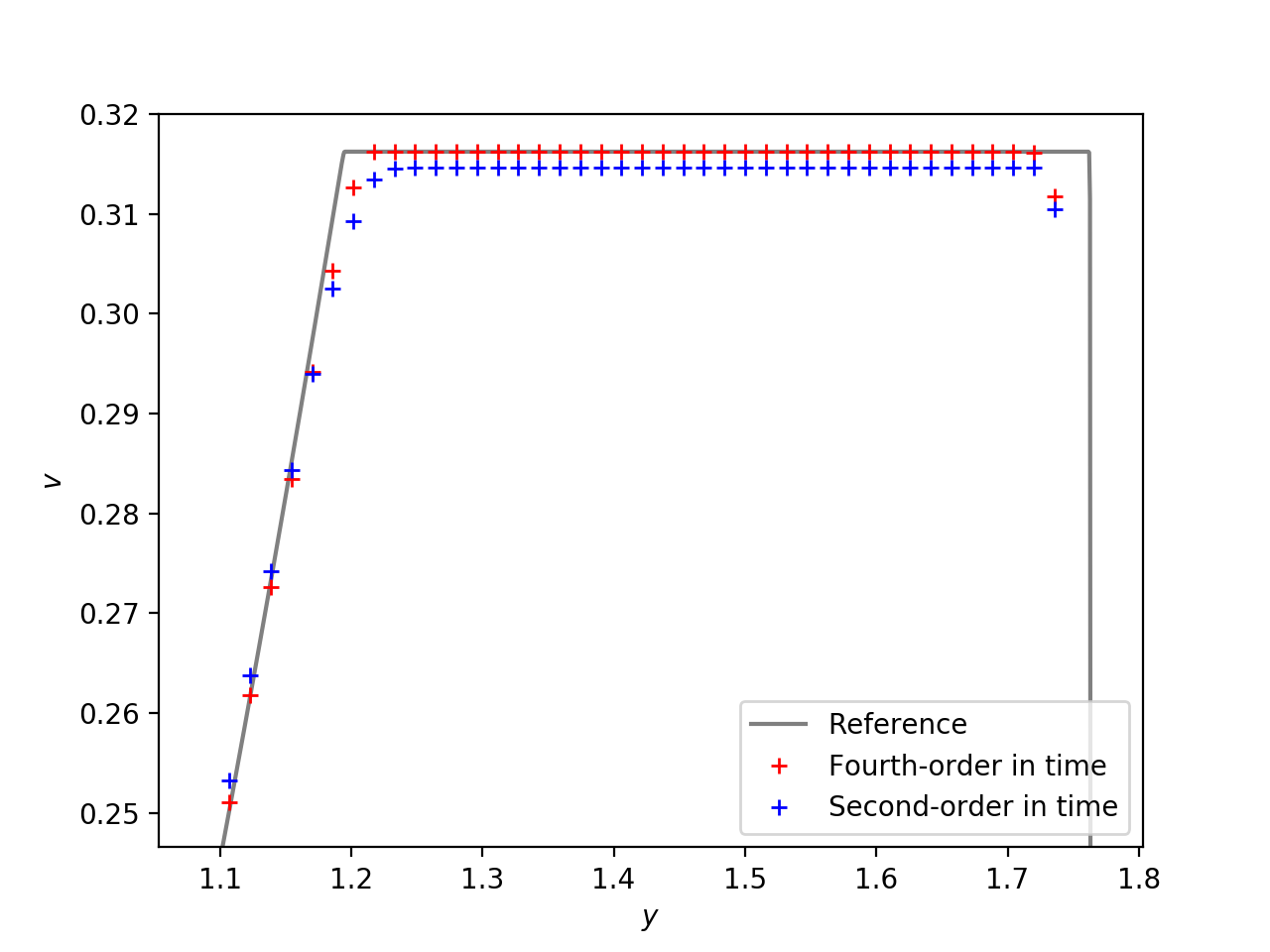}}
\caption{(a) Second-order in space second-order in time compared to second-order in space and fourth-order in time. (b) 
Magnified solution.}
\label{FIG-15}
\end{figure}

\begin{figure}[htbp]
\centering
 \includegraphics[width=3in, trim={0 0 0.0 1.1cm},clip]{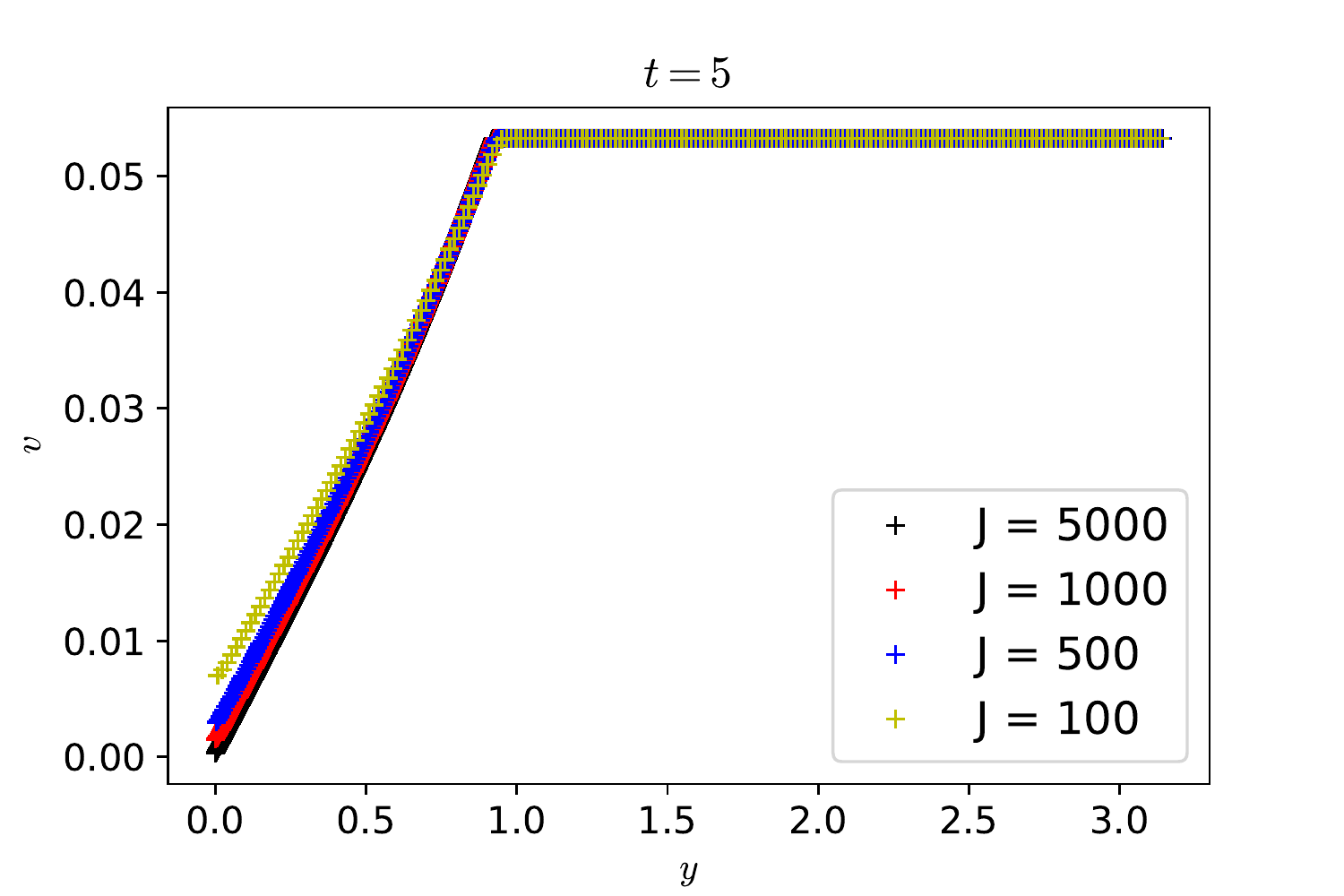}
\caption{Comparing solutions for different grid cells in space at the time $\tau=5$.}
\label{FIG-Constant1}
\end{figure}
\begin{figure}[htbp]
\centering
  \subcaptionbox{$\tau=2$}{\includegraphics[width=2.8in, trim={0 0 0.0 1.1cm},clip]{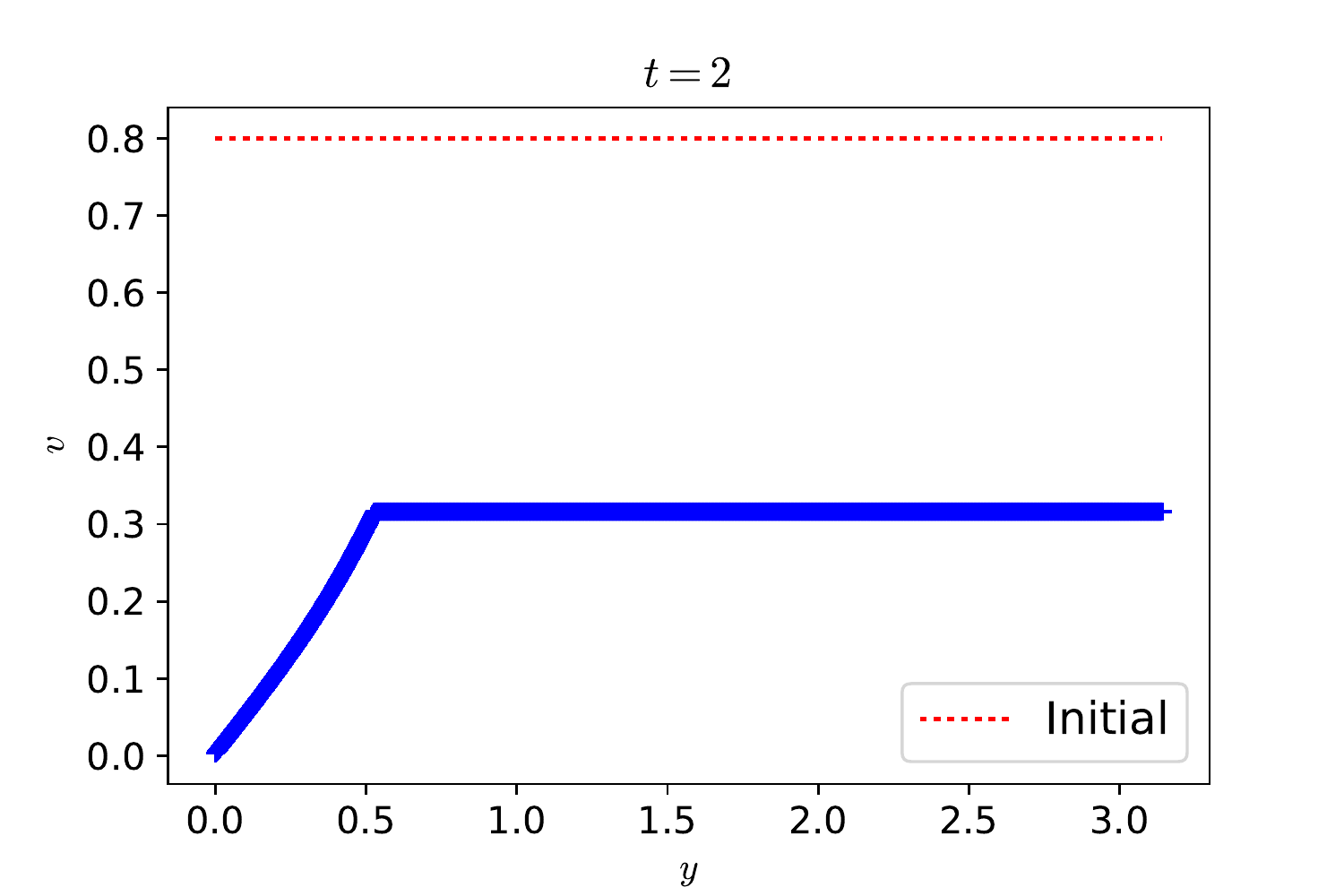}}
  \subcaptionbox{$\tau=5$}{\includegraphics[width=2.8in, trim={0 0 0.0 1.1cm},clip]{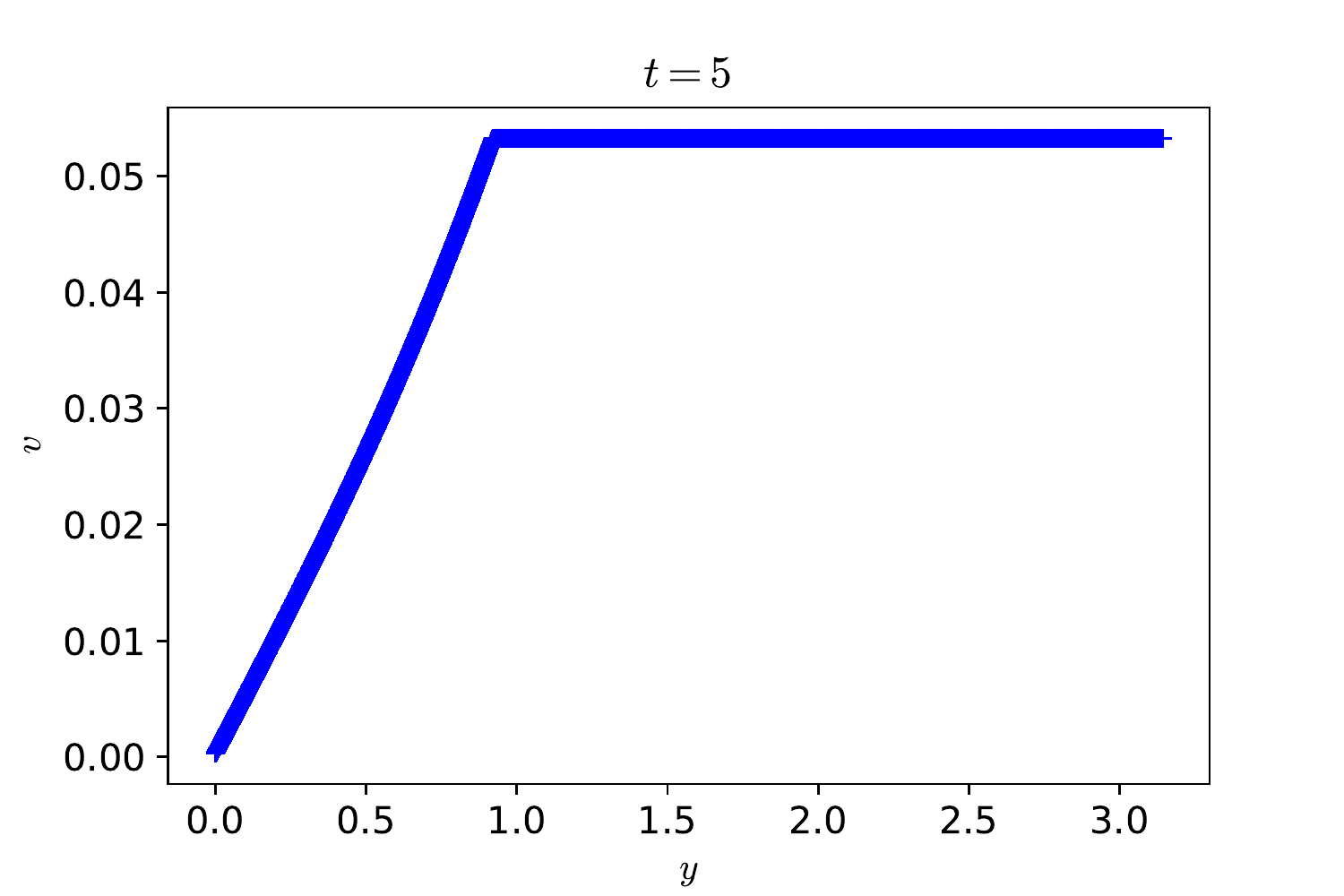}}
\centering
\caption{Numerical solution from a constant initial data.}
\label{FIG-Constant}
\end{figure}

We now consider the solutions to the cosmological Burgers model from a {\sl constant initial condition}, denoted by $v_0$, obtained with the second-order in space and fourth-order in time discretization.
We choose $\kappa = 2$ and the initial value $v_0 = 0.8$.  In the numerical tests,  the CFL number is taken to be $0.7$.  In the expanding case, $\tau_0 = 1$ is chosen. In Figure \ref{FIG-Constant1}, the solutions $y \mapsto v(\tau, y)$ if presented at the time $\tau = 5$
and 
for  $J = 100, 500, 1000, 5000$,  respectively.
Clearly, the results demonstrate that the approximate solution approches  our reference solution as $J$ increases. See Figure \ref{FIG-Constant}, where the  evolution of the reference solution is included as $\tau$ increases.


\subsection{Asymptotic behavior on an expanding background}

\begin{figure}[htbp]
\centering
  \subcaptionbox{$\tau=1.5$}{\includegraphics[width=2.8in, trim={0 0 0.0 1.1cm},clip]{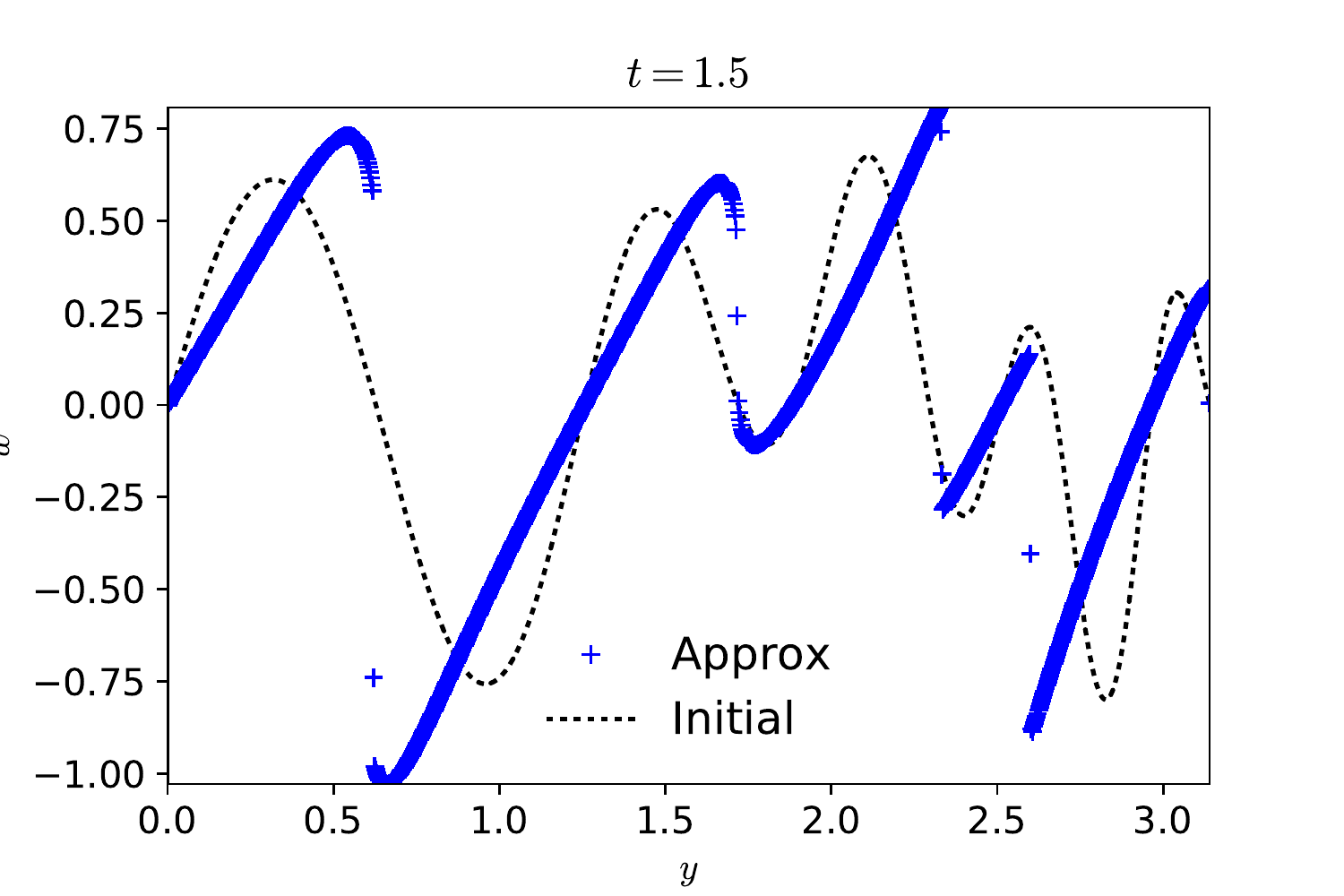}}
  \subcaptionbox{$\tau=2$}{\includegraphics[width=2.8in, trim={0 0 0.0 1.1cm},clip]{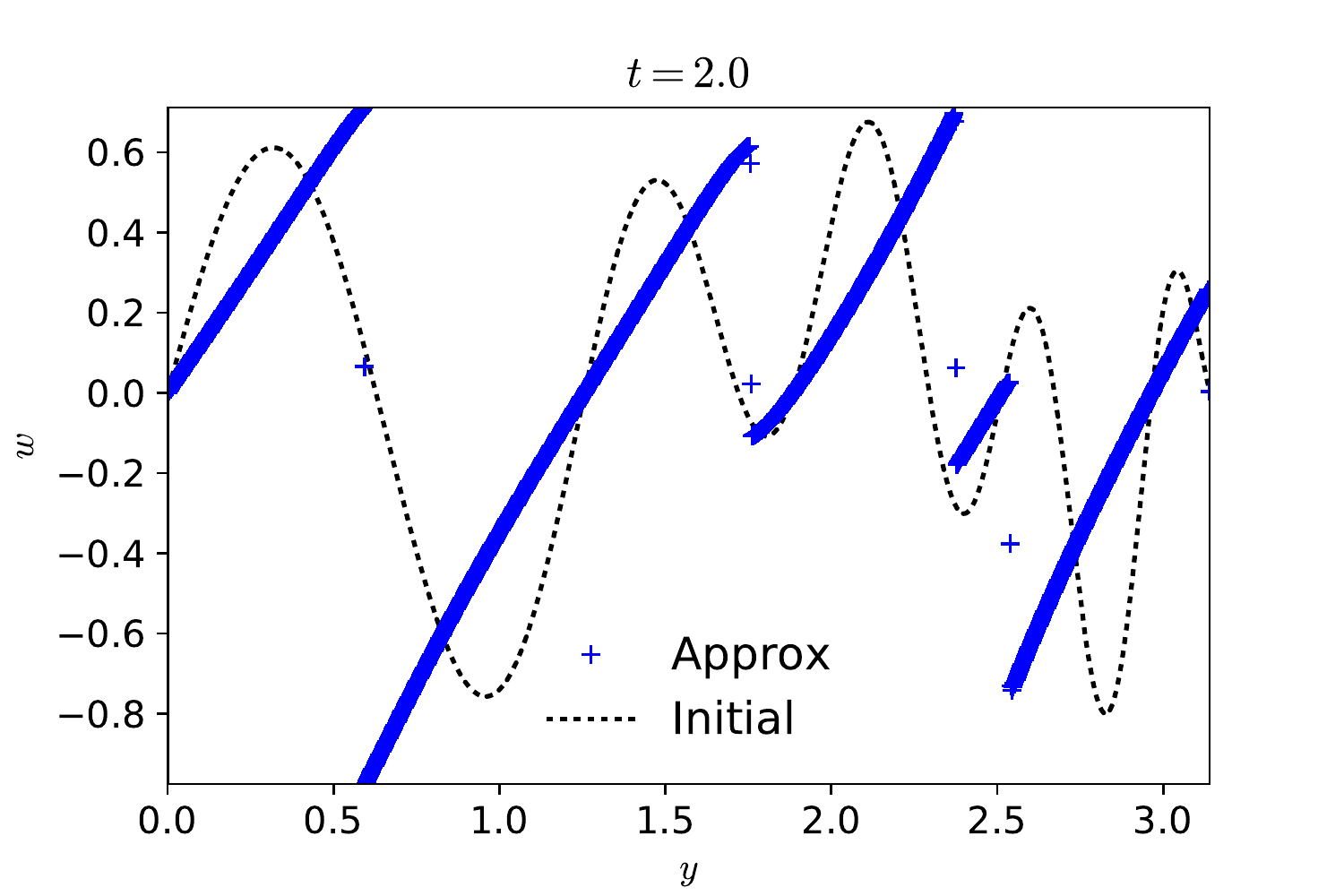}}
\centering
  \subcaptionbox{$\tau=16$}{\includegraphics[width=2.8in, trim={0 0 0.0 1.1cm},clip]{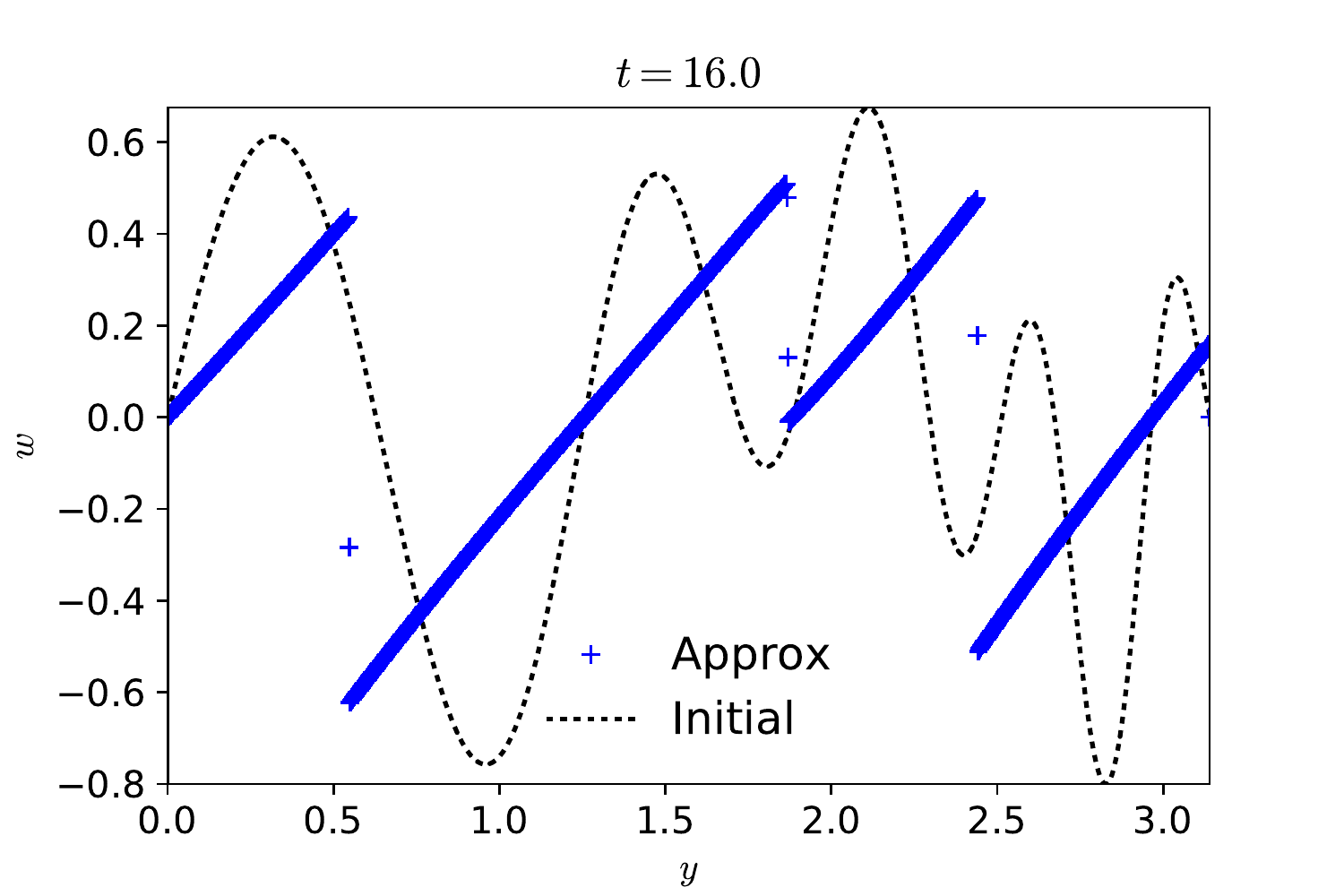}}
  \subcaptionbox{$\tau=64$}{\includegraphics[width=2.8in, trim={0 0 0.0 1.1cm},clip]{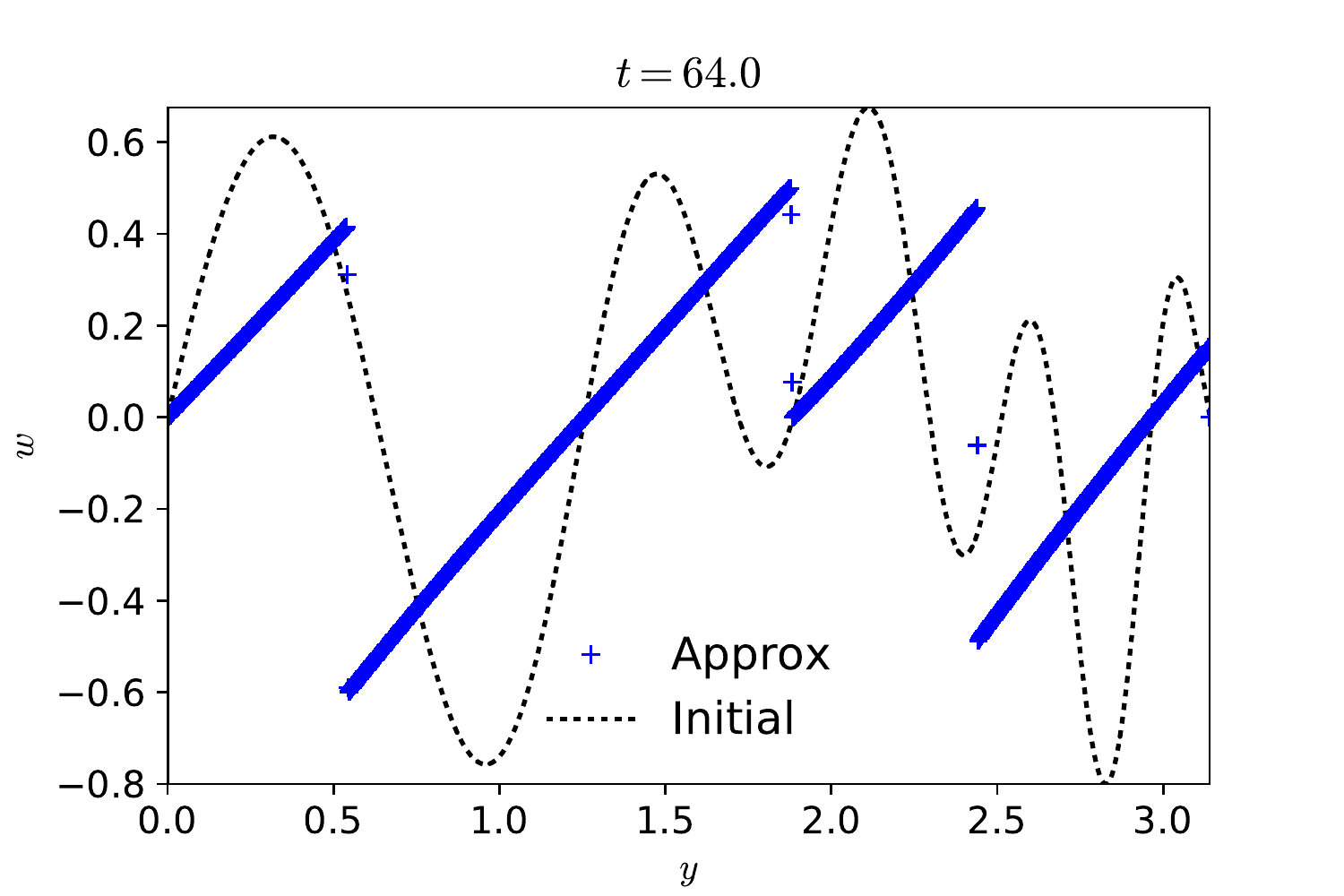}}
   
\caption{Rescaled solution for the model with $\kappa = 2$.}
\label{FIG-22}
\end{figure}

\begin{figure}[htbp]
\centering
  \subcaptionbox{$\tau=4$}{\includegraphics[width=2.8in, trim={0 0 0.0 1.1cm},clip]{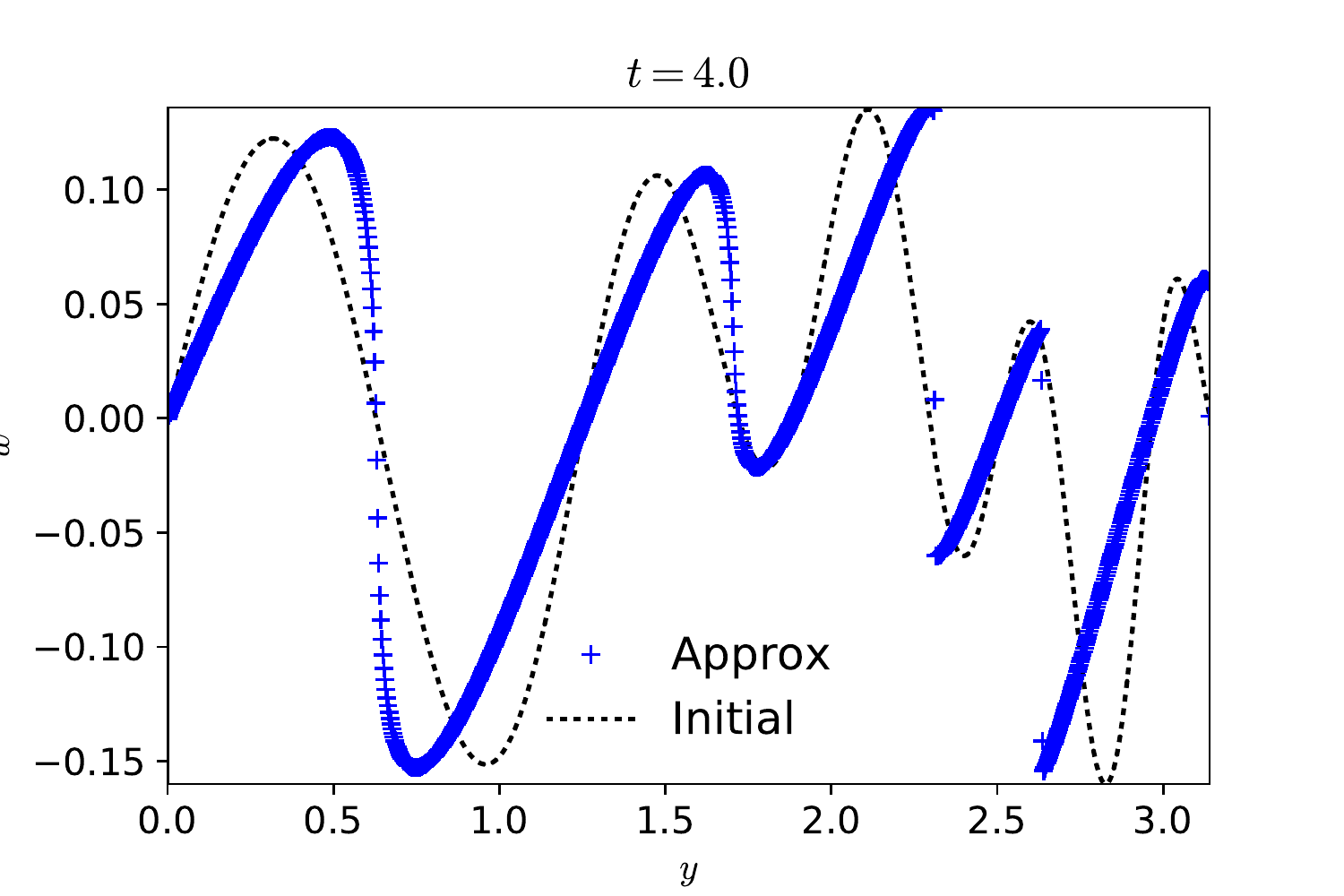}}
  \subcaptionbox{$\tau=64$}{\includegraphics[width=2.8in, trim={0 0 0.0 1.1cm},clip]{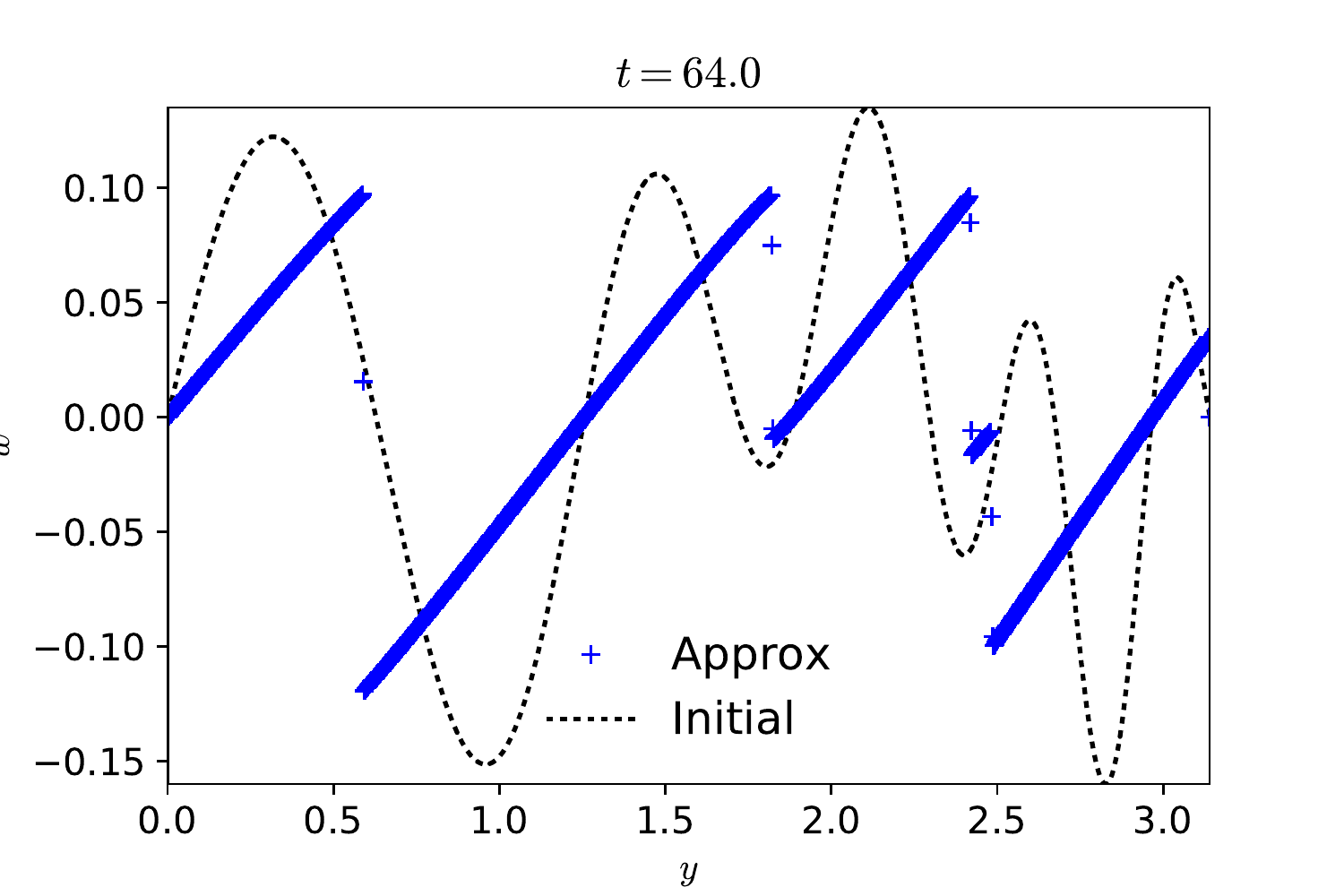}}
\centering
  \subcaptionbox{$\tau=256$}{\includegraphics[width=2.8in, trim={0 0 0.0 1.1cm},clip]{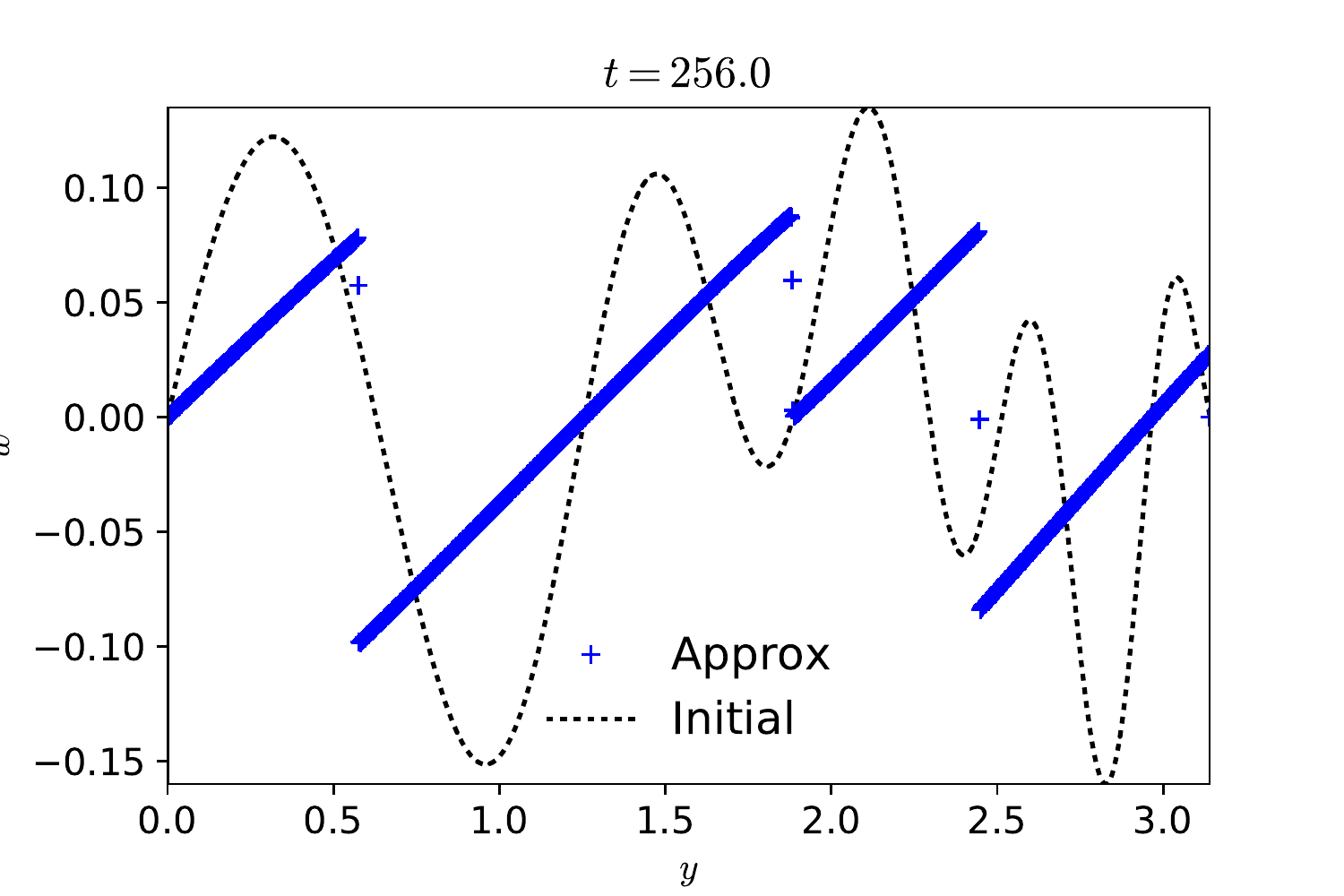}}
  \subcaptionbox{$\tau=512$}{\includegraphics[width=2.8in, trim={0 0 0.0 1.1cm},clip]{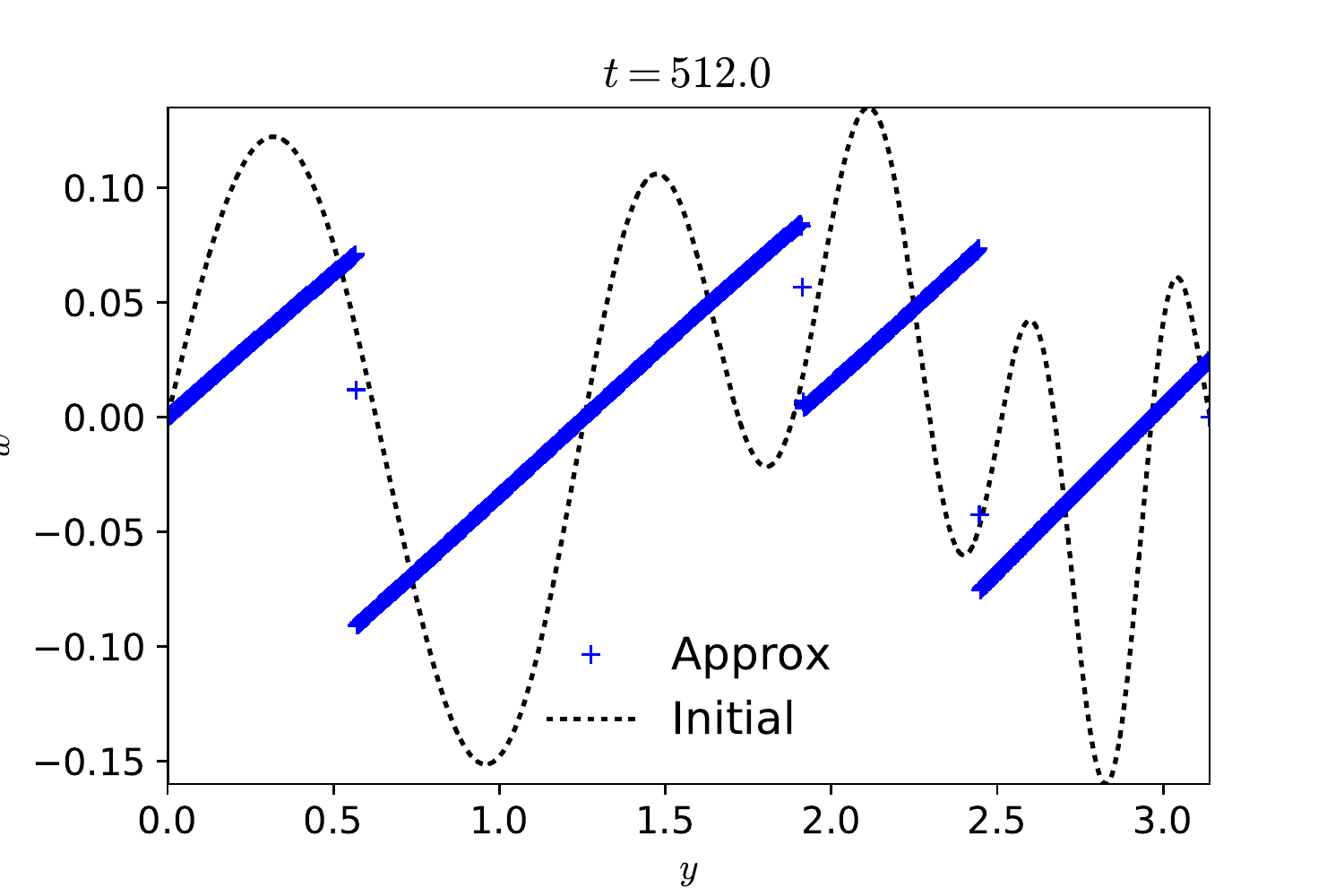}}
   
\caption{Rescaled solution for the model with $\kappa = 1$.}
\label{FIG-23}
\end{figure}

We now study the asymptotic behavior of the solutions using the proposed scheme at second-order in space and fourth-order in time.
The solution is expected to approach zero as time increases, and we propose to work with the following {\bf rescaled solution}
\bel{wy}
w(\tau, y) = \tau^\kappa v(\tau, y).
\ee
The asymptotic behavior of this function is thus computed in the expanding case when $\tau \to + \infty$. We take here $J = 1024$ and the CFL number $0.7$. At the initial time, $\tau_0 = 1$ the initial data is set to be 
$v_0(y) = 0.8 \sin(5y) \cos({\pi y^3-3 \over 7})$ when $􏰀\kappa =2$, 
and $v_0(y) = 0.16 \sin(5y) \cos({\pi y^3-3 \over 7})$ when $􏰀\kappa =1$. 
The evolution of the rescaled solution $w$ as $\tau$ increases is shown in Figures~\ref{FIG-22} and \ref{FIG-23}. We observe that the solutions, $y \mapsto w(\tau, u)$ eventually reaches a limit at a sufficiently large time $\tau > \tau_0$. Our numerical investigations lead us to state the following conclusion and conjecture. 

\begin{claim}[Cosmological Burgers flows on a future expanding background]
The asymptotic behavior of a solution to the cosmological Burgers model in the future expanding background is such that the solution $y \mapsto v = v(\tau, y)$ decays to zero uniformly in space:
$$
\lim_{\tau \to +\infty} v(\tau, y) = 0.
$$
Furthermore, the rescaled function $w = \tau^\kappa v$ approaches a (in general) non-trivial limit as $\tau \to +\infty$, which is a piecewise affine function with finitely many jumps.
\end{claim}


\subsection{Asymptotic behavior on  contracting background}

\begin{figure}[htbp]
\centering
  \subcaptionbox{$\tau=-0.5$}{\includegraphics[width=2.8in, trim={0 0 0.0 1.1cm},clip]{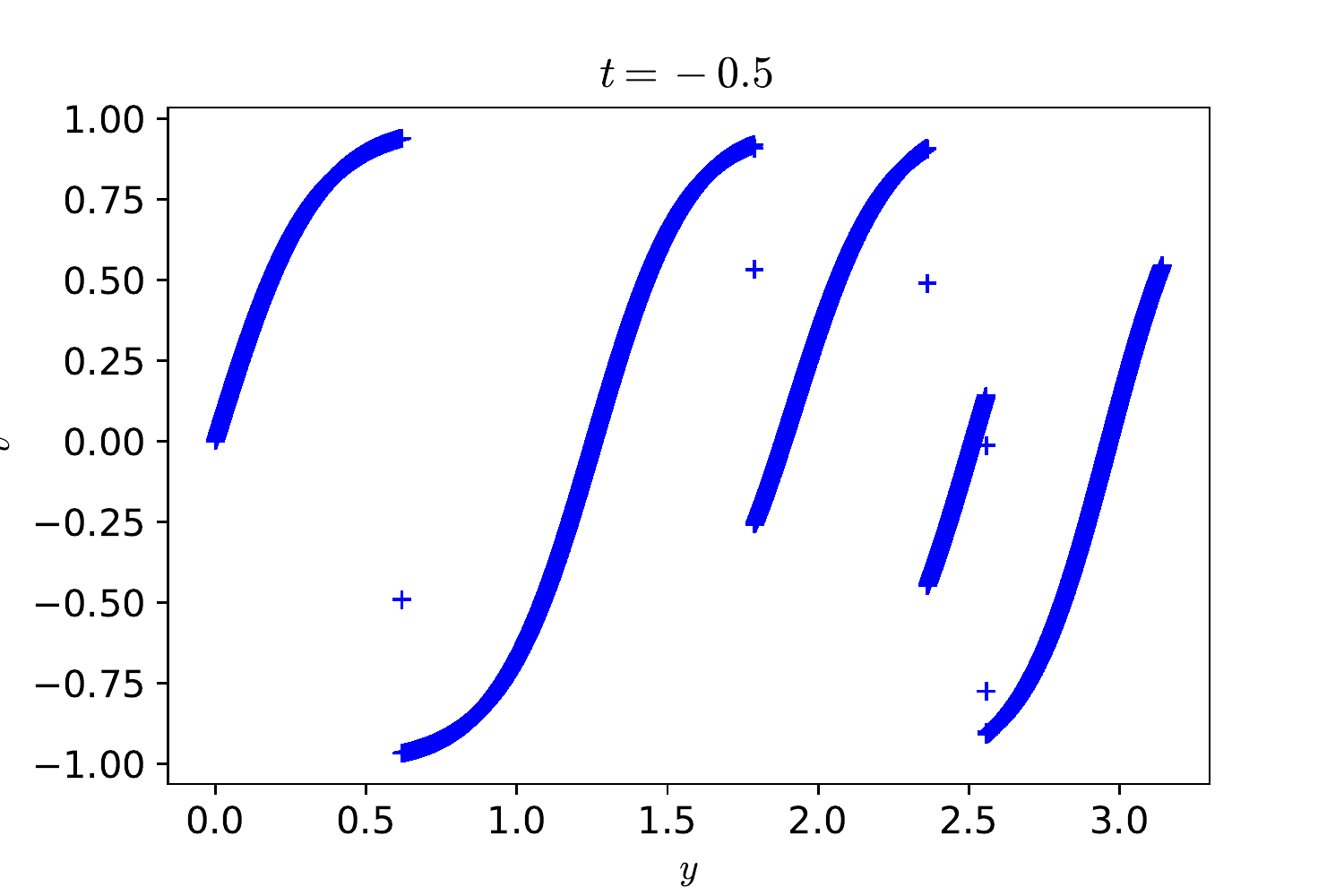}}
  \subcaptionbox{$\tau=-0.1024$}{\includegraphics[width=2.8in, trim={0 0 0.0 1.1cm},clip]{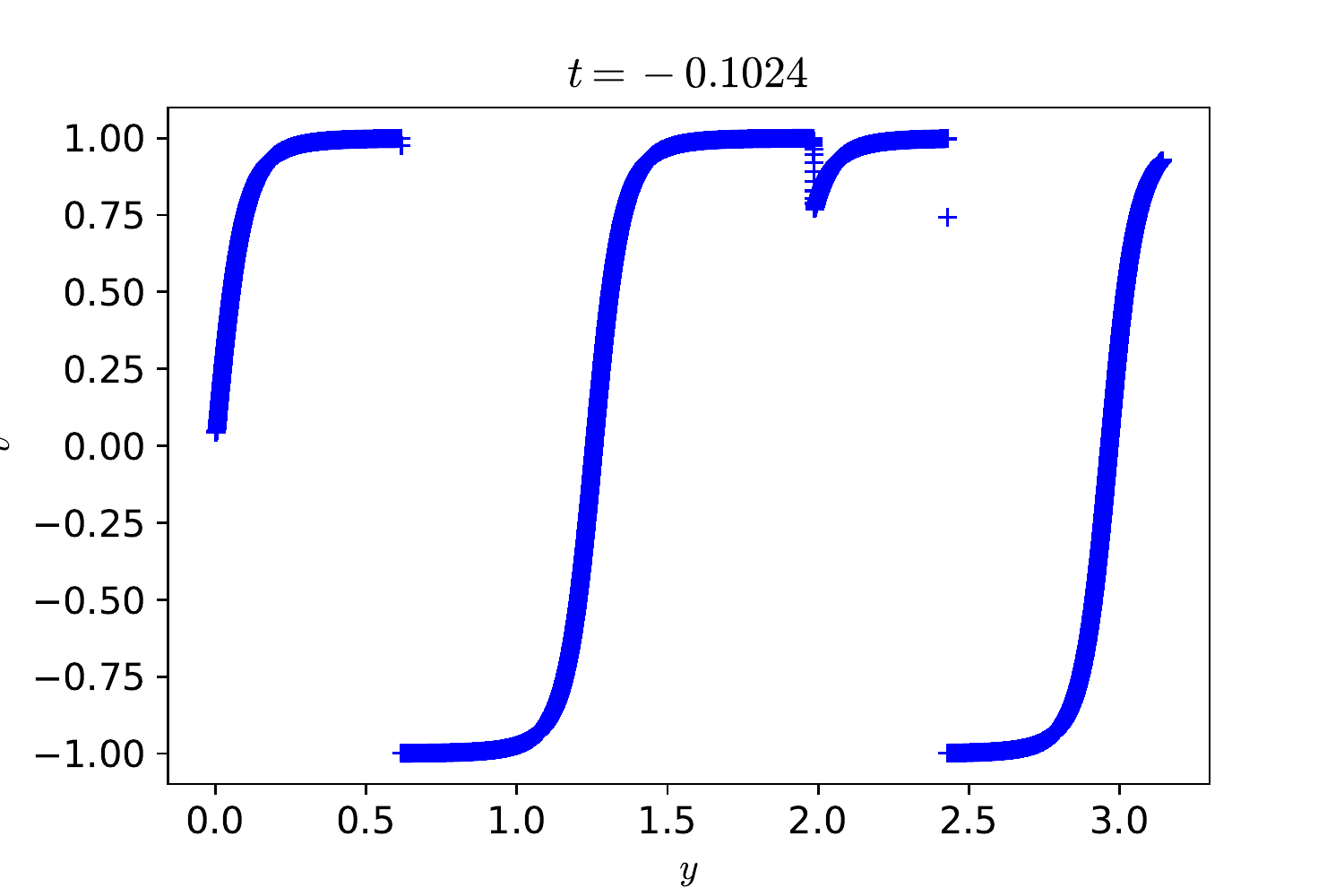}}
\centering
  \subcaptionbox{$\tau=-0.0016$}{\includegraphics[width=2.8in, trim={0 0 0.0 1.1cm},clip]{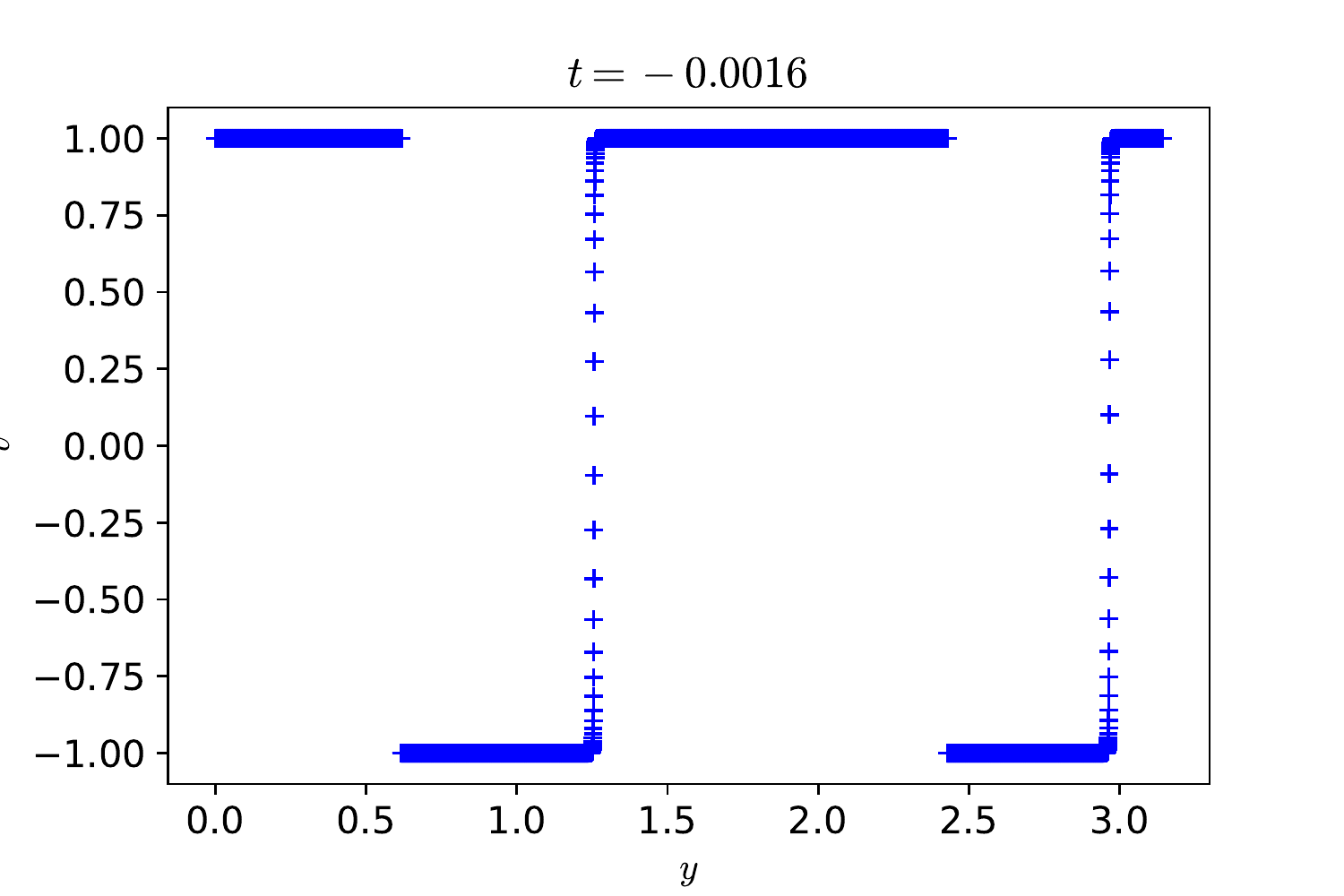}}
  \subcaptionbox{$\tau=-0.0001$}{\includegraphics[width=2.8in, trim={0 0 0.0 1.1cm},clip]{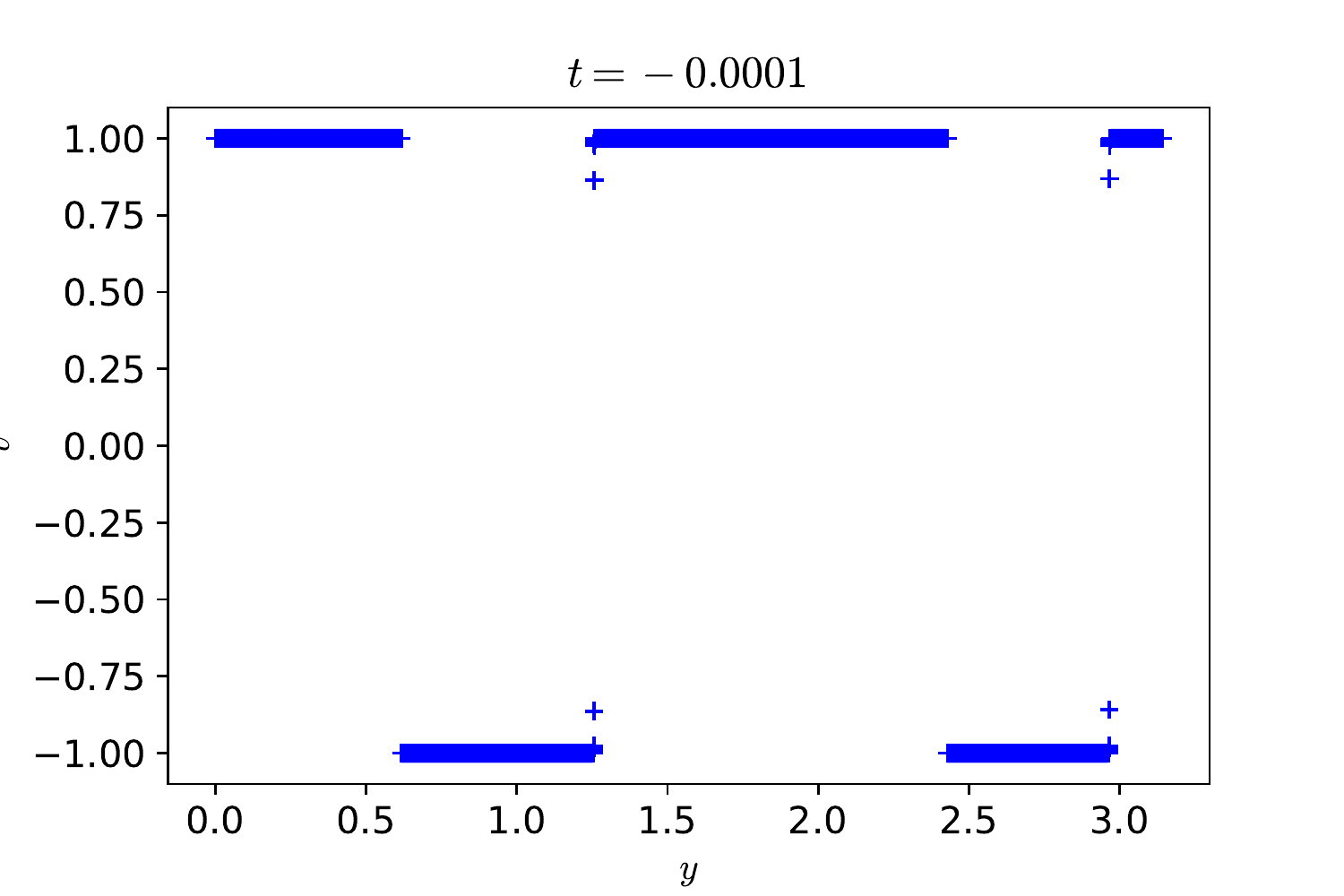}}
   
\caption{Asymptotic behavior on a contracting background.}
\label{FIG-con}
\end{figure}

\begin{figure}[htbp]
\centering
\epsfig{figure = 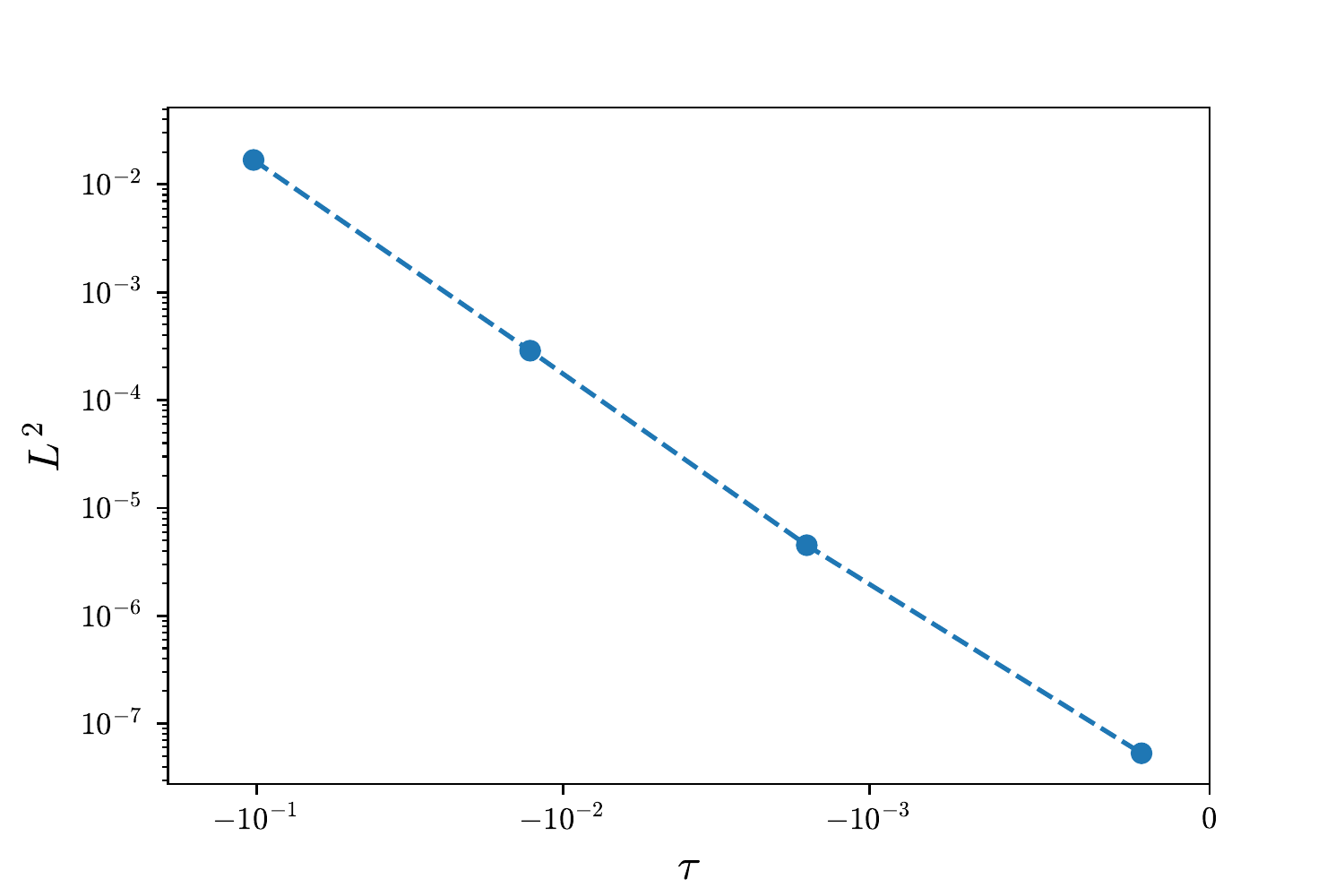, width = 2.9in}

\hspace{0.1in}

\caption{\label{FIG-convegence} Convergence of the solution on a contracting background.}
\end{figure}

The behavior of the solution in the future-contracting case as $\tau \to 0$ is next investigated. We take $J = 10000$ and the CFL number is $0.7$ and the initial value is prescribed at the initial time $\tau_0= -1$ to be $v_0(y) = 0.8 \sin(5y) \cos({\pi y^3-3 \over 7})$ and $􏰀\kappa =2$. 
The evolution of the solution as $\tau \to 0$ is shown in Figure \ref{FIG-con}. 
The $L^2$ norm of the numerical solutions is plotted in Figure \ref{FIG-convegence} at the times $\tau = -0.1024, -0.0128, -0.0016, -0.0002$,  respectively. We observe that the solutions converge to our reference solution as $\tau$ increases and the geometry is contracting.

\begin{figure}[htbp]
\centering
  \subcaptionbox{$\tau=-0.5$}{\includegraphics[width=2.8in, trim={0 0 0.0 1.1cm},clip]{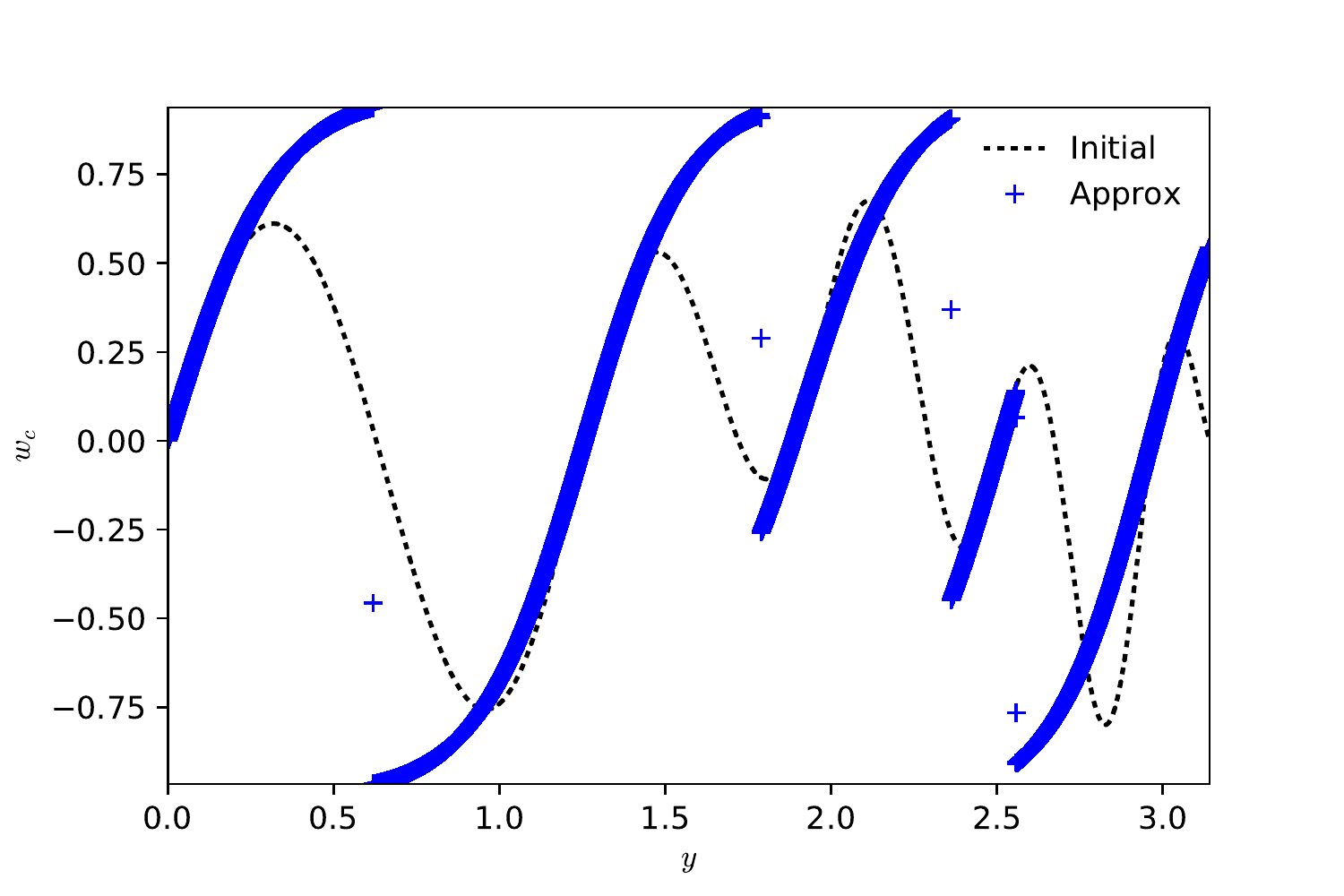}}
  \subcaptionbox{$\tau=-0.1$}{\includegraphics[width=2.8in, trim={0 0 0.0 1.1cm},clip]{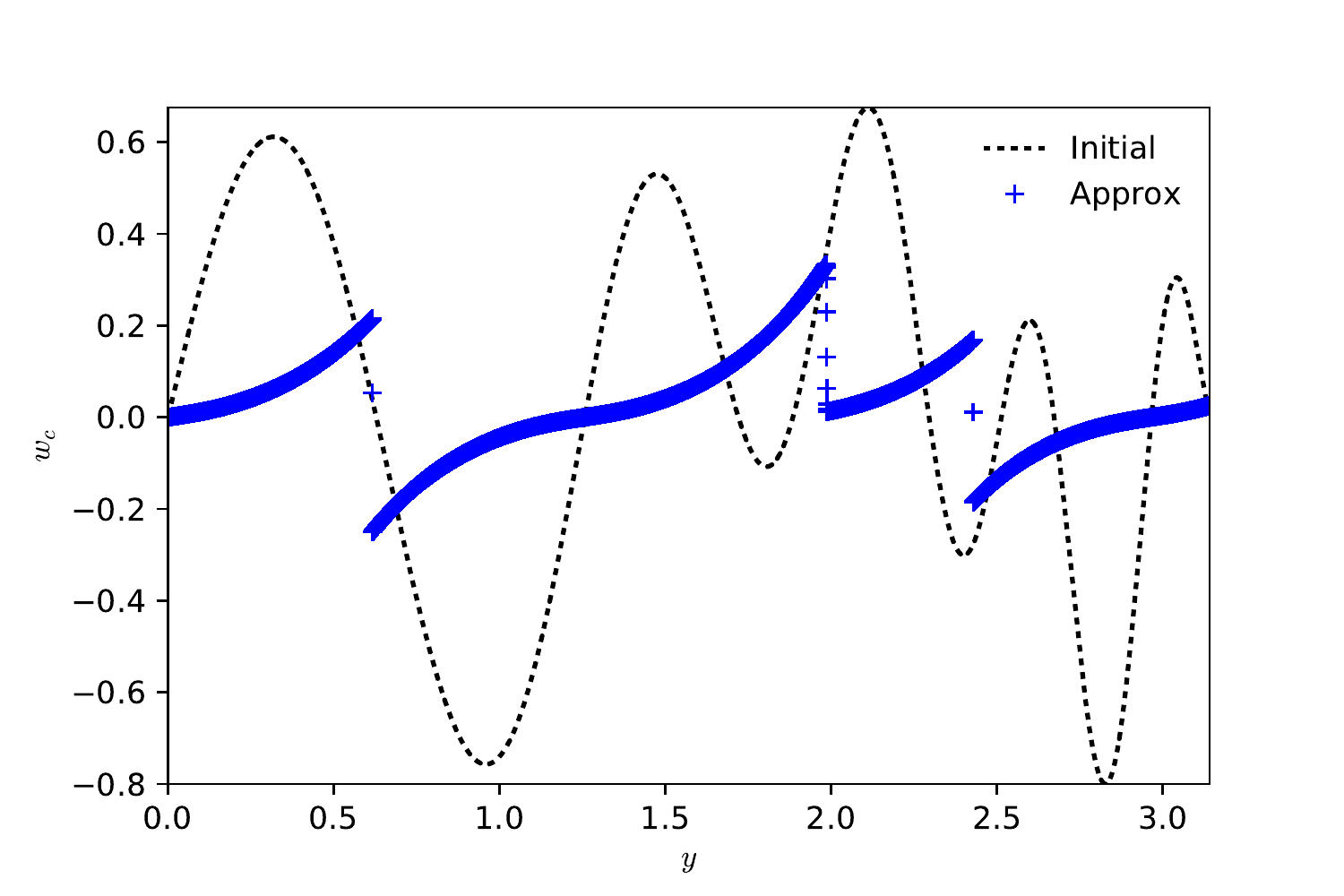}}
\centering
  \subcaptionbox{$\tau=-0.001$}{\includegraphics[width=2.8in, trim={0 0 0.0 1.1cm},clip]{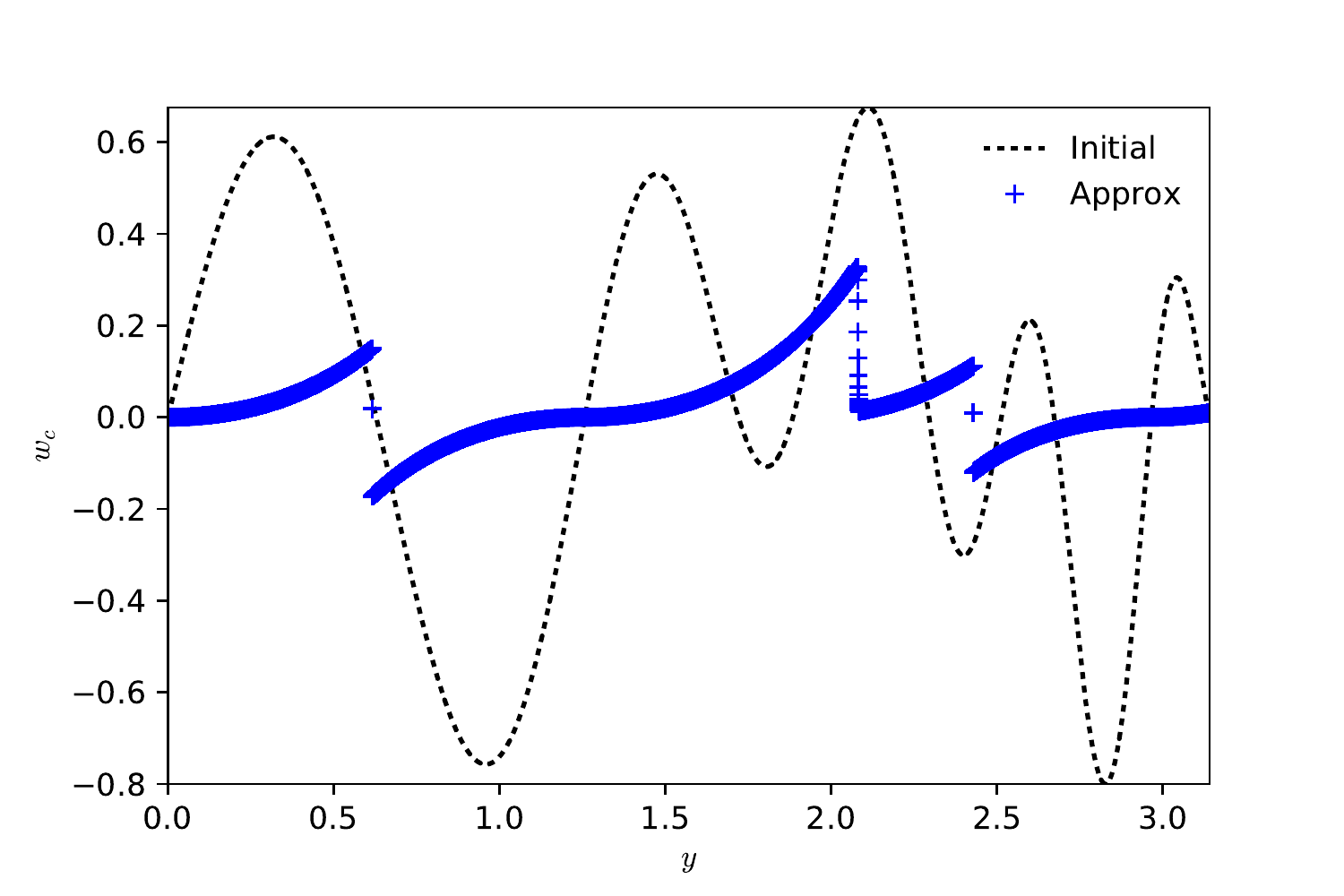}}
  \subcaptionbox{$\tau=-0.0001$}{\includegraphics[width=2.8in, trim={0 0 0.0 1.1cm},clip]{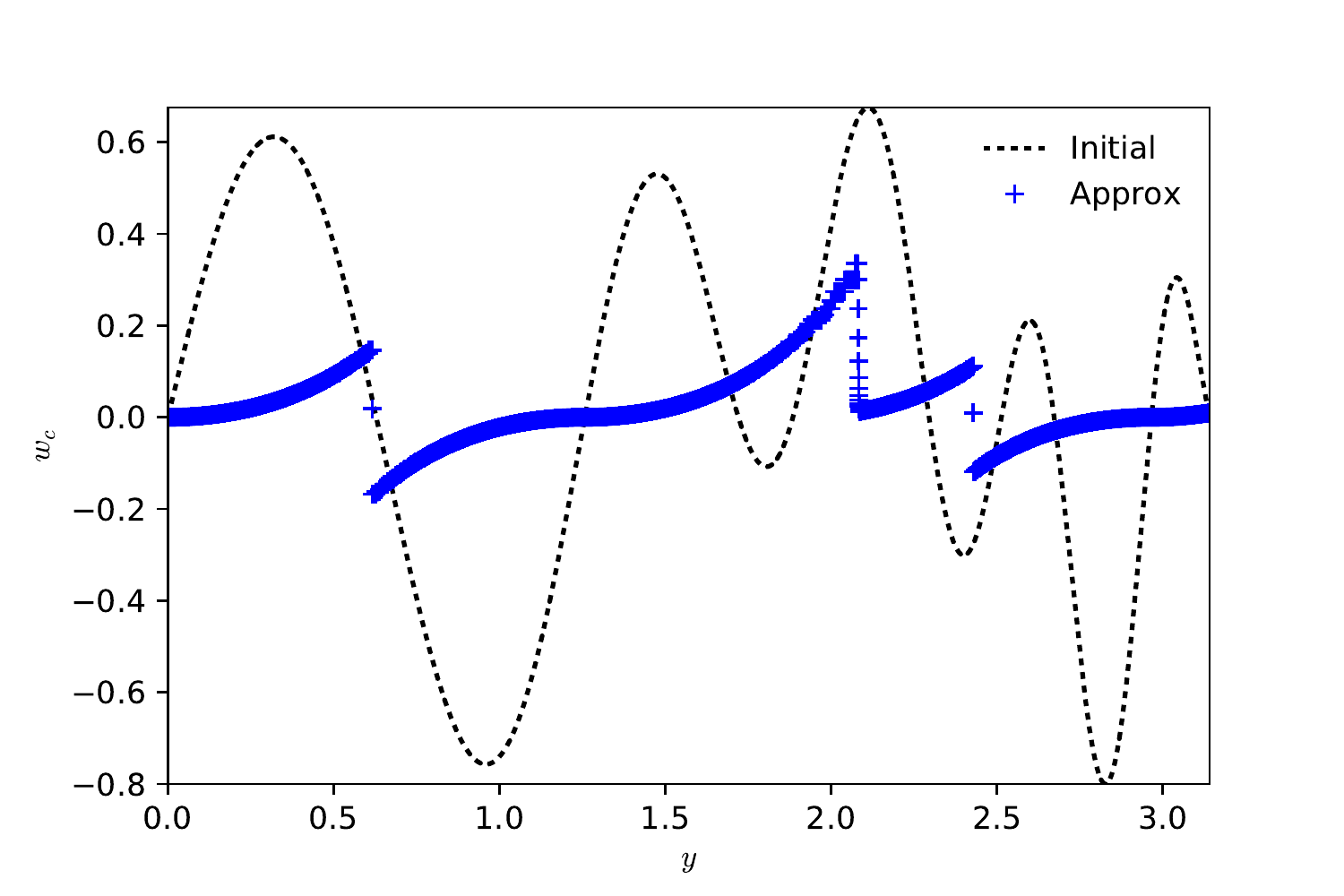}}
   
\caption{Rescaled solution on a contracting background with $\kappa =2$.}
\label{FIG-Wccon}
\end{figure}

It is convenient also to introduce the {\bf rescaled solution} $w$ defined, on a contracting background, by 
\bel{eq:wc}
w(\tau, y) = \sgn(v) {(-\tau)^{\kappa}\over \sqrt{1-v^2}}.
\ee
The evolution of this rescaled function $w(\tau, y)$ is presented in Figure~\ref{FIG-Wccon}, and our numerical investigations lead us to state the following conclusion and conjecture. 

\begin{claim}[Cosmological Burgers flows on a future-contracting background]
The asymptotic behavior of solutions to the cosmological Burgers model in the future-contracting case is such that the solutions $v = v(\tau, y)$ approach the light speed value $\pm 1$, that is, 
$$
\lim_{\tau \to 0} v(\tau, y) = \pm 1.
$$
Furthermore, the rescaled solution $w = \sgn(v) (-\tau)^{\kappa}/ \sqrt{1-v^2}$ approaches a non-trivial limit as $\tau \to 0$, which is a piecewise continuous function with finitely many jumps.
\end{claim}

\begin{figure}[htbp]
\centering
 \includegraphics[width=3in, trim={0 0 0.0 1.1cm},clip]{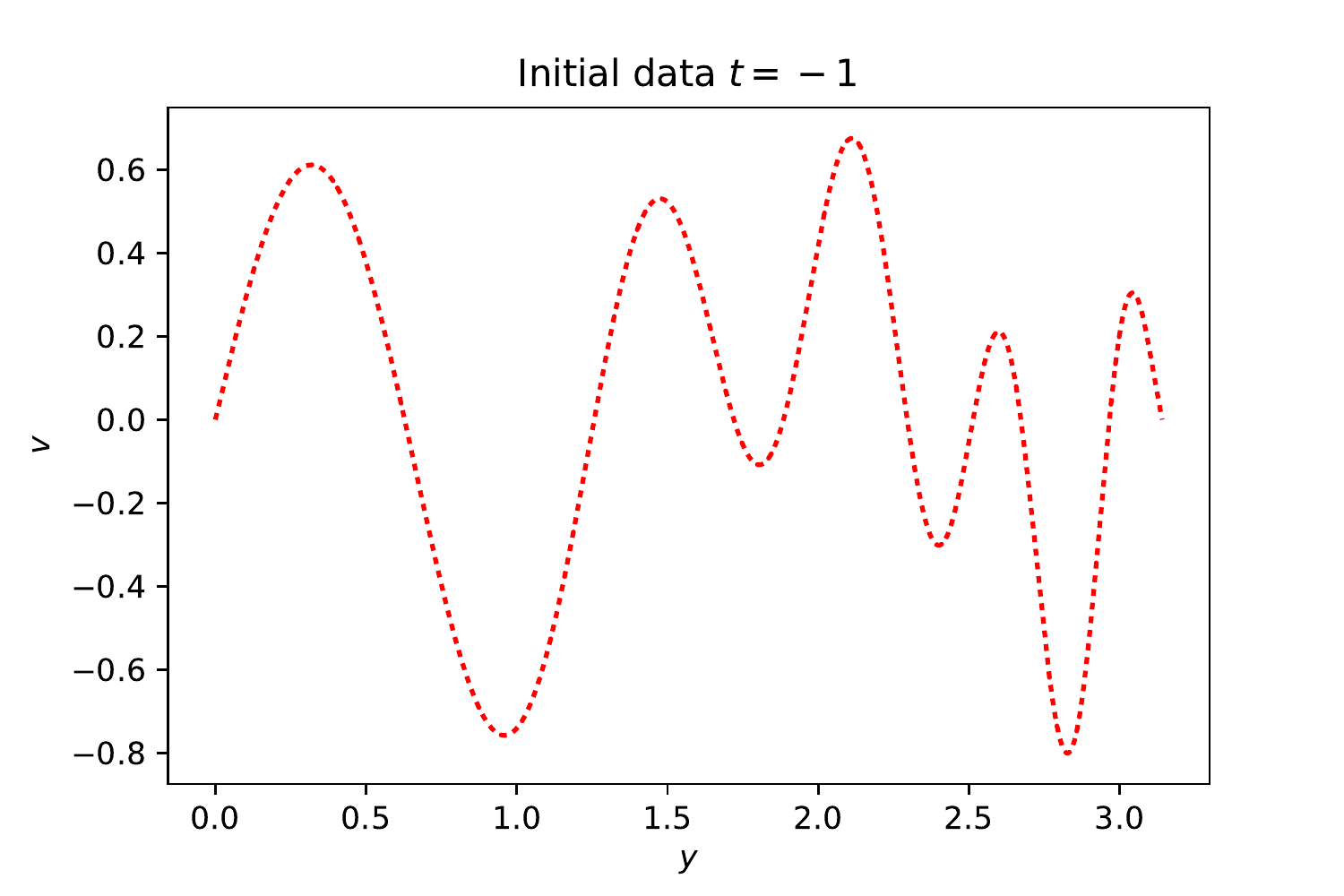}
\caption{Initial data at $\tau = -1$ for a contracting background.}
\label{FIG-in}
\end{figure} 


\begin{figure}[htbp]
\centering
 \subcaptionbox{}{\includegraphics[width=2.9in]{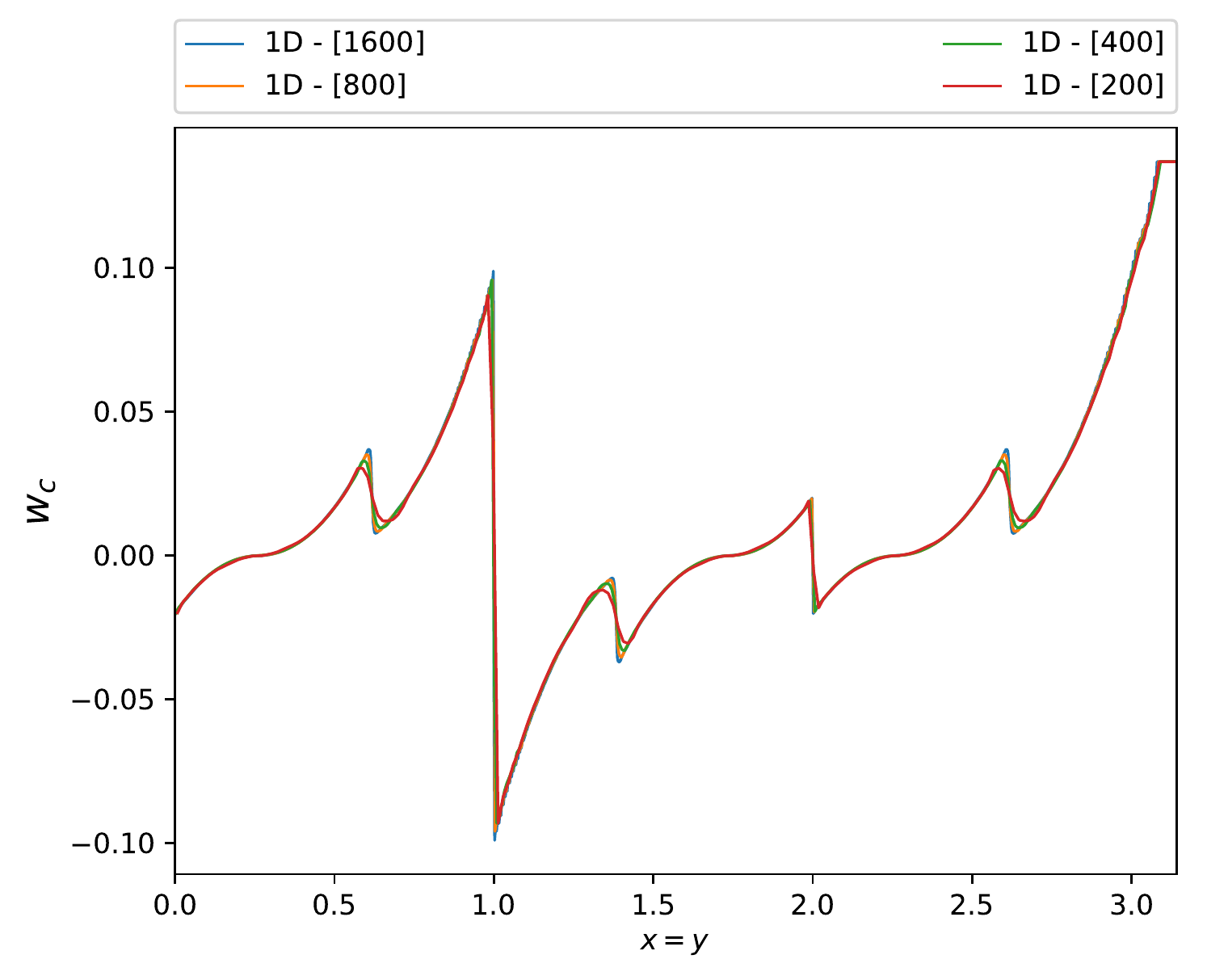}}
 \subcaptionbox{}{\includegraphics[width=2.9in]{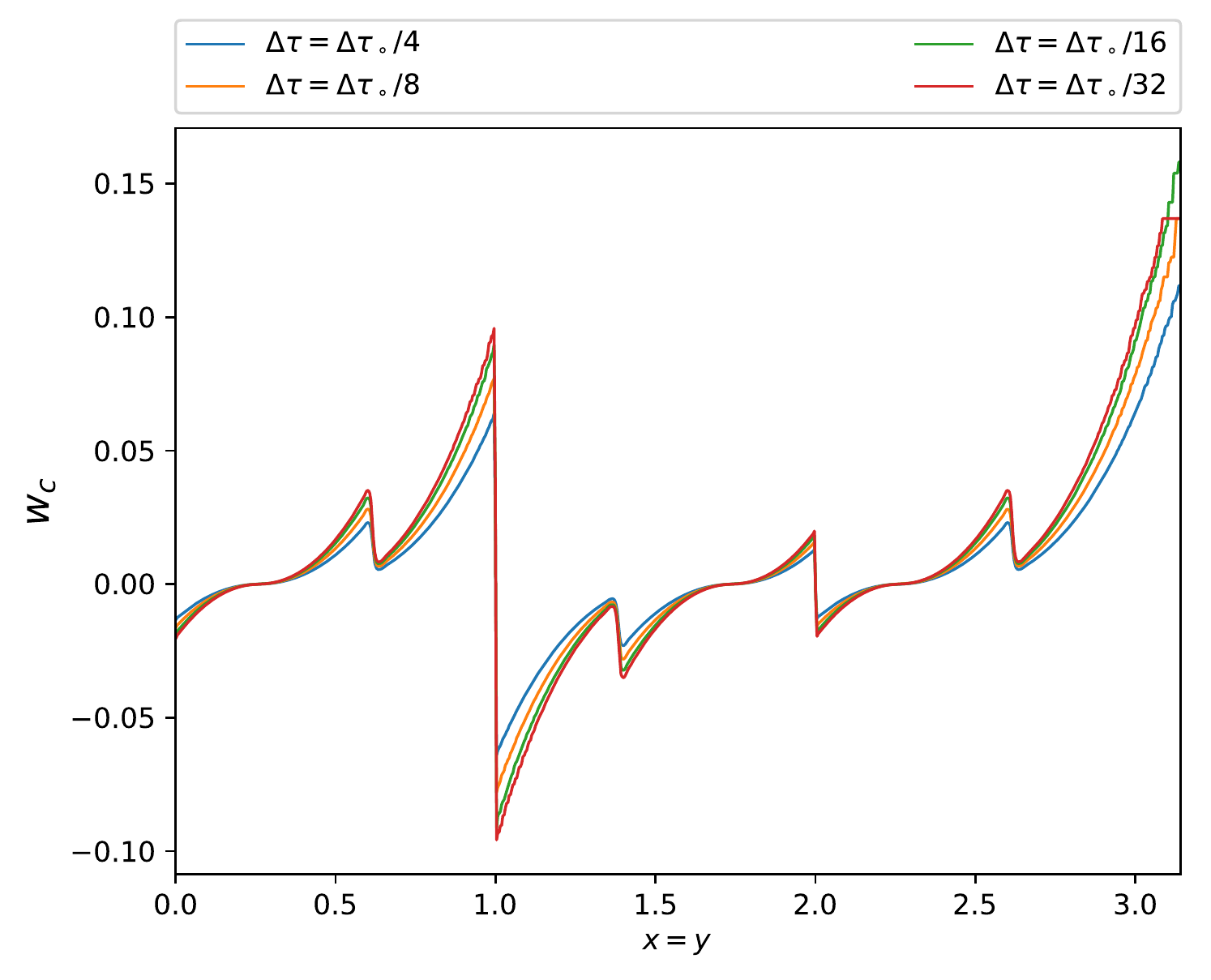}}

\caption{Comparing the rescaled solutions on a contracting background with $\kappa = 2$.
 (a) Different grid resolutions at $\tau = -10^{-4}$. 
(b) High-resolution grid for different $\Delta\tau$ with $\tau$ smaller than the one given by a CFL number $0.25$.}
\label{fig:comparisonCo2}
\end{figure}

\begin{figure}[htbp]
\centering
\includegraphics[width=2.8in]{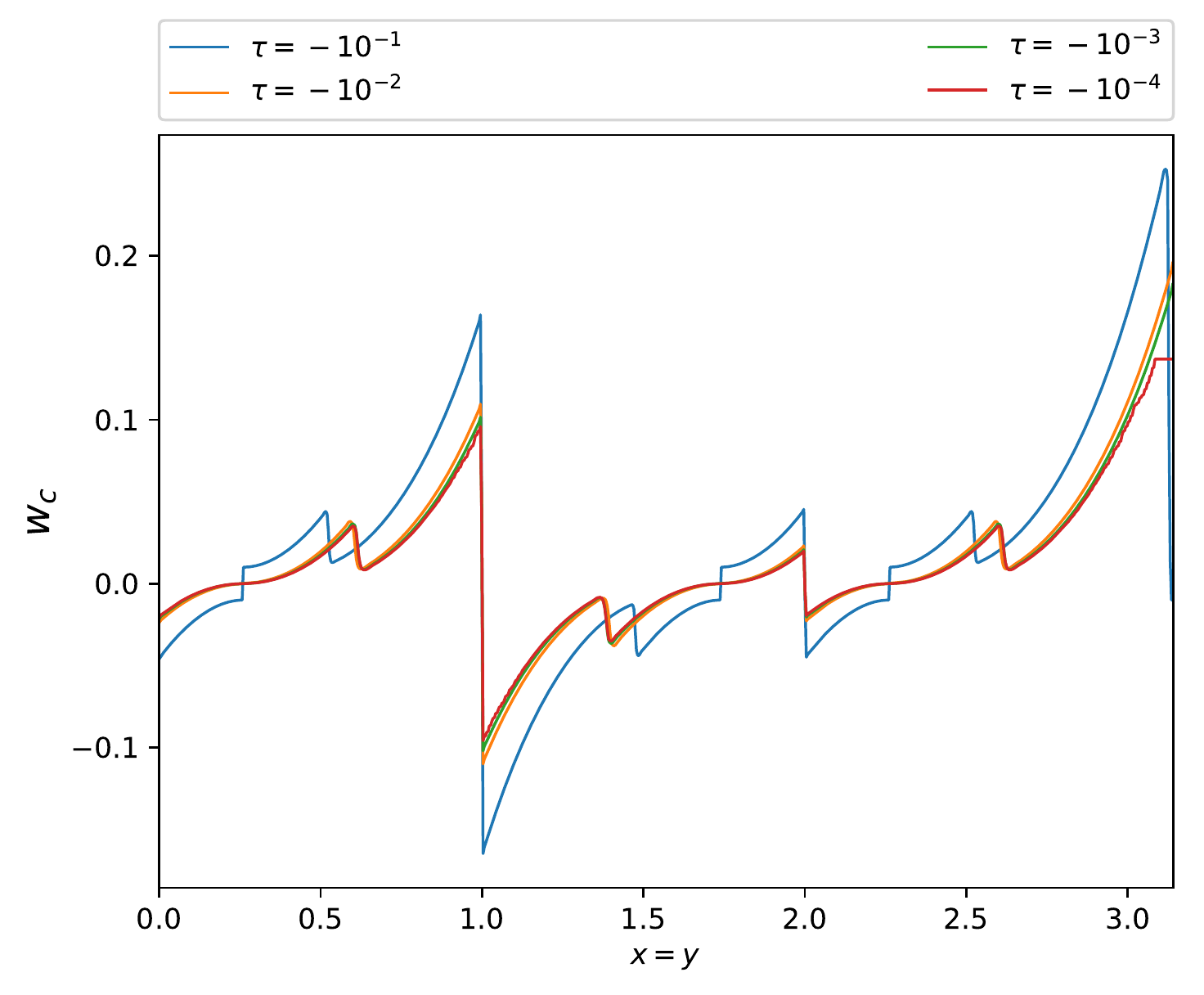}
 
\caption{Convergence of solutions at different $\tau$ with a grid of $800$ cells.}
\label{fig:5.17}
\end{figure}

\begin{figure}[ht!]
\centering
 \subcaptionbox{}{\includegraphics[width=2.8in]{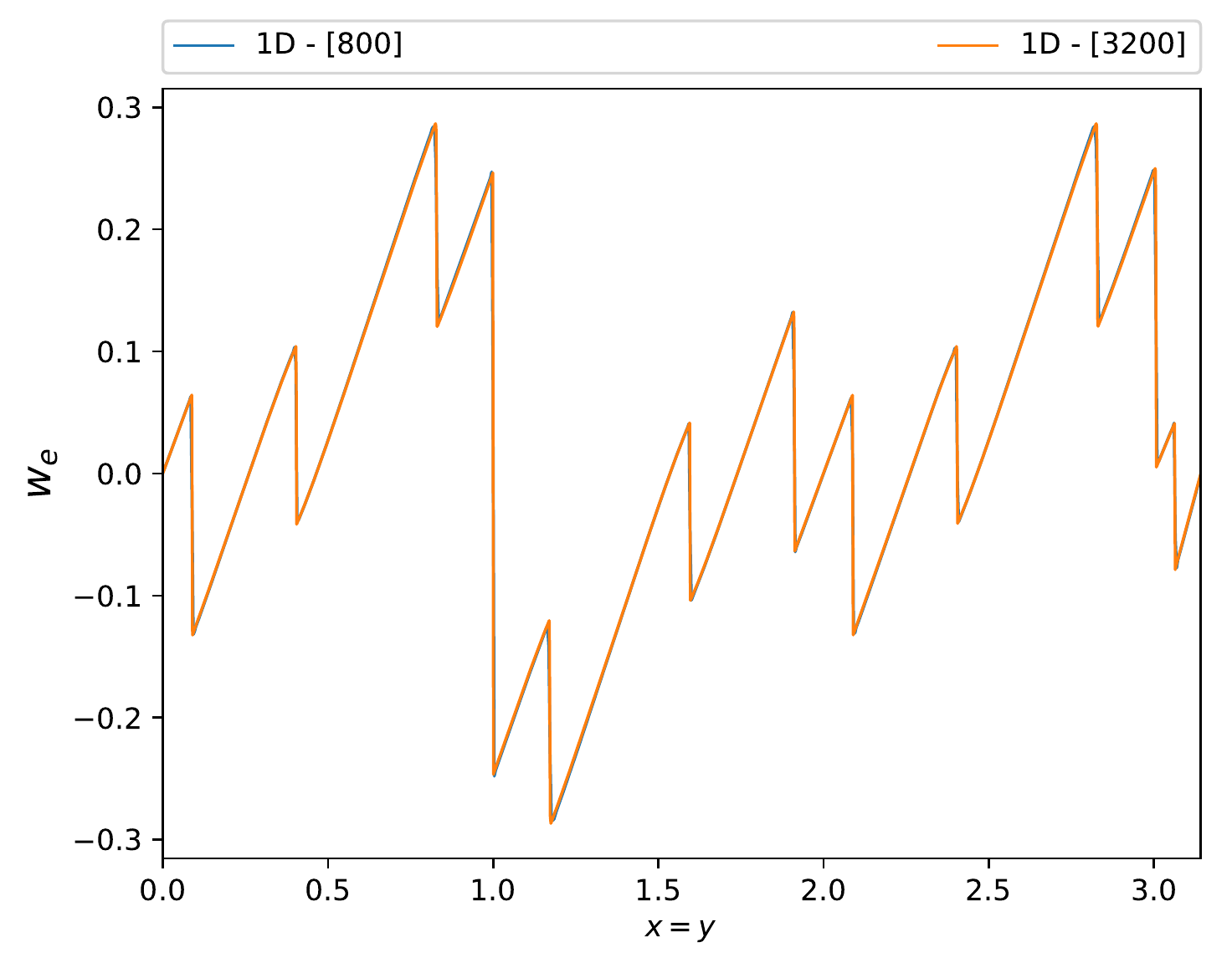}}
 \subcaptionbox{}{\includegraphics[width=2.8in]{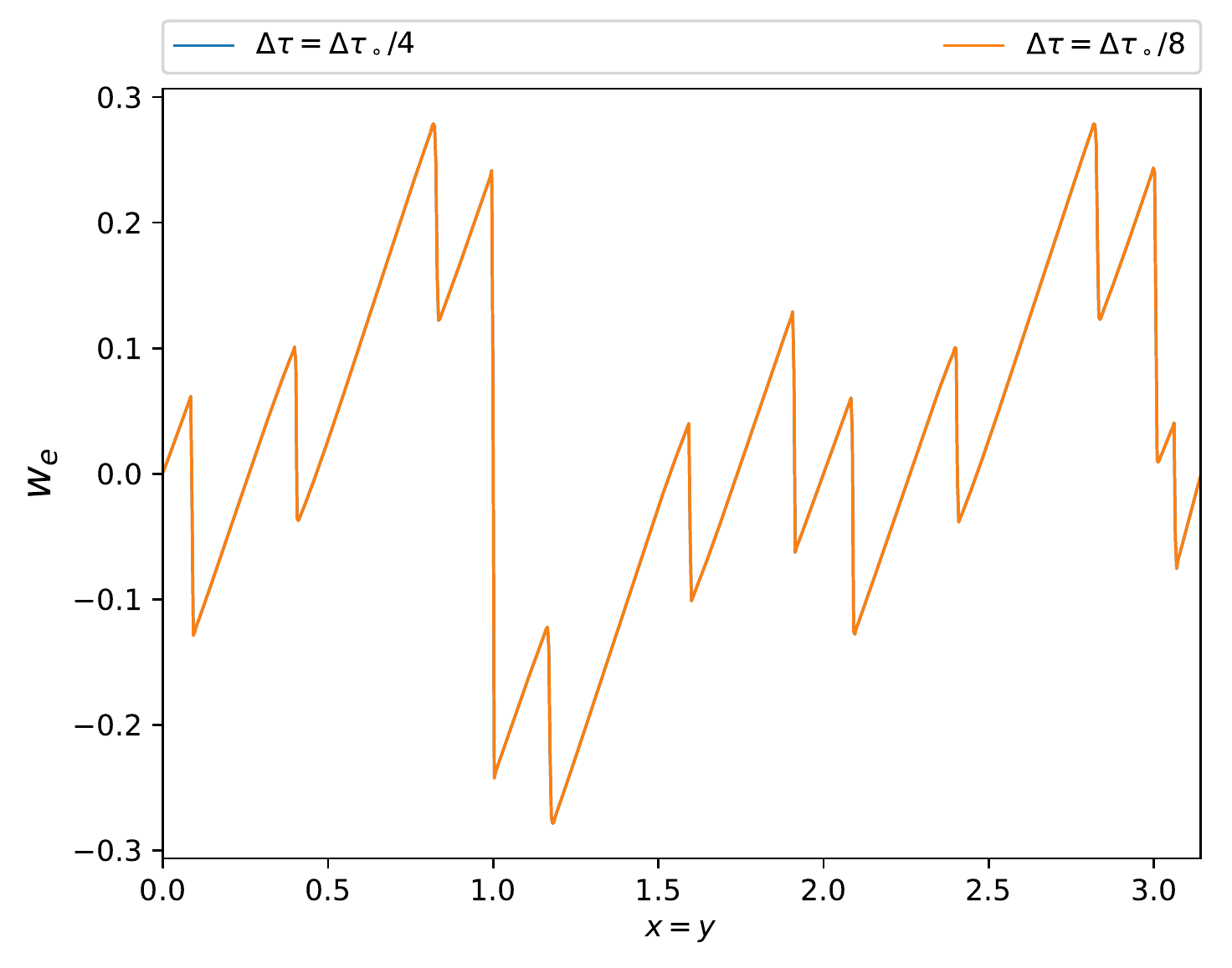}}

\caption{Comparing different rescaled solutions on an expanding background with $\kappa = 2$. 
(a) Different grid resolutions at $\tau = 1024$. 
(b) High resolution grid with $800$ cells and different $\Delta\tau$.}
\label{fig:comparisonEx2}
\end{figure}

\begin{figure}[htbp]
\centering
\includegraphics[width=2.8in]{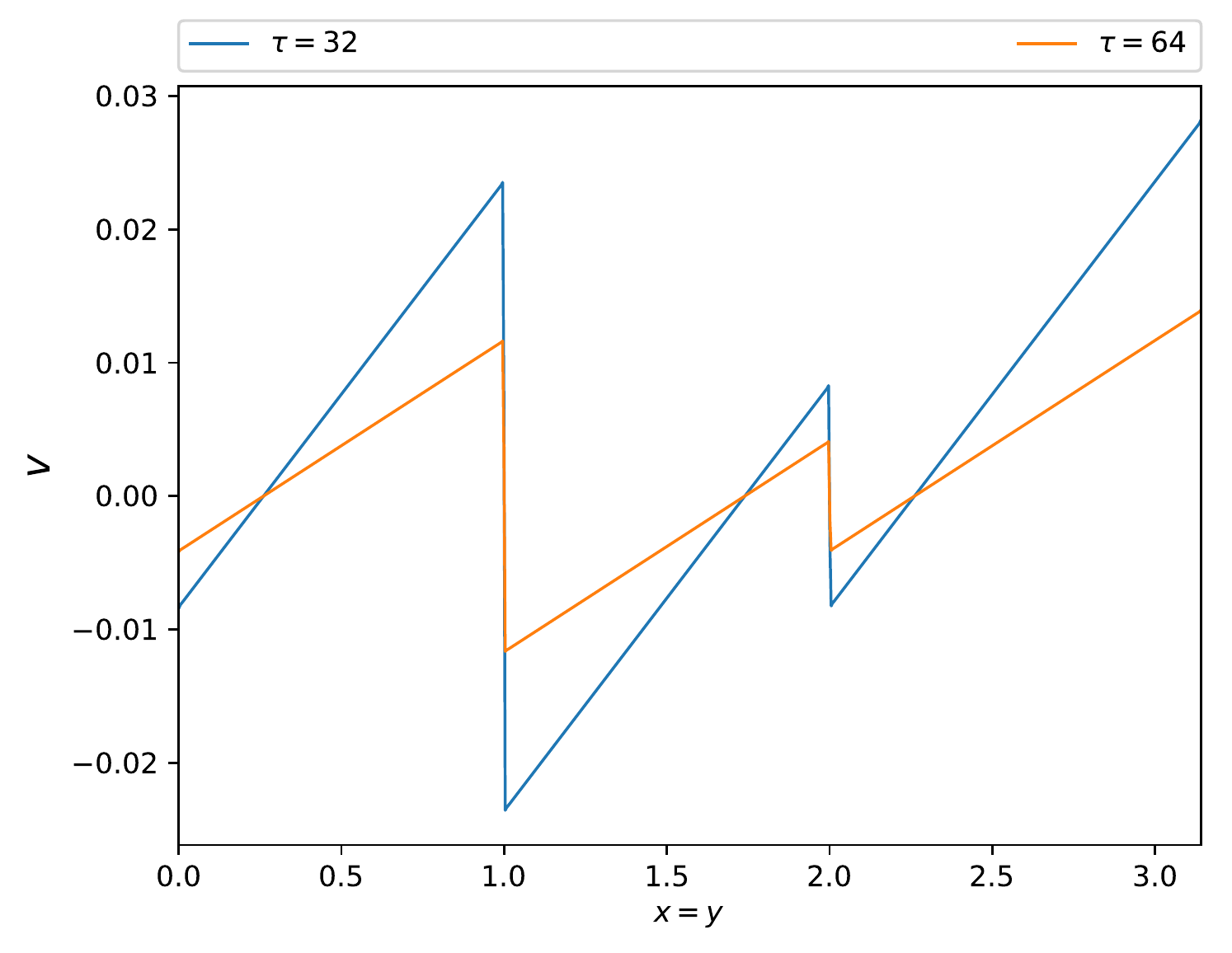}
 
\caption{Standard Burgers equation at  the times $\tau = 32$ and $64$.}
\label{fig:classicburgers}
\end{figure}

Further tests in $(1+1)$ dimensions are presented in Figures~\ref{fig:comparisonCo2} and \ref{fig:5.17}. In Figure~\ref{fig:comparisonCo2} we have plotted (in (a)) four different  rescaled solutions at $\tau = -10^{-4}$ with different grid sizes. The solutions follow a similar evolution and establish the convergence of our scheme. In 
Figure~\ref{fig:comparisonCo2} we have plotted in 
(b) four different solutions with different $\Delta\tau$, CFL, and $\Delta\tau_0$ where 
the CFL number is  $0.25$, a grid containing 800 cells where $\Delta\tau_0$ is 
\begin{align}
&\text {When } \kappa > 1,  \quad \Delta \tau_0 = \frac{\vert \tau_n \vert}{\kappa}.
\label{eq:Riemann 12}          
\end{align}
Figure~\ref{fig:5.17} also compares different solutions with the same grid and different $\tau$.
Similarly to the contracting background, another test is presented in Figure~\ref{fig:comparisonEx2} in $(1+1)$ dimensions where $\kappa = 2$. For this test $\Delta\tau_0$ is given by 
\be
\Delta\tau_0 =  0.9 \min_{j} \left( \frac{\Delta y}{|v_j|}, \frac{\tau}{\kappa\left(1 - v_{j}^{2} \right)}\right),
\label{eq:StabilityEx}  
\ee
where the time indices are omitted.
The standard Burgers equation (i.e.~with $a(t)=1$) is solved with the same initial conditions \eqref{eq:IC2D} in $(1+1)$ dimensions. The solutions in Figure~\ref{fig:classicburgers} show the the number of shocks are significantly less than the number of shocks in the solutions of the expanding case.

\section{Global dynamics of $(2+1)$--cosmological Burgers flows}
\label{FVforBurgers}

\subsection{The algorithm in $(2+1)$-dimensions}

\paragraph*{First-order finite volume discretization.}

We now turn our attention to the model in $(2+1)$-dimensions. We use our finite volume Godunov-type scheme with second-order accuracy in space, fourth-order (expanding background) or a third-order (contracting background) accuracy in time, 
and we solve the cosmological Burgers model \eqref{RBEat1}, that is, written as 
\bel{eq:2+1Burger's - Flux-Source}
v_\tau + f(v)_x + g(v)_y = m(\tau) \, h(v).
\ee
with flux-functions $f(v)= g(v) = v^2/2$ and source given by 
$m(\tau) = \kappa/\tau$ and $h(v) = - v(1-v^2)$. The scheme is based on a uniform grid of intervals $[x_{j-1/2},x_{j+1/2}]$ and $[y_{k-1/2},y_{k+1/2}]$. Here, $j$ and $k$ are integers describing the $x$ and $y$ directions, respectively, and we use the same notation $\Delta x, \Delta y$ as in 1D.
 The cell averages of the main variable and the source are  
\be
\barv_{j,k}(\tau) \approx \frac{1}{\Delta{x}\Delta{y}}\int\limits_{x_{j-1/2}}^{x_{j+1/2}}\int\limits_{y_{k-1/2}}^{y_{k+1/2}} v(\tau, x, y) dxdy,
\quad
\overline{\boldsymbol{S}}_{j,k}(\tau )\approx \frac{1}{\Delta{x}\Delta{y}}\int\limits_{x_{j-1/2}}^{x_{j+1/2}}\int\limits_{y_{k-1/2}}^{y_{k+1/2}} m(\tau)\,h(v) \, dxdy.         
\ee
The semi-discrete version of the first-order Godunov-type scheme then reads 
\bel{eq:Burger's discretization}
\dfrac{d}{dt}\barv_{j,k}(\tau) 
= -\frac{\boldsymbol{H}_{j+1/2,k} (\tau) - \boldsymbol{H}_{j-1/2,k} (\tau)}{\Delta{x}} -\dfrac{\boldsymbol{H}_{j,k+1/2} (\tau) - \boldsymbol{H}_{j,k-1/2} (\tau)}{\Delta{y}} + \overline{\boldsymbol{S}}_{j,k}(\tau).       
\ee

\bse

Introducing a time step on the interval $[\tau_n, \tau_{n+1}]$, we arrive at the fully-discrete first-order finite volume scheme:
\bel{eq:Burger's discretization2}  
v_{j,k}^{n+1} = v_{j,k}^{n} - \frac{\Delta{\tau}}{\Delta{x}} \left(\boldsymbol{H}_{j+1/2,k}^n - \boldsymbol{H}_{j-1/2,k}^n \right) -\frac{\Delta{\tau}}{\Delta{y}} \left(\boldsymbol{H}_{j,k+1/2}^n - \boldsymbol{H}_{j,k-1/2}^n \right) + \Delta{\tau}\overline{\boldsymbol{S}}^n_{j,k}. 
\ee
It remains to specify the numerical discretization of the flux and source, and we set 
\bel{eq:source at n}
\overline{\boldsymbol{S}}^n_{j,k} = m(\tau_n) h(v^n_{j,k}),
\ee
\bel{eq:Godunov flux}
\begin{split}
&\boldsymbol{H}_{j-1/2,k} = f_G(v_{j-1/2,k}^-, v_{j-1/2,k}^+), \qquad \boldsymbol{H}_{j+1/2,k} = f_G(v_{j+1/2,k}^-, v_{j+1/2,k}^+), \\
&\boldsymbol{H}_{j,k-1/2}^n = f_G(v_{j,k-1/2}^-, v_{j,k-1/2}^+), \qquad \boldsymbol{H}_{j,k+1/2}^n = f_G(v_{j,k+1/2}^-, v_{j,k+1/2}^+). 
\end{split}
\ee
Here, $f_G(v^n_l, v^n_r)$ and $f_G(v^n_b, v^n_t)$ in the $x$-direction are nothing but the standard Godunov fluxes and are obtained by solving a local Riemann problem, as explained earlier. The term $f_G(v^n_b, v^n_t)$ is defined similarly in the $y$-direction.

\ese


\paragraph*{Second-order finite volume discretization.}

\bse

Next, in order to improve the accuracy, the numerical solution is now based on stated reconstructed using a piecewise linear approximation, as follows:
\bel{eq:UPLR1}
{v}_{j-1/2,k}^+=\overline{{v}}_{j,k} - \frac{\Delta{x}}{2}\delta_{j,k}^x,
\quad\quad
{v}_{j+1/2,k}^-=\overline{{v}}_{j,k} + \frac{\Delta{x}}{2}\delta_{j,k}^x,
\ee
\be
{v}_{j,k-1/2}^+=\overline{{v}}_{j,k} - \frac{\Delta{y}}{2}\delta_{j,k}^y,
\quad\quad
{v}_{j,k+1/2}^-=\overline{{v}}_{j,k} + \frac{\Delta{y}}{2}\delta_{j,k}^y, 
\label{eq:UPLR2} 
\ee
where the following limiters are used:
\bel{eq:LMTRx} 
\delta_{j,k}^x\Delta{x} =
\begin{cases}
\sgn\left(v_{j+1,k} - v_{j-1,k}\right) 
\min \Big(2|v_{j,k}-v_{j-1,k}|, 2|v_{j+1,k}-v_{j,k}|, \left| v_{j+1,k}-v_{j-1,k}\right|/2 \Big),
\\
0,\qquad\qquad\,\,\,\,\,  \text{otherwise}, 
\end{cases}
\ee
\bel{eq:LMTRy} 
\delta_{j,k}^y\Delta{y} =
\begin{cases}
\sgn\left(v_{j,k+1} - v_{j,k-1}\right) \min
\Big(2|v_{j,k}-v_{j,k-1}|, 2|v_{j,k+1}-v_{j,k}|, \left| v_{j,k+1}-v_{j,k-1} \right|/2 \Big),
\\
0,\qquad\qquad\,\,\,\,\,  \text{otherwise}.
\end{cases}
\ee

\ese
Furthermore, the time-dependent ODE  implied by \eqref{eq:2+1Burger's - Flux-Source} is integrated in time, 
by using a stable and accurate
ODE solver. We use a fourth-order Runge-Kutta solver for the expanding case, 
and a third-order strong stability preserving (SSP) Runge-Kutta solver for the contracting case. The time step restriction is constrained by the CFL condition ($\Delta x = \Delta y$):
\be
\frac{\Delta{\tau}}{\Delta{x}} \max_{j,k} |v_{j\pm\frac{1}{2}, k}^{\pm}|
\leq \frac{1}{2}.
\label{eq:CFL}  
\ee

\subsection{Asymptotic behavior on an expanding background} 
\label{ExpandingBurgers}

\begin{figure}[htbp]
\centering
  \subcaptionbox{$\tau=1$}{\includegraphics[height=2.2in]{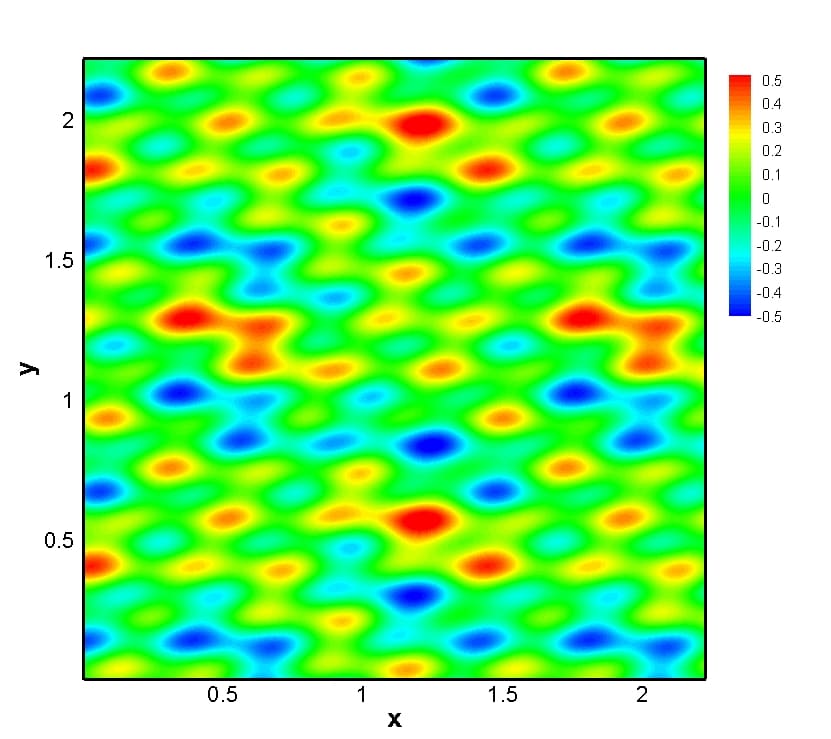}}
  \subcaptionbox{$\tau=16$}{\includegraphics[height=2.2in]{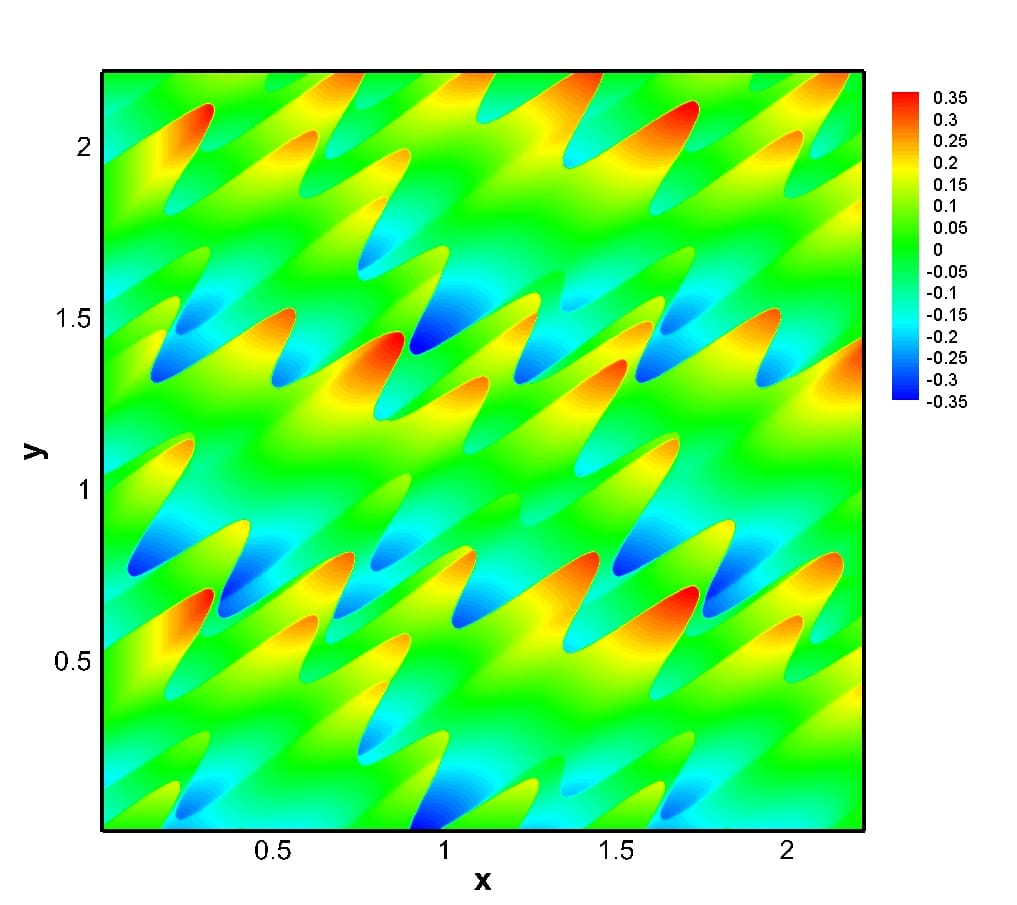}} \\
  
\centering
  \subcaptionbox{$\tau=128$}{\includegraphics[height=2.2in]{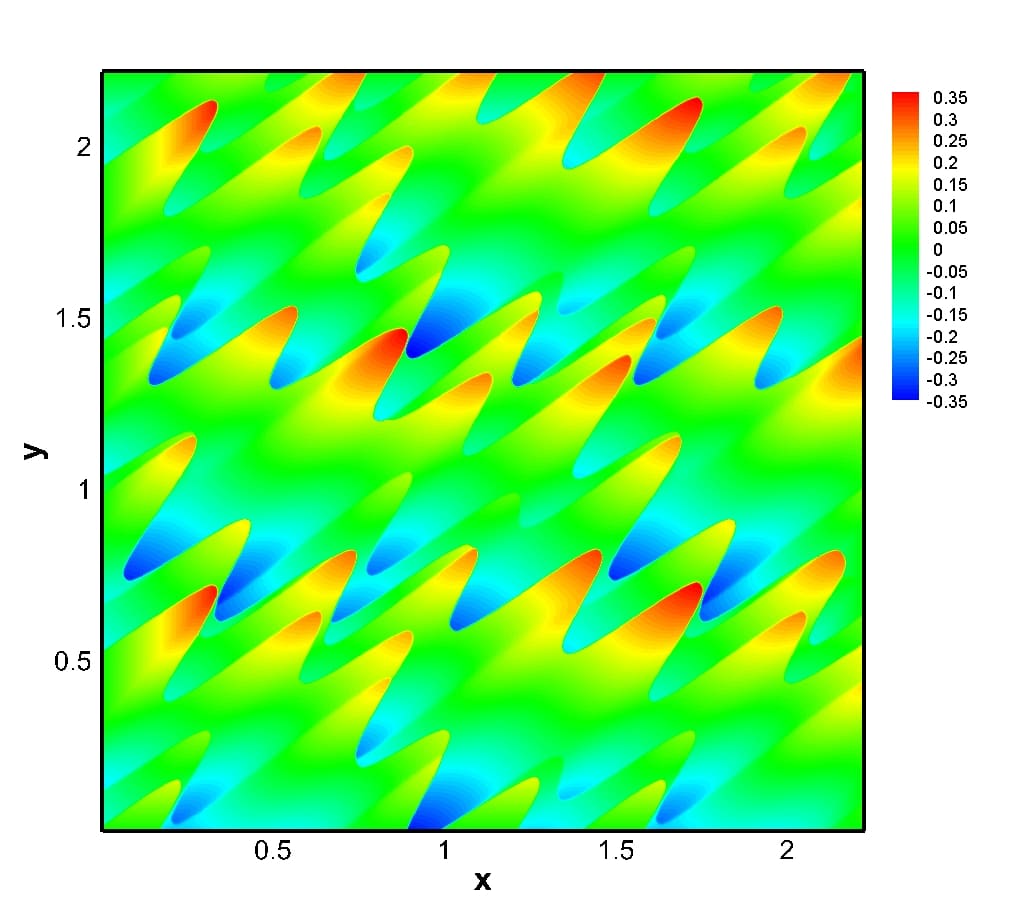}}
  \subcaptionbox{$\tau=512$}{\includegraphics[height=2.2in]{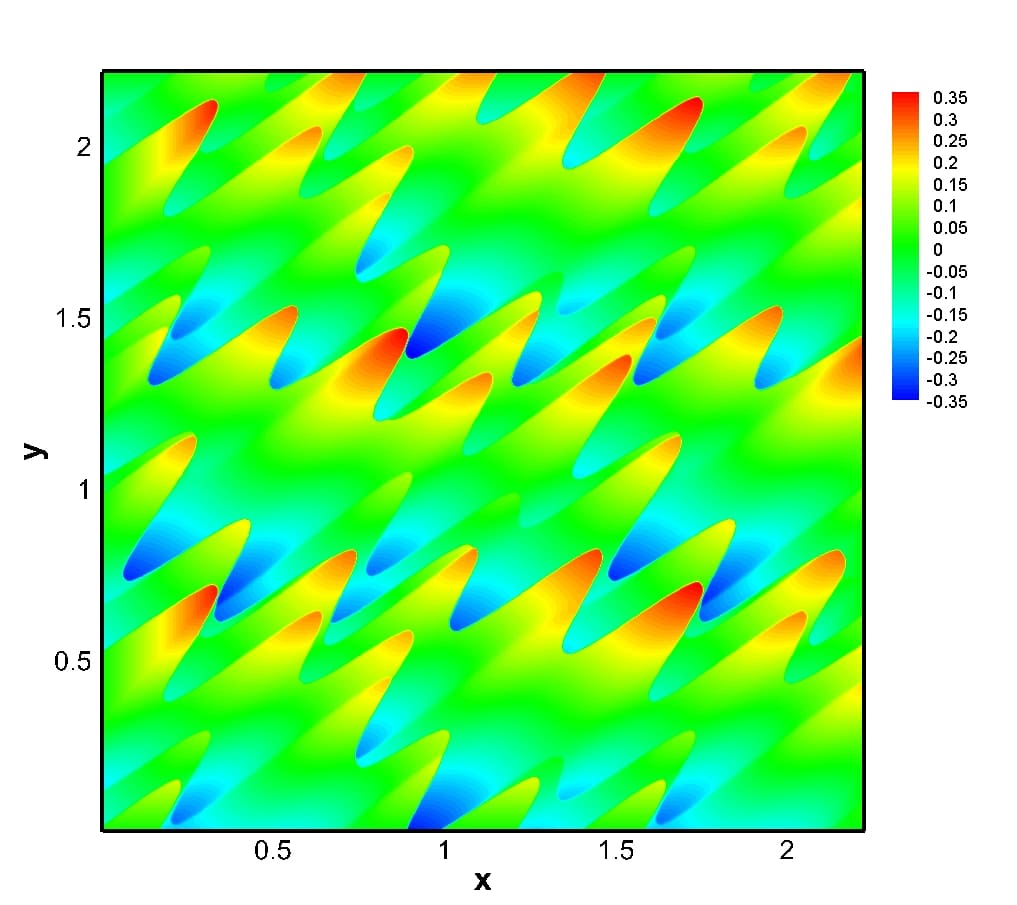}}
    
\caption{2-D contours of the rescaled velocity $w$ with an $[800\times800]$ grid and $\kappa = 2$.
 (a) Initial condition at $\tau=1$.
 (b) $\tau=16$.
 (c) $\tau=128$.
 (d) $\tau=512$.}
\label{fig:ex2D}
\end{figure}

\begin{figure}[htbp]
\centering
  \subcaptionbox{$\tau=1$}{\includegraphics[height=2.4in]{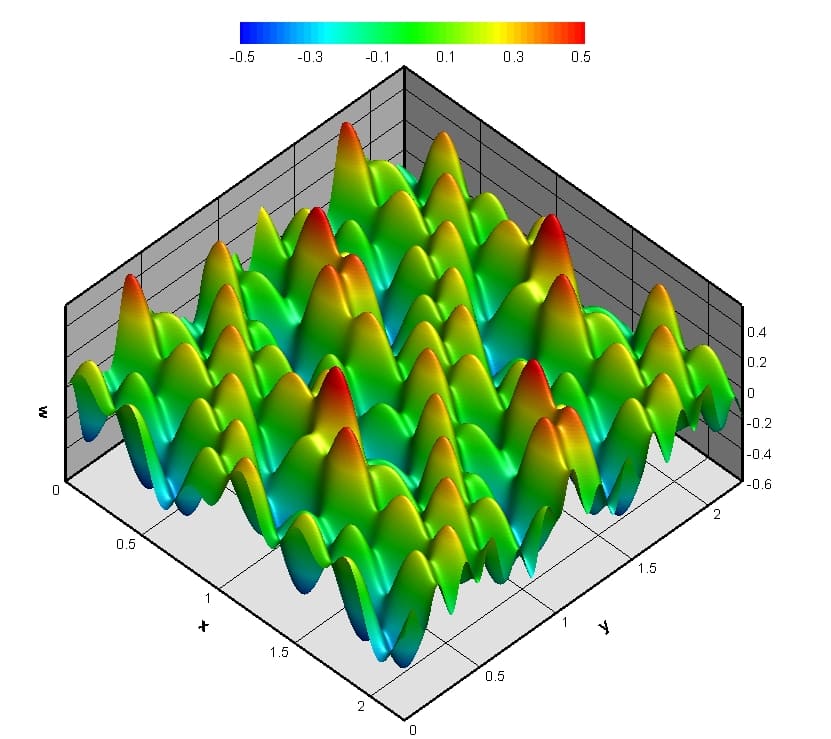}}
  \subcaptionbox{$\tau=512$}{\includegraphics[height=2.4in]{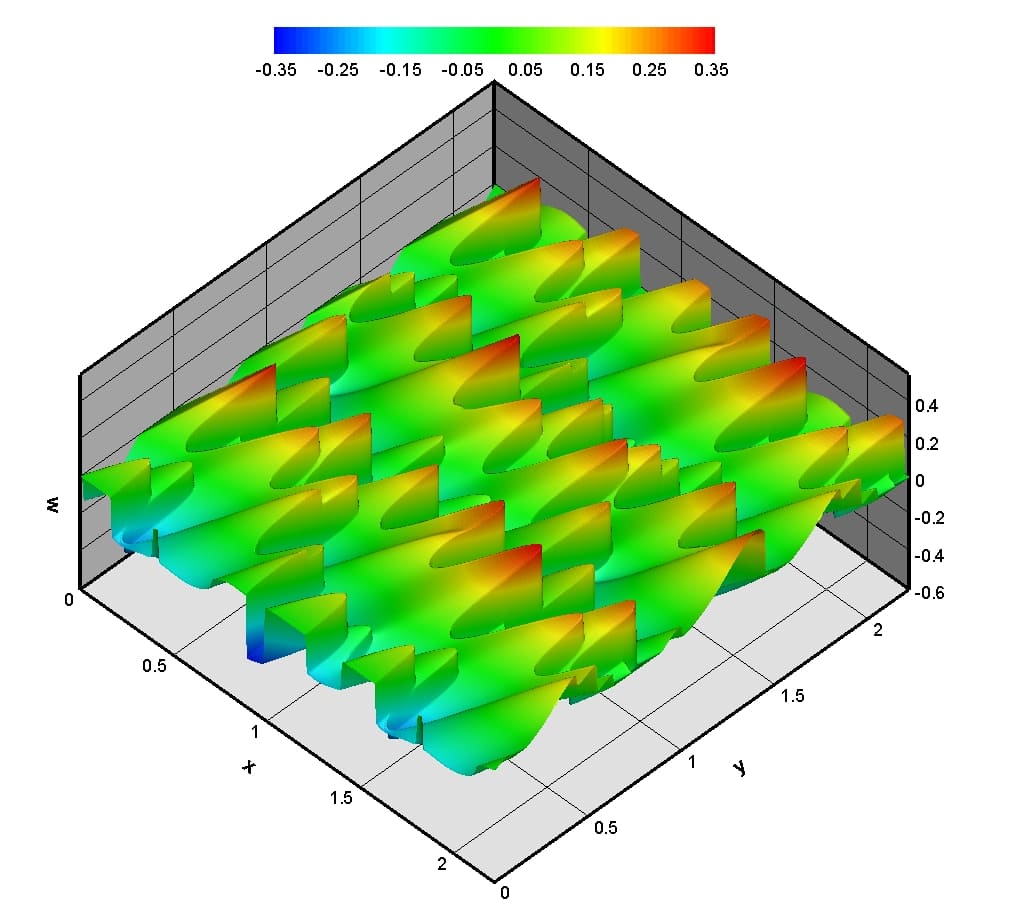}}
  
\caption{3-D contours of rescaled velocity $w$ with an $[800\times800]$ grid and $\kappa = 2$. 
(a) Initial condition at $\tau=1$.
 (b) $\tau=512$.}
\label{fig:ex3D}
\end{figure}

In this section and the following one, we present numerical tests in which the dynamics of asymptotic solutions for the expanding and contracting cases. The computational domain is $[0,\pi/\sqrt{2}] \times [0,\pi/\sqrt{2}]$. For both tests, we set $\kappa=2$ and $4$ and choose some ``arbitrary'' initial condition at $\tau_0$, specifically 
\be
\begin{split}
v(\tau_0, x, y) =  &\frac{1}{8}\left(
 \sin \left(4\sqrt{2}\pi x - 3\sqrt{2} \pi y \right) + \cos \left(\sqrt{2}\pi x + 3\sqrt{2} \pi y \right)  +  \sin \left(3\sqrt{2}\pi x - 5\sqrt{2} \pi y \right) \right)
 \\
+ &\frac{1}{8} \left(
\sin \left(5\sqrt{2}\pi x + 3\sqrt{2} \pi y \right) - \cos \left(2\sqrt{2}\pi x + 2\sqrt{2} \pi y \right)  \right).
\label{eq:IC2D} 
\end{split}
\ee
This example displays the dynamics of the $(2+1)$-dimensional cosmological Burgers in an expanding spacetime. Three different grid refinements $[200\times200]$, $[400\times400]$, $[800\times800]$ are chosen to be able to compare the error of different solutions (Section \ref{Convergence}).
We observe that for relatively large $\tau$ the source term of the cosmological Burgers becomes unstable, as $v(\tau, x, y) \to 0$. This increases $\tau$ and $\Delta\tau$ significantly. Hence, we introduce our second stability condition as follows:
\be
\Delta\tau_n \leq {1 \over \kappa}
 \min_{j,k} 
\Big( \frac{\tau_n}{(1 - (v_{j,k}^n)^2)}\Big). 
\label{eq:StabilityEx-deux}  
\ee
The solutions decay to zero uniformly at the rate of $\tau^{-\kappa}$ and, therefore, the rescaled solution $w=v\tau^{\kappa}$ approaches a non-trivial limit as $\tau$ goes $+\infty$. This can be seen in Figures~\ref{fig:ex2D} and~\ref{fig:ex3D}, in 2-D and 3-D respectively.


\begin{figure}[htbp]
\centering
  \subcaptionbox{ 1S1T: $ L^1= 5.3663 \times 10^{-5} $}{\includegraphics[height=2.2in]{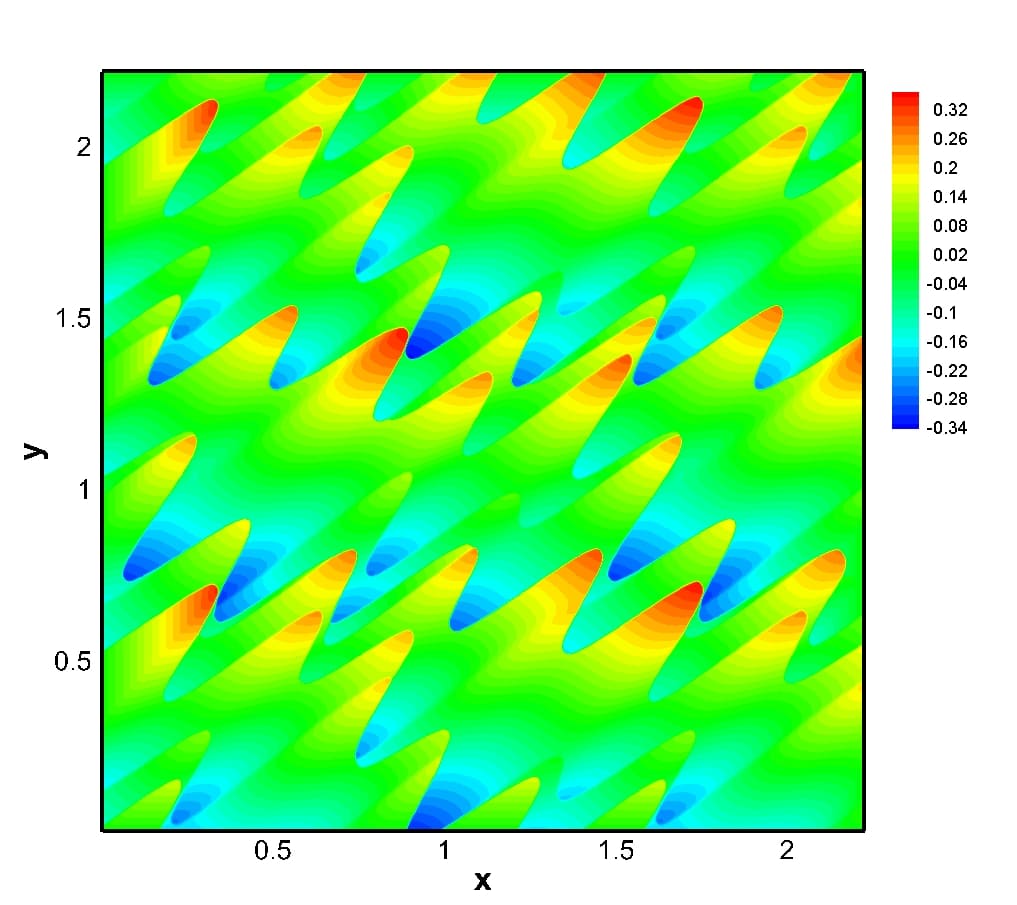}}
  \subcaptionbox{1S4T: $ L^1= 1.0549  \times 10^{-6}$}{\includegraphics[height=2.2in]{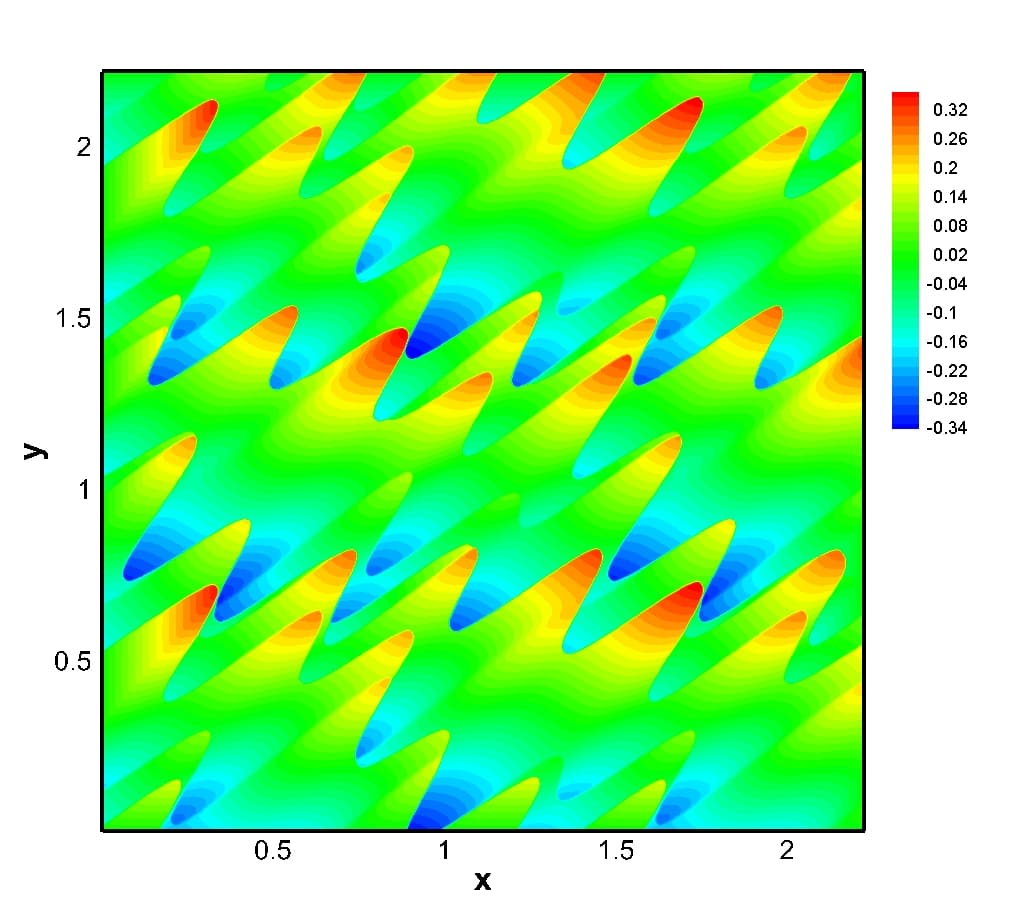}}
 
\centering
  \subcaptionbox{ 2S1T: $L^1= 5.2826 \times 10^{-5}$}{\includegraphics[height=2.2in]{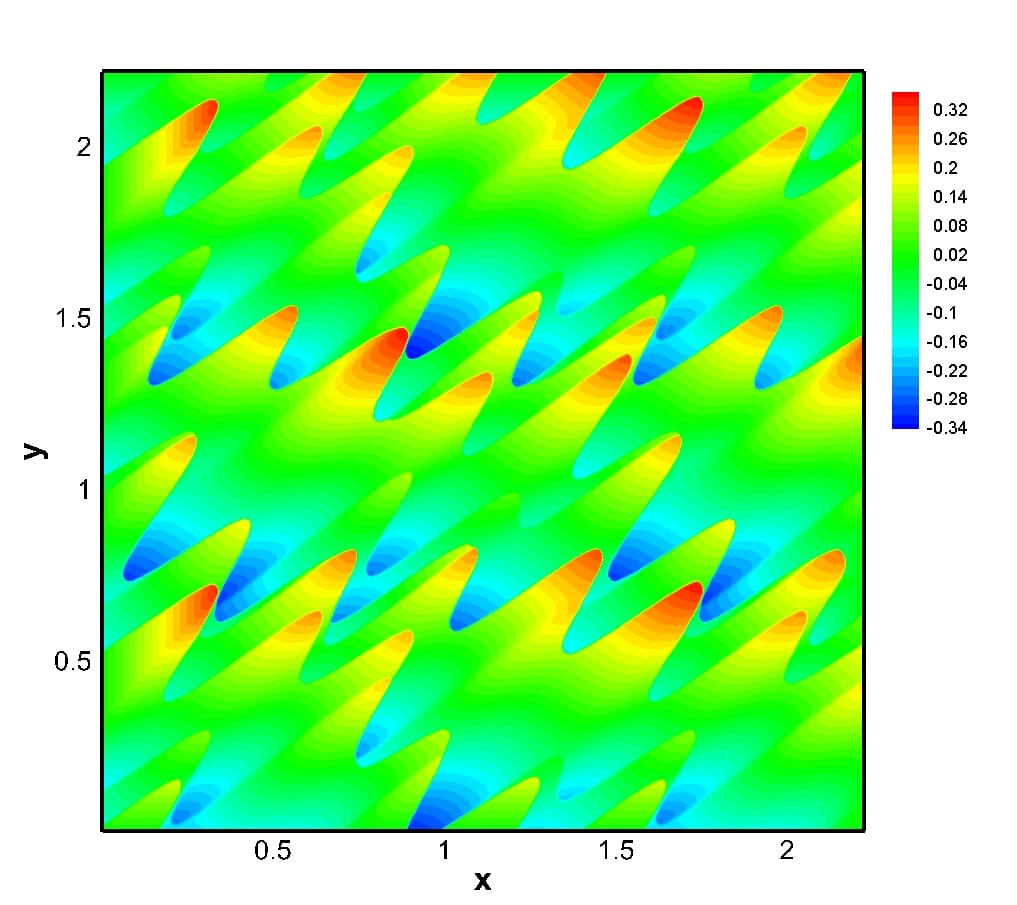}}
  \subcaptionbox{2S4T: $L^1=0$}{\includegraphics[height=2.2in]{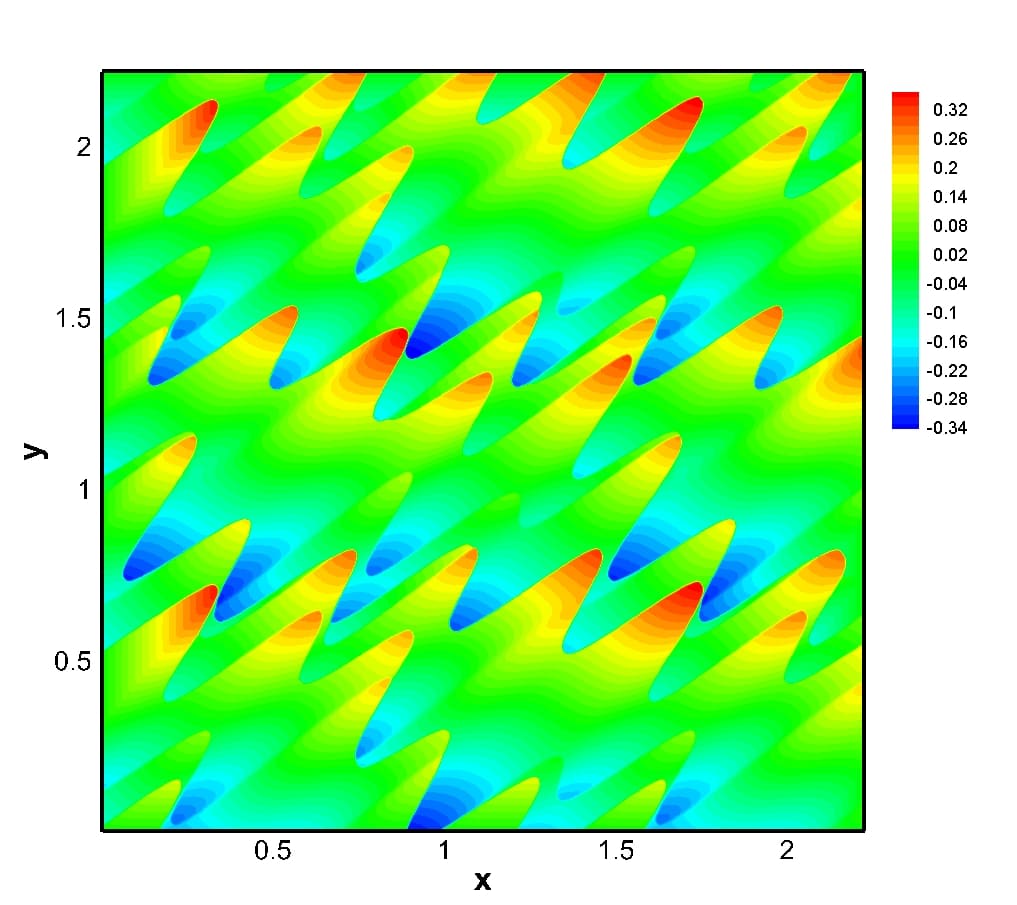}}
   
\caption{2-D contours of rescaled velocity $w$ with an $[800\times800]$ grid, $\kappa = 2$ at $\tau = 1024$ with different orders of spatial and temporal discretization.}
\label{fig:ex2D-SCH}
\end{figure}

Based on the 1-D tests, we choose a second-order spatial and fourth-order temporal discretization (2S4T). In addition,  in order to analyze the effect on the solution we run this example with the following discretizations: first-order space and first-order time (1S1T), first-order space and fourth-order time (1S4T), and second-order space and first-order time (2S1T). Figure~\ref{fig:ex2D-SCH} shows the solutions of this test with above-mentioned schemes. Furthermore, we compute the $L^1$ norm for these schemes based on the best scheme (2S4T). The results show that the $L^1$ norms are small. Increasing the order of temporal discretization is more effective than increasing the spatial order. Moreover, it  can be concluded that the lower-order spatial schemes can be used to be able to reduce computational cost. 


\begin{figure}[htbp]
\centering
  \subcaptionbox{1S1T: $ L^1(v)= 1.5056 \times 10^{-4}$}{\includegraphics[height=2.2in]{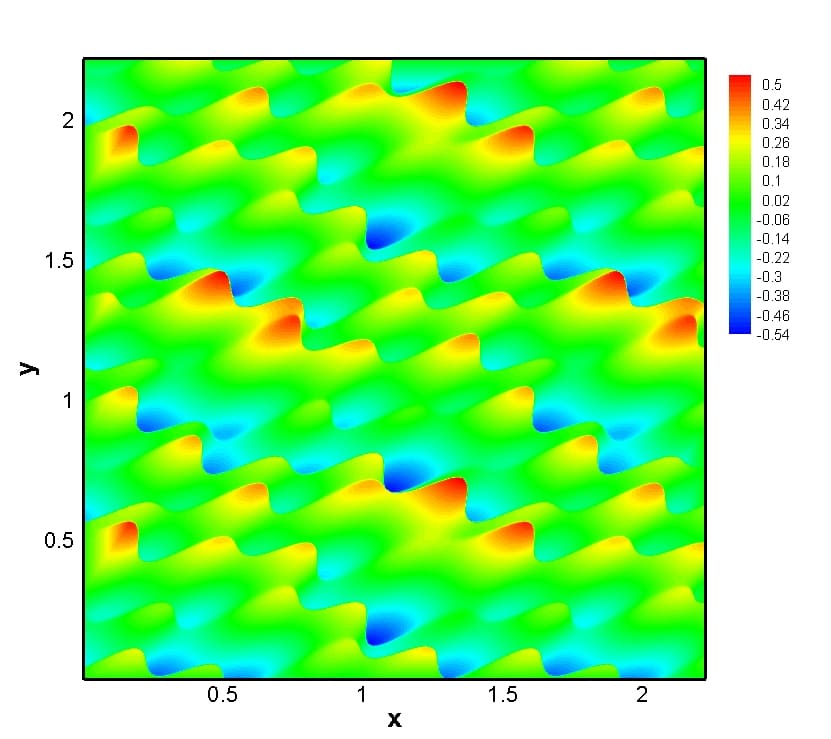}}
  \subcaptionbox{1S4T: $ L^1(v)= 6.0099 \times 10^{-7}$}{\includegraphics[height=2.2in]{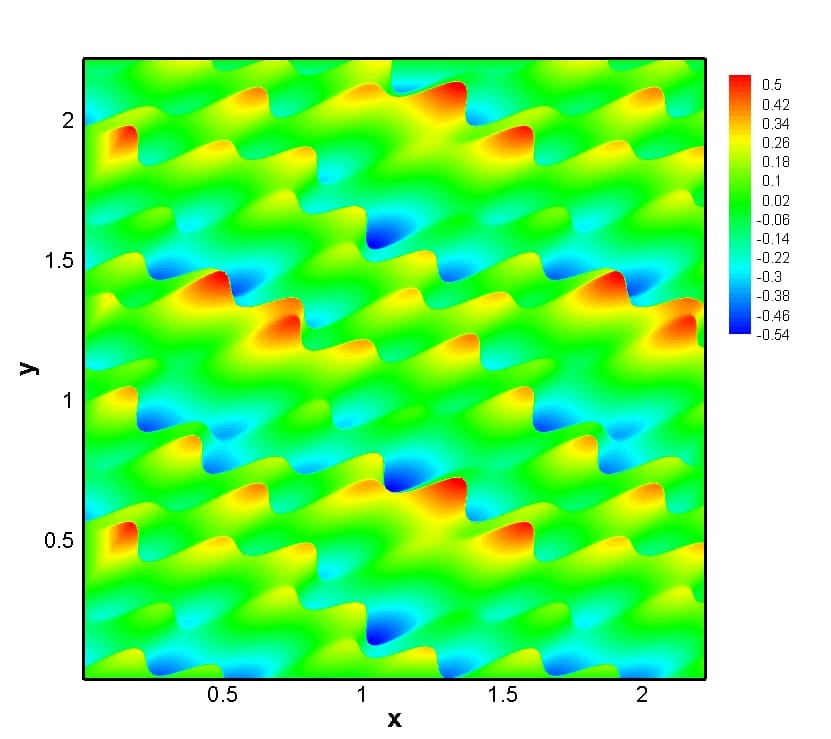}}
  
\centering
  \subcaptionbox{2S1T: $L^1(v)= 1.5038 \times 10^{-4}$}{\includegraphics[height=2.2in]{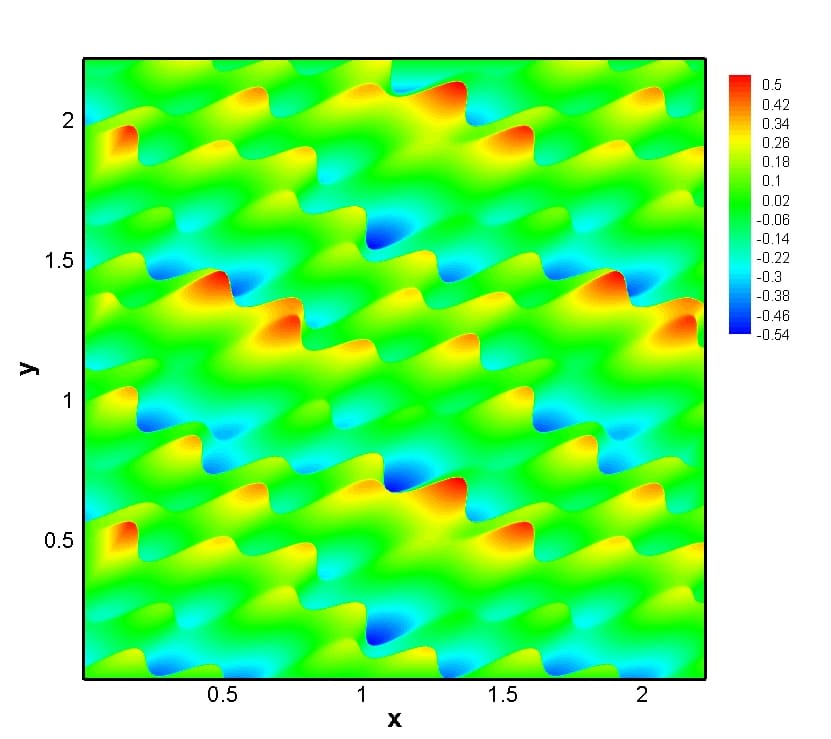}}
  \subcaptionbox{2S4T: $L^1(v)= 0$}{\includegraphics[height=2.2in]{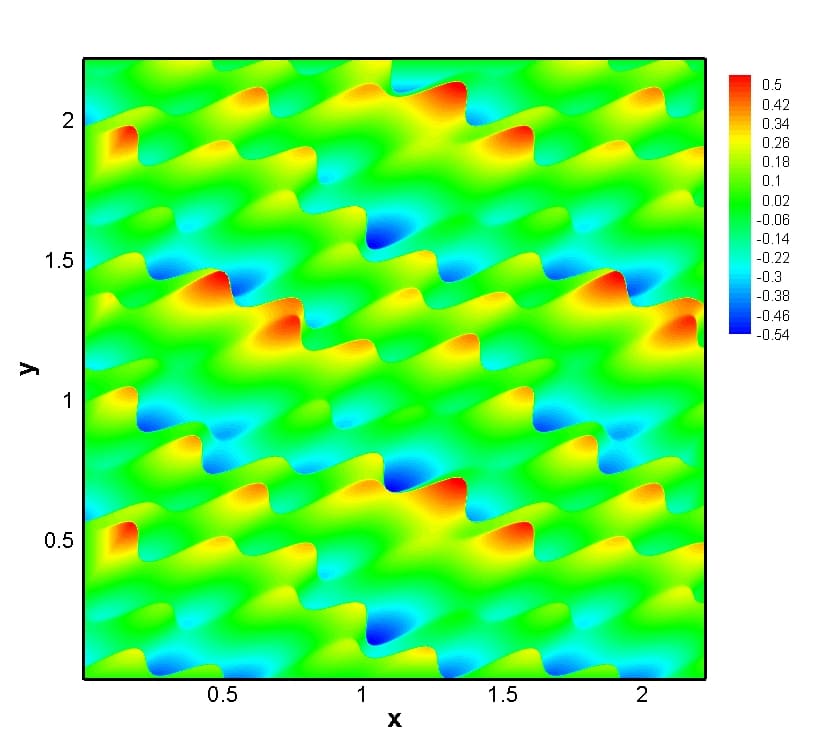}}
   
\caption{2-D contours of the rescaled velocity $w$ for an $[800\times800]$ grid, $\kappa = 4$ at $\tau = 1024$ with different orders of spatial and temporal discretization.}
\label{fig:exk4}
\end{figure}

\begin{figure}[htbp]
\centering
  \subcaptionbox{$\tau = 16$}{\includegraphics[height=2.2in]{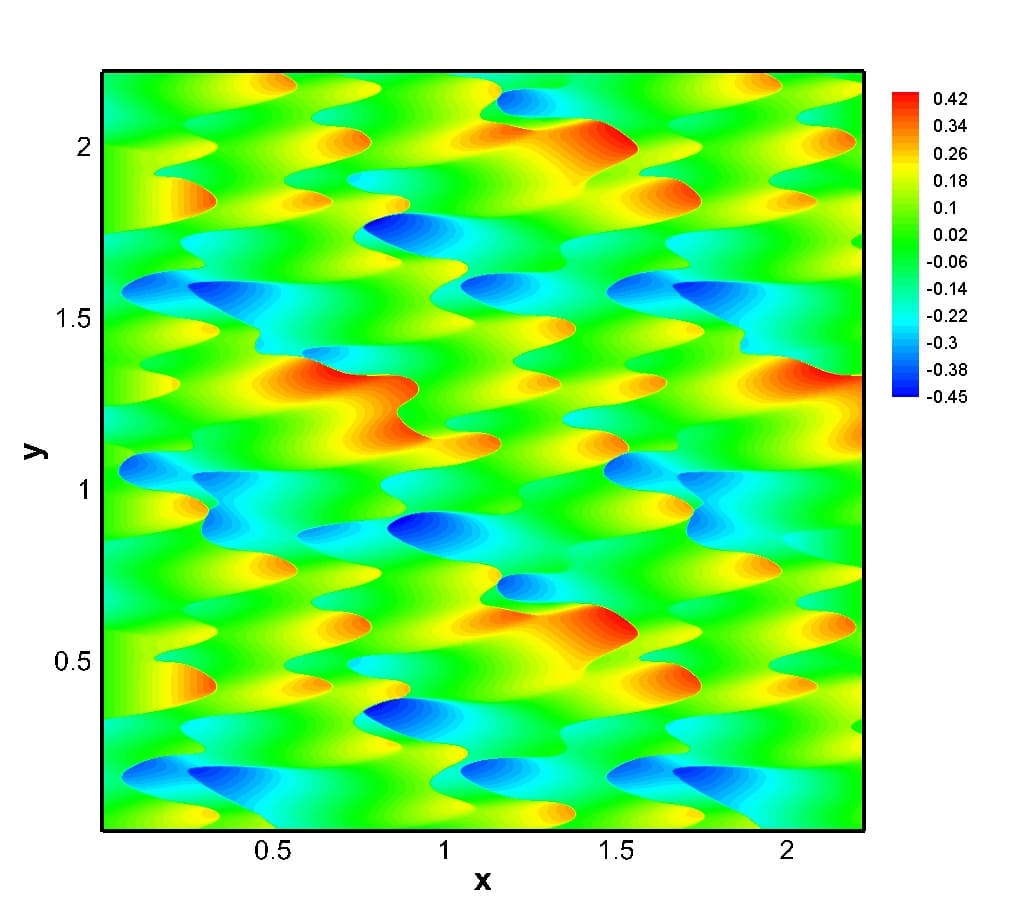}}
  \subcaptionbox{$\tau = 1024$}{\includegraphics[height=2.2in]{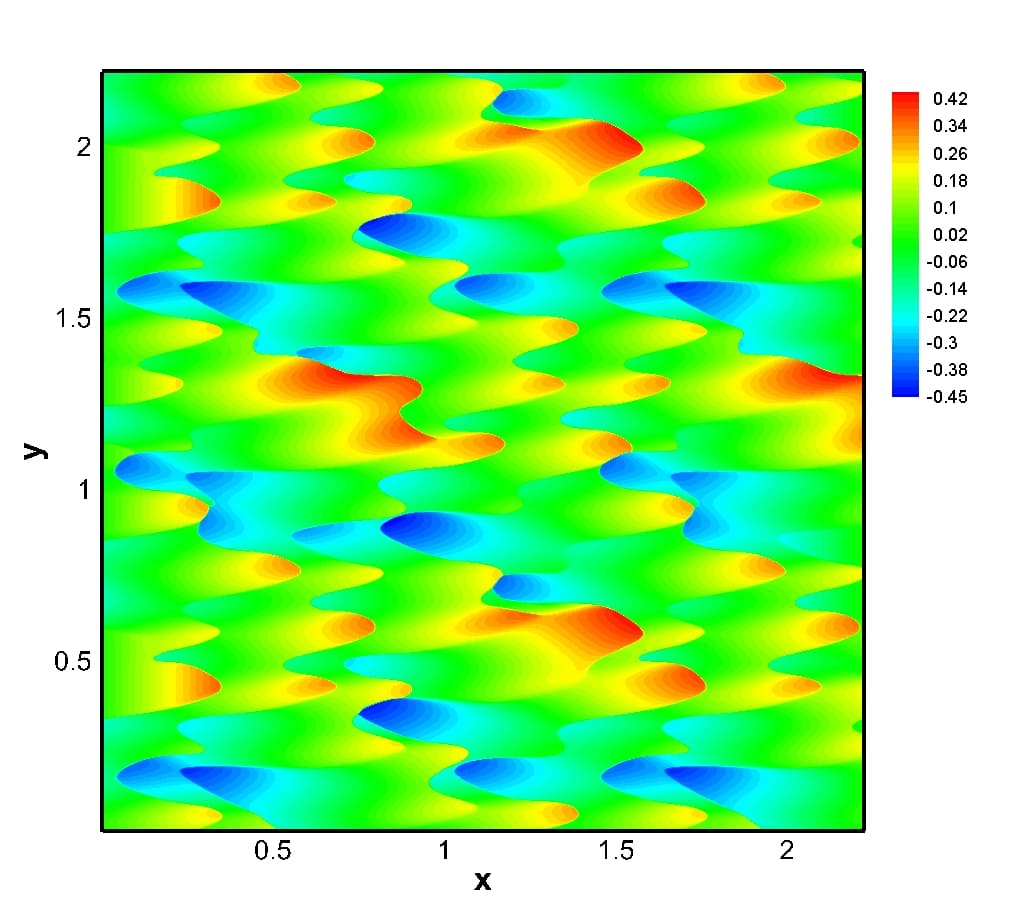}}
  
\centering
  \subcaptionbox{$\tau = 16$}{\includegraphics[height=2.2in]{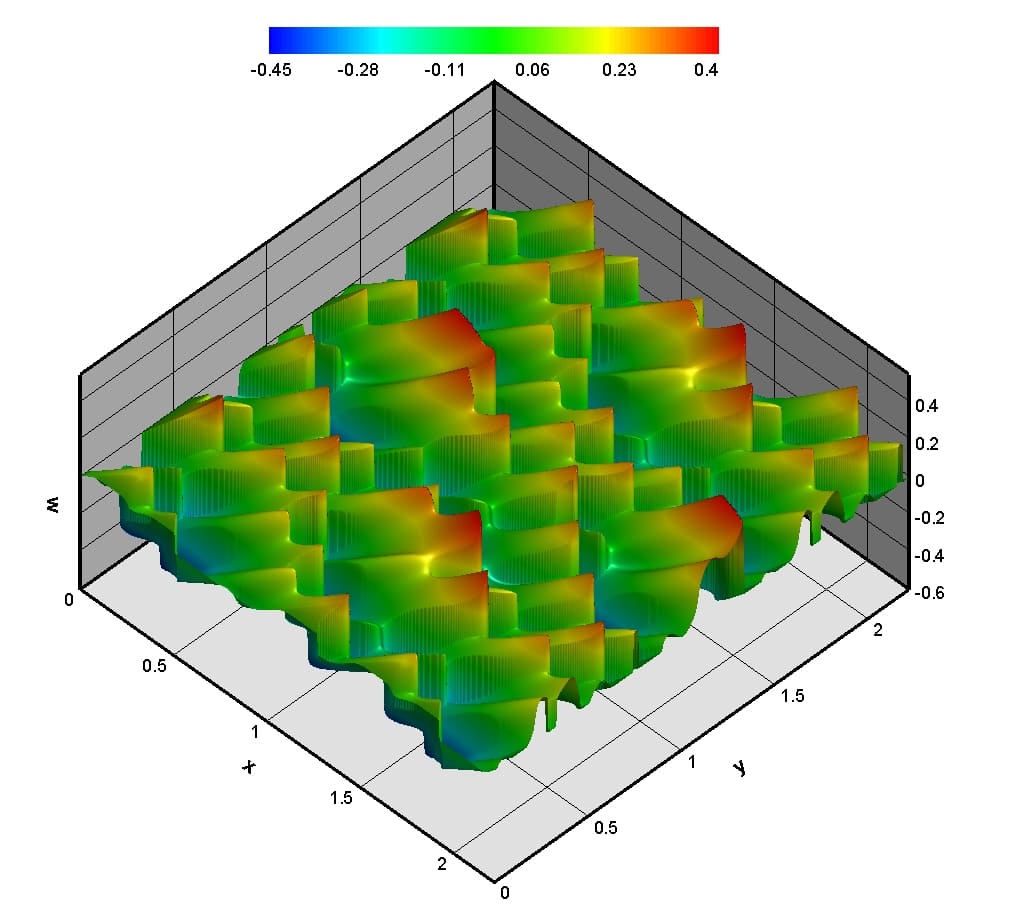}}
  \subcaptionbox{$\tau =1024$}{\includegraphics[height=2.2in]{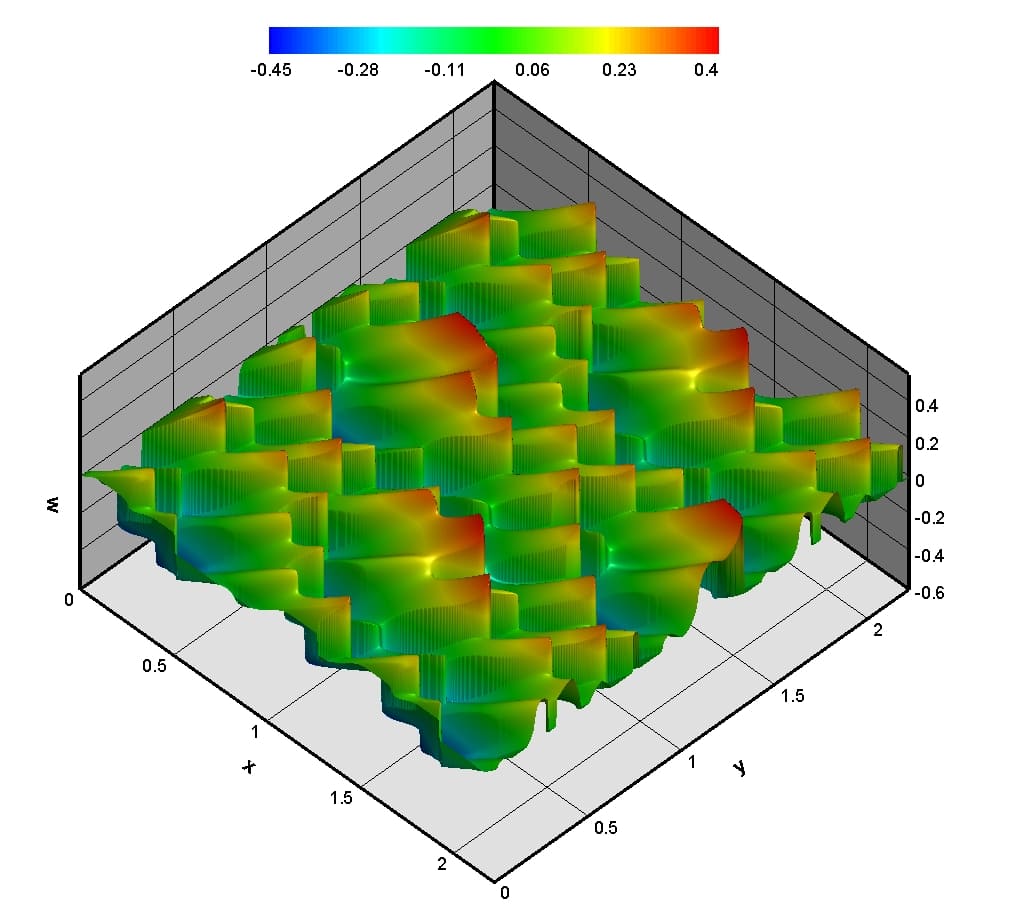}}
   
\caption{2-D and 3-D contours of the rescaled solution $w$ with flux $g(v) = v^3 /2$ at two different times.}
\label{fig:ex-g=u31}
\end{figure}

\begin{figure}[htbp]
\centering
  \subcaptionbox{$\tau = 16$}{\includegraphics[height=2.2in]{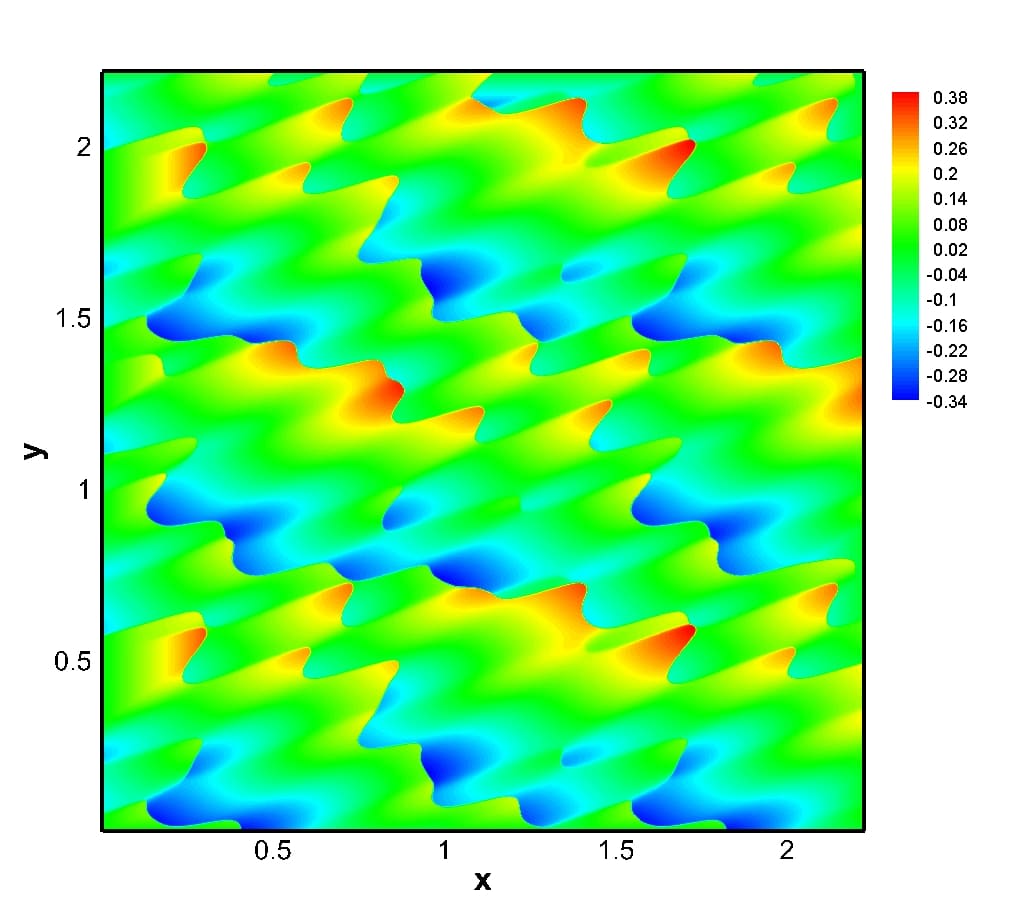}}
  \subcaptionbox{$\tau = 1024$}{\includegraphics[height=2.2in]{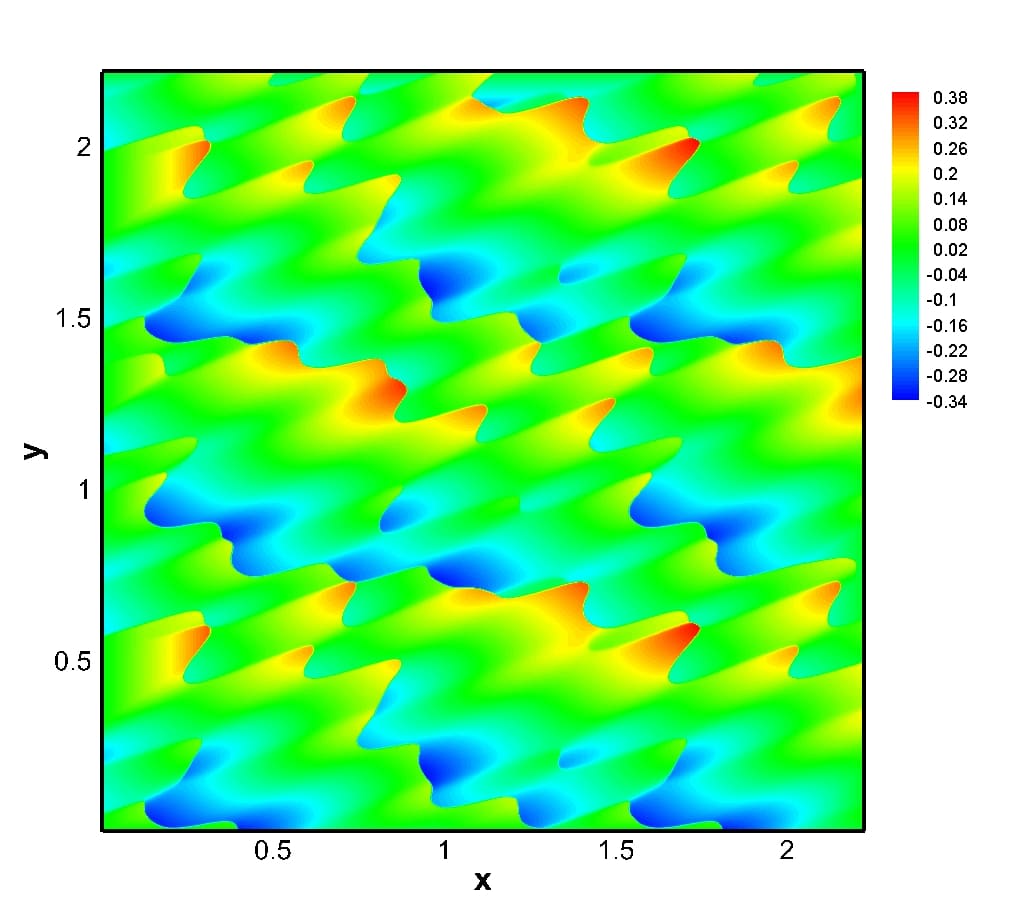}}
  
\centering
  \subcaptionbox{$\tau = 16$}{\includegraphics[height=2.2in]{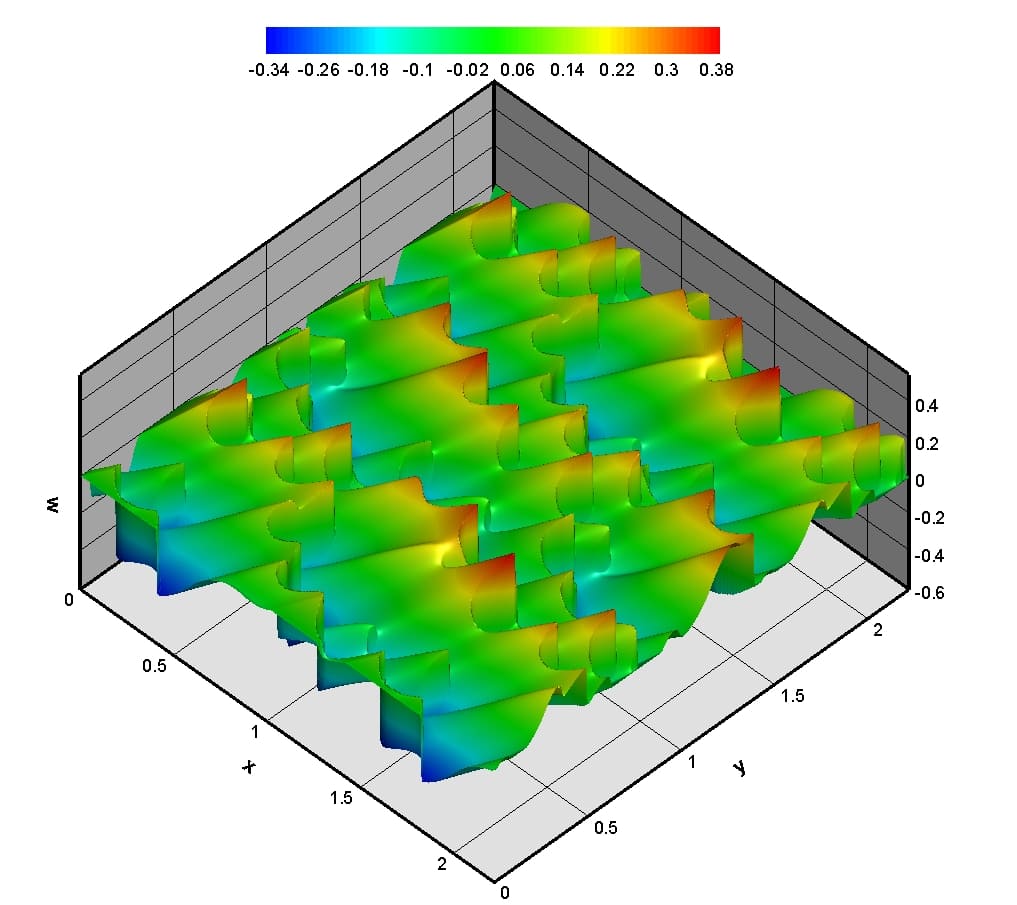}}
  \subcaptionbox{$\tau =1024$}{\includegraphics[height=2.2in]{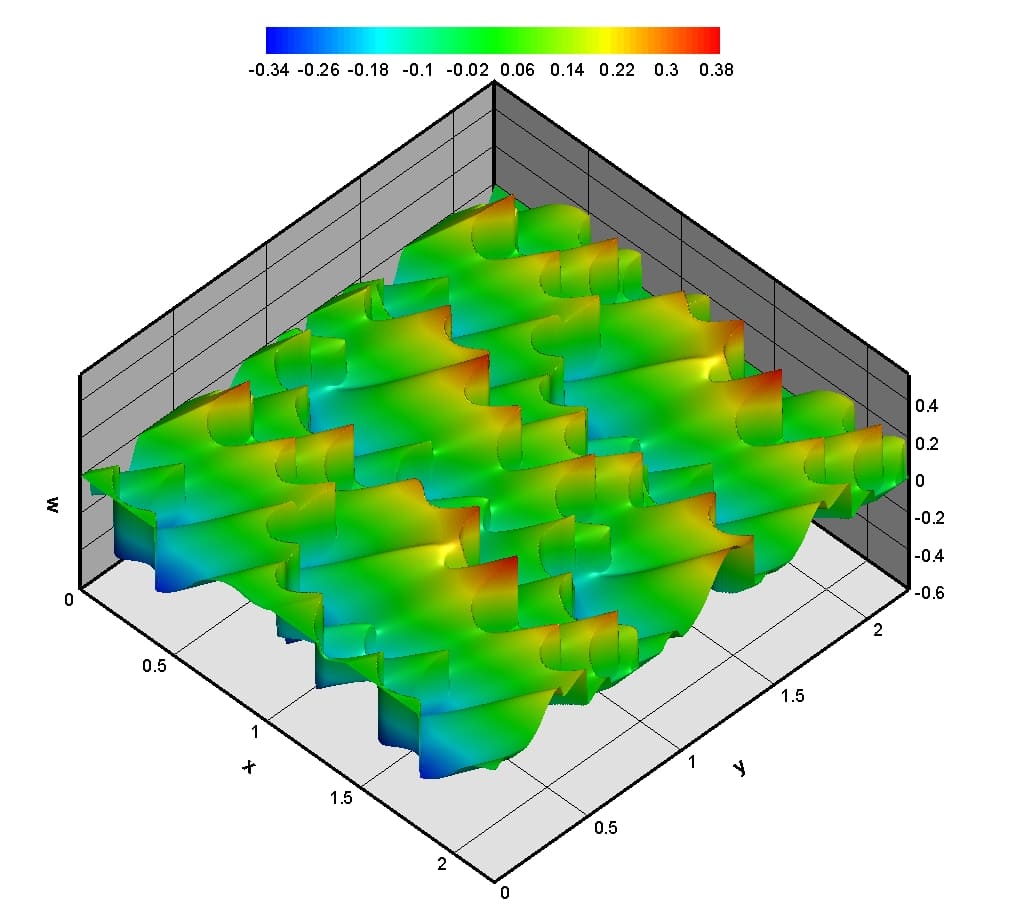}}
   
\caption{2D and 3D contours of the solution with flux $g(v) = (1-\beta) v^2/2 + \beta v^3/3$  and $\beta=1/2$ at two different times.}
\label{fig:ex-g=u32}
\end{figure}

The expanding background test is also solved with $\kappa = 4$. The results show that velocity $v$ approaches to zero. Figure~\ref{fig:exk4} illustrates the rescaled velocity $w$ in four different schemes. Once again the $L^1$ norm for these solutions are very small and the order order temporal discretization is more effective than the order of spatial one. 
We also study another choice of flux functions in \eqref{eq:2+1Burger's - Flux-Source} as follows:
\bel{eq:gv3}
f(v) =  {1 \over 2} v^2, \qquad g(v) = {1 \over 2} v^3, 
\ee
and
\bel{eq:gv32}
f(v) =  {1 \over 2} v^2, 
\qquad g(v) = {(1-\beta) \over 2} v^2 + {\beta \over 3} v^3, 
\ee
where $\beta \in (0, 1)$. $\beta = 1/2$ is chosen for the numerical tests. Observe that the Godunov fluxes in the $y$ direction change which is not presented for the sake of brevity. The effect of different fluxes can be seen in Figures~\ref{fig:ex-g=u31} and~\ref{fig:ex-g=u32}. Observe also that shock waves are formed in the $x$ direction.


\subsection{Asymptotic behavior on a contracting background} 
\label{ContractingBurgers}

\begin{figure}[htbp]
\centering
  \subcaptionbox{$\tau = -10^{-1}$}{\includegraphics[height=2.2in]{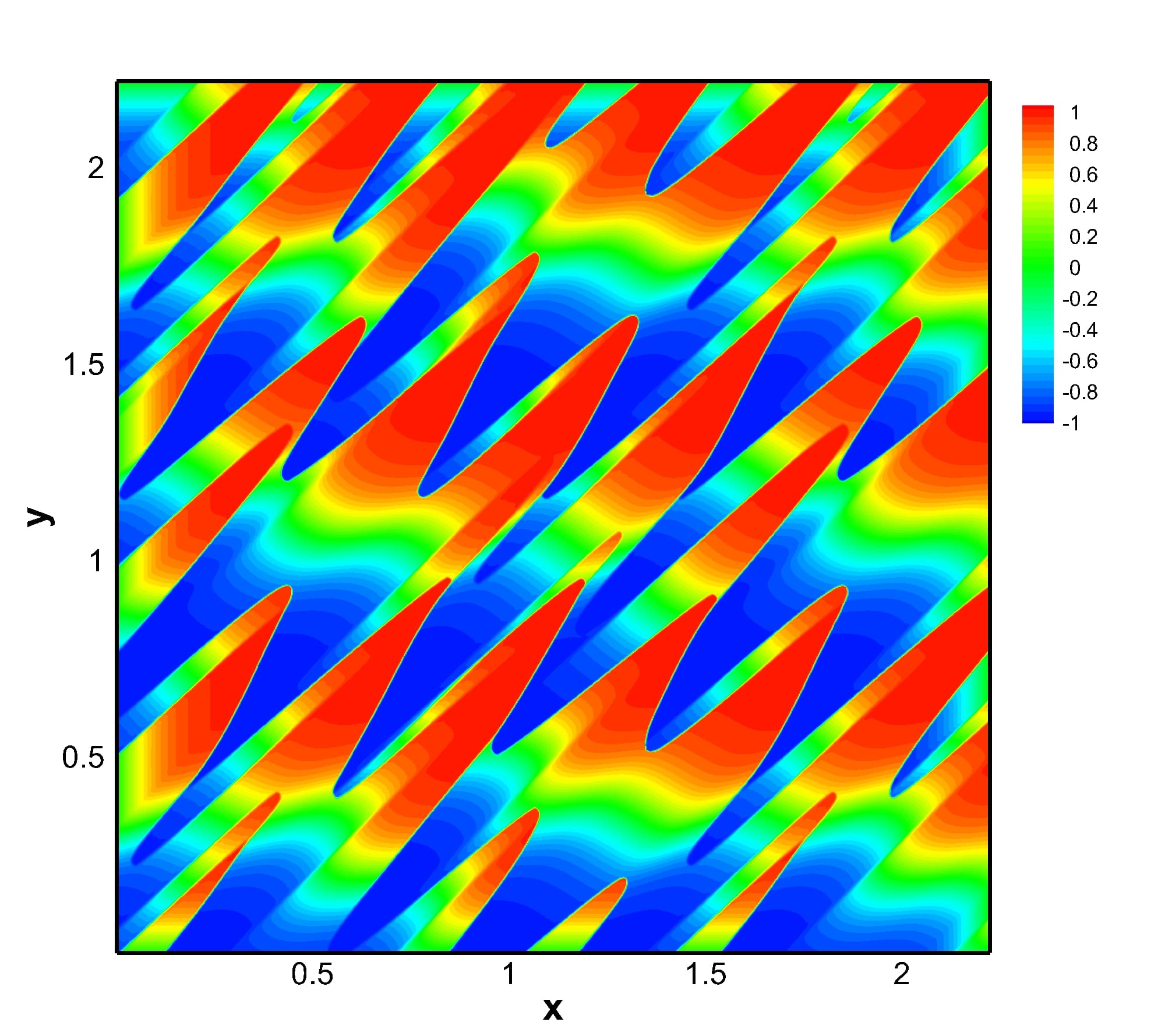}}
  \subcaptionbox{$\tau = -10^{-2}$}{\includegraphics[height=2.2in]{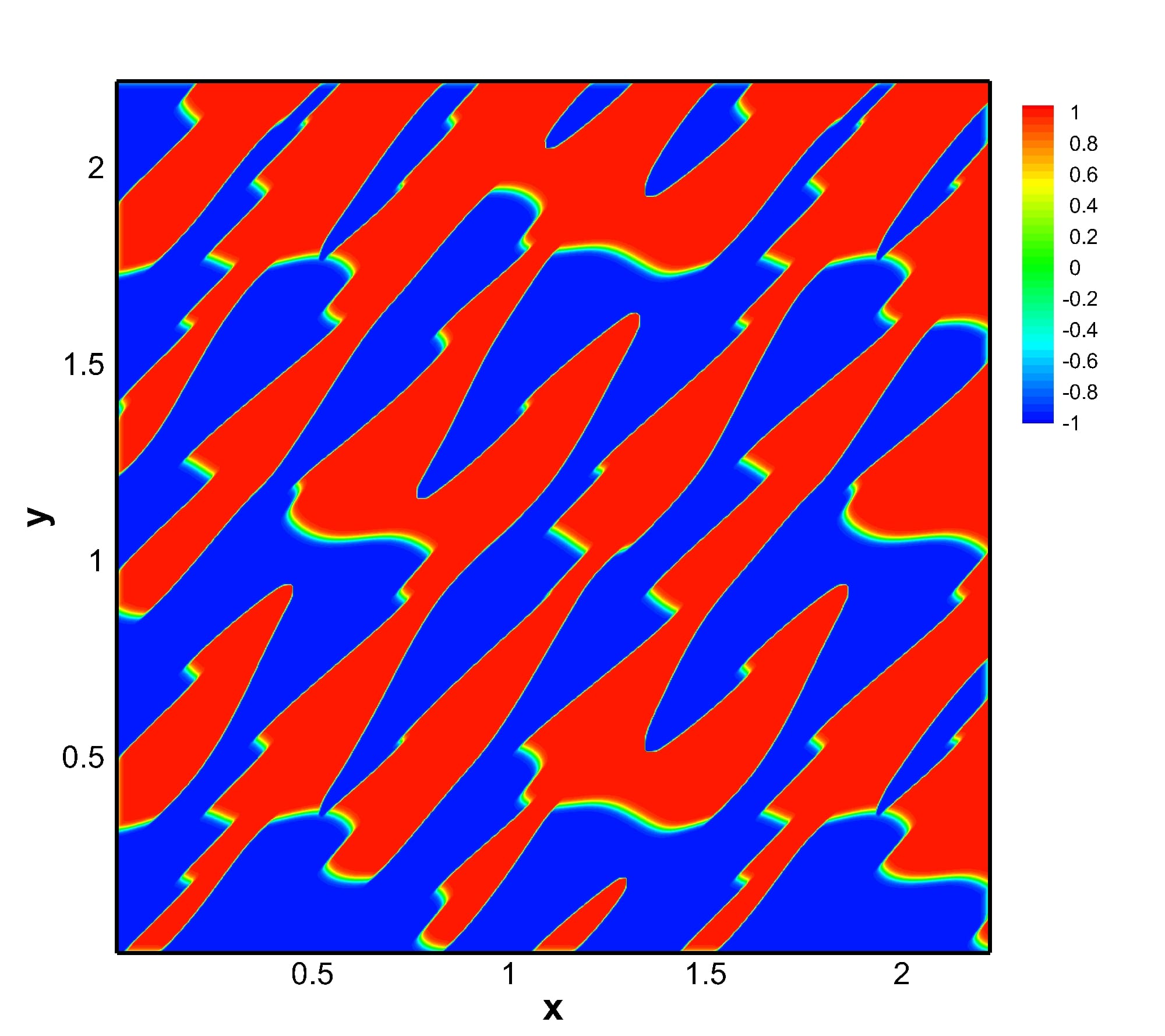}}
  
\centering
  \subcaptionbox{$\tau = -10^{-3}$}{\includegraphics[height=2.2in]{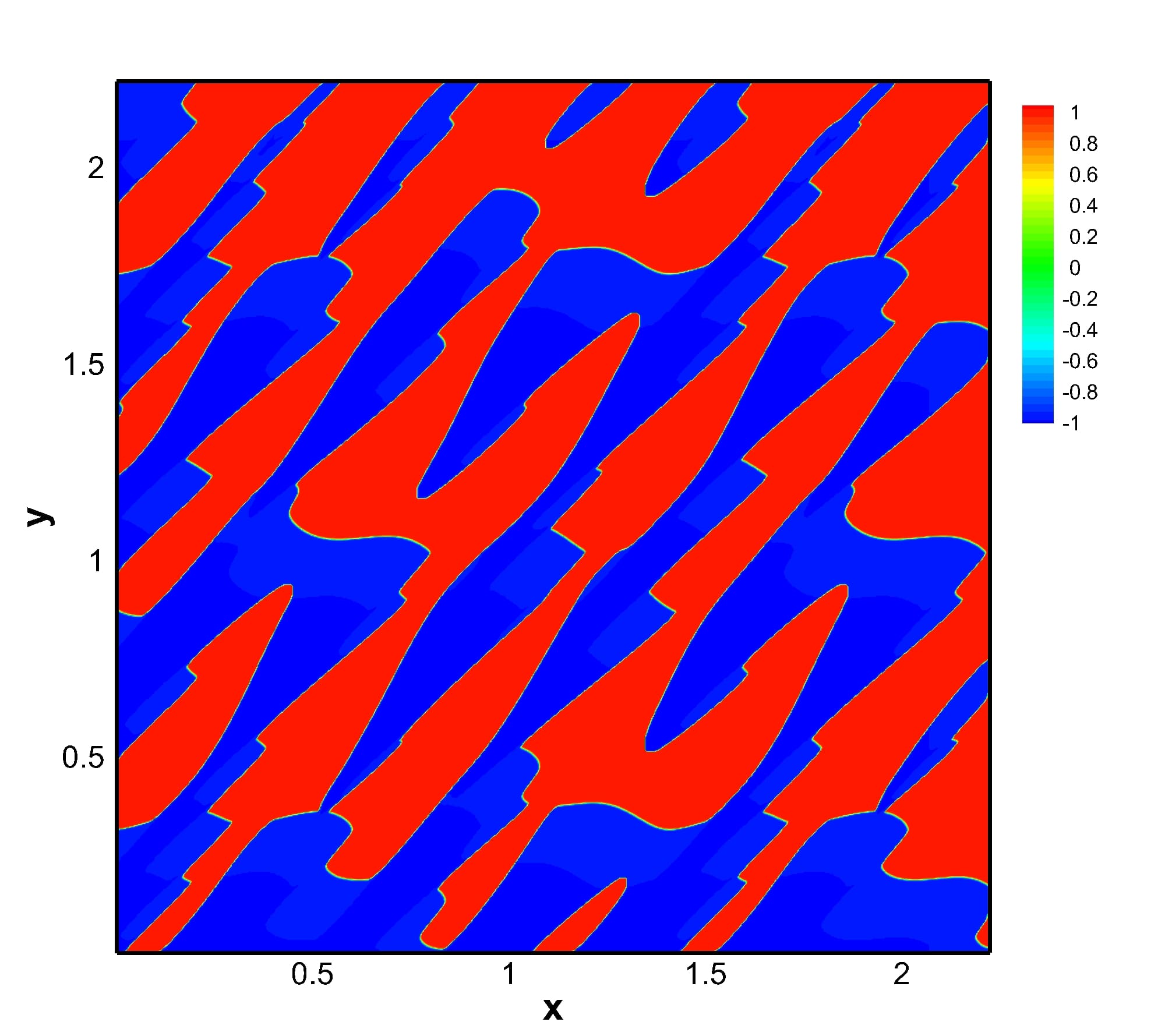}}
  \subcaptionbox{$\tau = -10^{-4}$}{\includegraphics[height=2.2in]{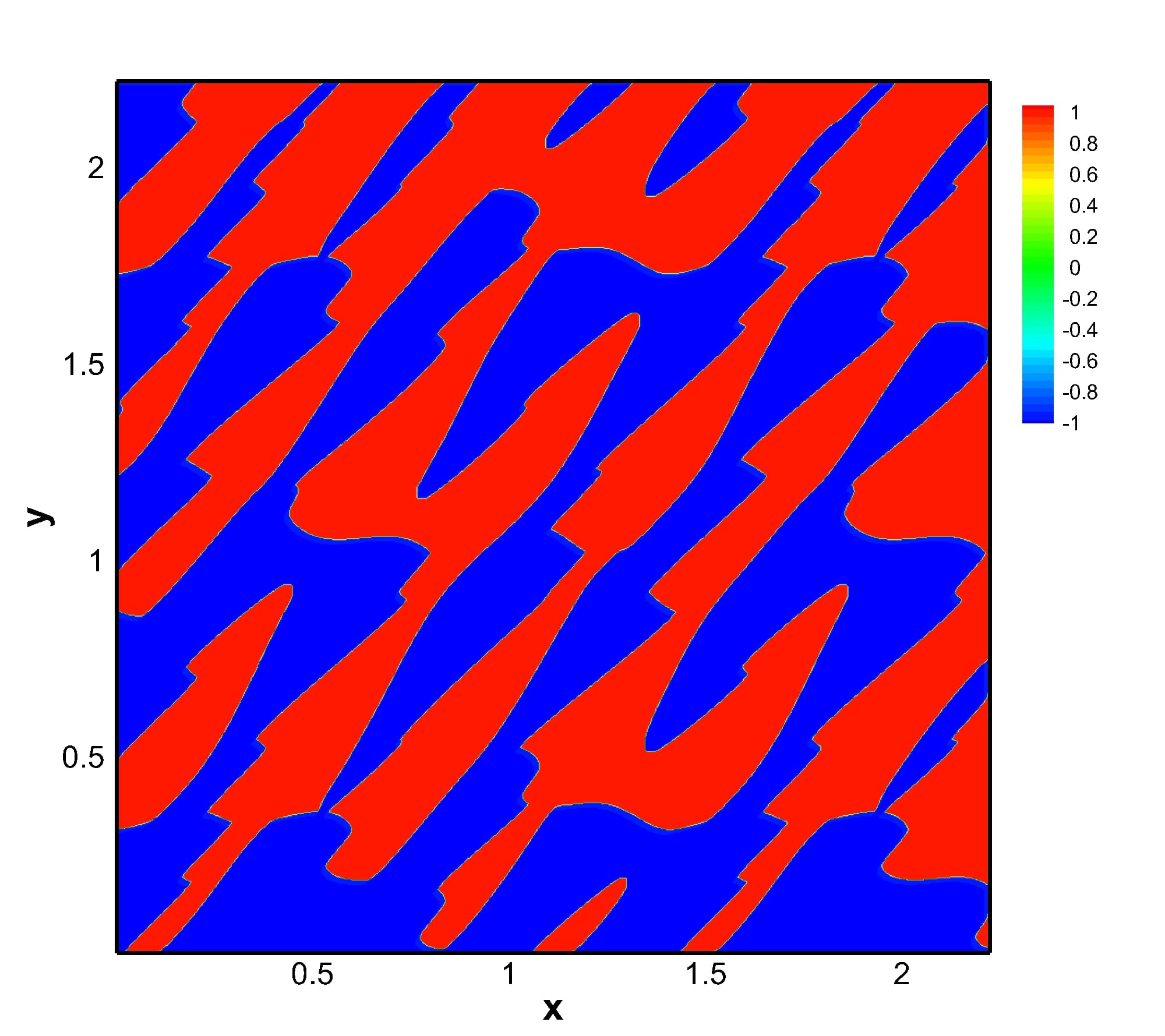}}
   
\caption{2-D contours of the velocity $v$ with an $[800\times800]$ grid and $\kappa = 2$ at different times. 
(a) $\tau = -10^{-1}$,
 (b) $\tau = -10^{-2}$. (c) $\tau =-10^{-3}$.
 (d) $\tau =-10^{-4}$.}
\label{fig:co2D}
\end{figure}

\begin{figure}[htbp]
\centering
  \subcaptionbox{$\tau = -10^{-1}$}{\includegraphics[height=2.2in]{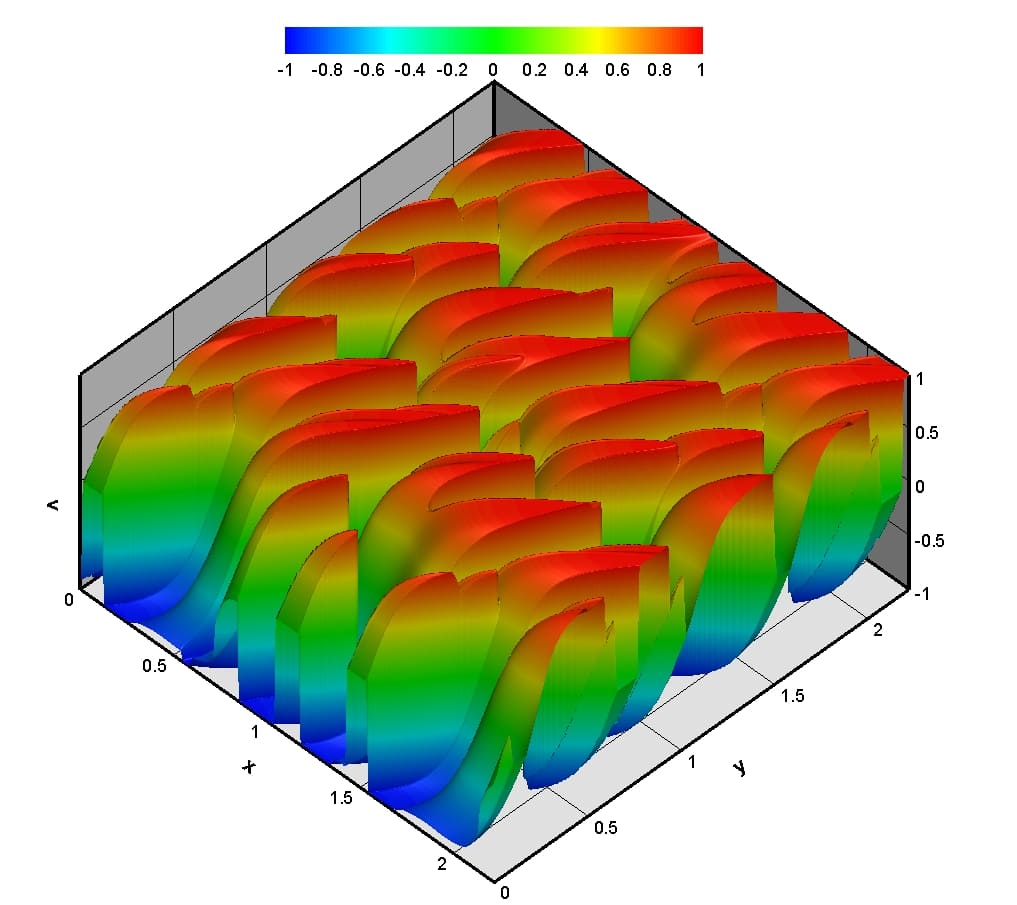}}
  \subcaptionbox{$\tau = -10^{-4}$}{\includegraphics[height=2.2in]{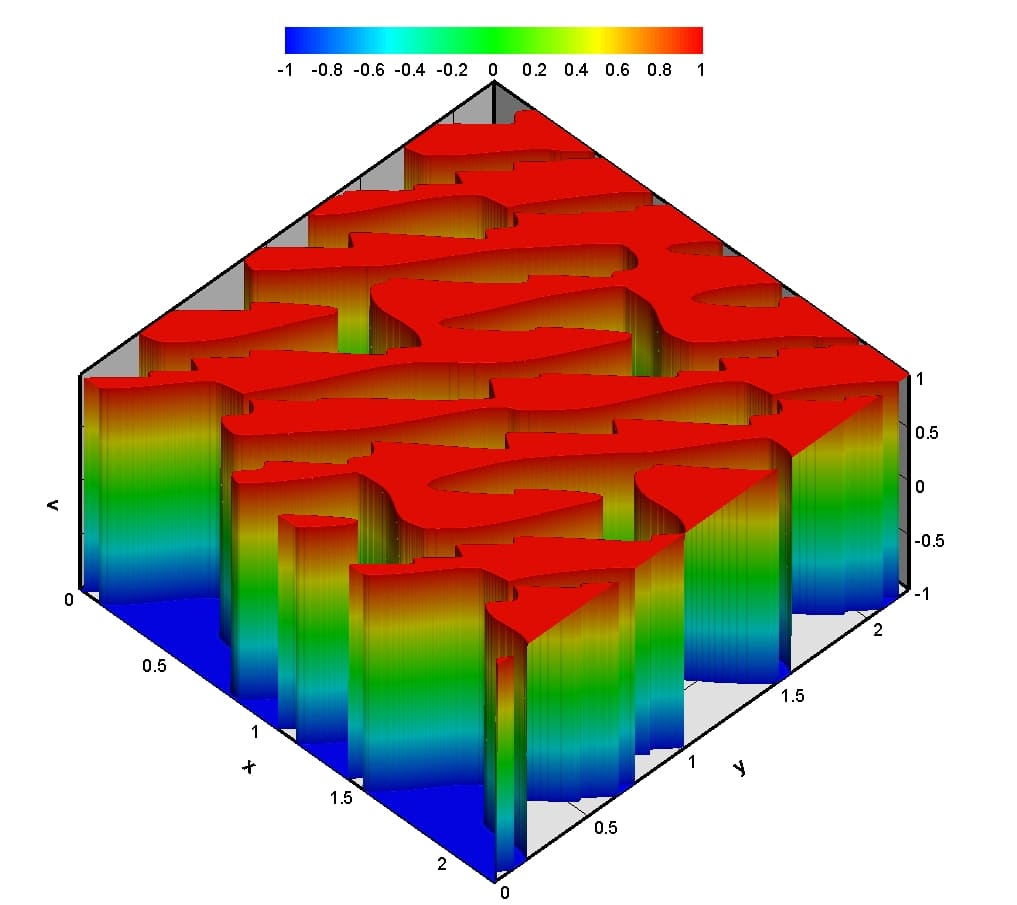}}
  
\caption{3-D contours of velocity $v$ with an $[800\times800]$ grid and $\kappa = 2$ at different times.
(a) $\tau = -10^{-1}$.
 (b) $\tau = -10^{-4}$.}
\label{fig:co3D}
\end{figure}

This example shows the dynamics of the $(2+1)$ dimensional cosmological Burgers in a contracting spacetime with the same initial condition \eqref{eq:IC2D}. Three different grid refinements $[200\times200]$, $[400\times400]$, $[800\times800]$ are chosen. 
Similarly to the previous example, stability is chosen based on the CFL condition; however, we notice that when $\tau\rightarrow 0^-$ the time step $\Delta\tau$ becomes close to a constant number due to the fact that $v(\tau, x,y)\rightarrow \pm1$; hence, $\tau$ cannot asymptotically approach zero. Therefore, the following stability condition is used to compute smaller $\Delta\tau$ as $\tau$ advances to $0$:
\be
\Delta\tau_n < \frac{\tau_n}{\tau_{n - 1}} \Delta\tau_{n - 1}.
\label{eq:StabilityEx-trois}  
\ee
The model is solved with many different conditions and schemes (different spatial and temporal discretization); however, it is noticed that at a very small $\tau$, the solution $v_{j,k}$ goes slightly (order of $10^{-5}$, $10^{-6}$, or smaller, depending on the grid refinement and scheme) above $1$ and below $-1$ where shocks exist. The solution approach to $\pm 1$ as the sign of the source term changes subsequently. Since stability condition \eqref{eq:StabilityEx} is too strong, we propose the following ones for this test ($\Delta x = \Delta y$):
\begin{align}
&\text {When } \kappa \le 1,  \quad \quad \Delta \tau_n \le\frac{1}{2}
\left(
\frac{\Delta y}{\max_{j}\left\{ \vert v_{j}^n\vert \right\}}, 2\vert \tau_n \vert \right).
\\
&\text {When } \kappa > 1  \quad  \quad \Delta \tau_n 
\le\frac{1}{2\kappa}\left(
\frac{\Delta y}{\max_{j}\left\{ \vert v_{j}^n \vert \right\}}, 2\vert \tau_n \vert  \right).
\label{eq:Riemann 12-deux}          
\end{align}
The solutions $v(\tau, x, y)$ converge to $\pm 1$, as is seen in Figures~\ref{fig:co2D} and~\ref{fig:co3D} in 2-D and 3-D, respectively.


\begin{figure}[htbp]
\centering
  \subcaptionbox{1S1T: $L^1 = 8.5124 \times 10^{-10}$}{\includegraphics[height=2.2in]{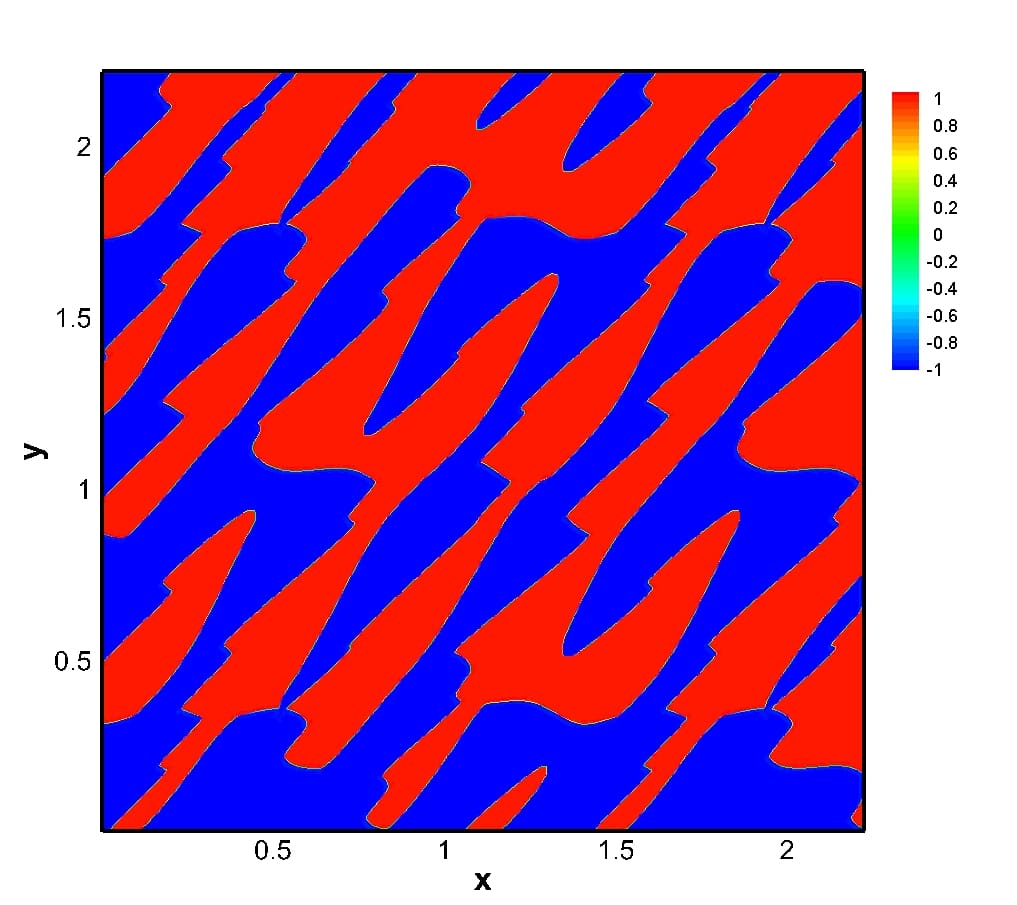}}
  \subcaptionbox{1S3T: $L^1 = 2.2043 \times 10^{-11}$}{\includegraphics[height=2.2in]{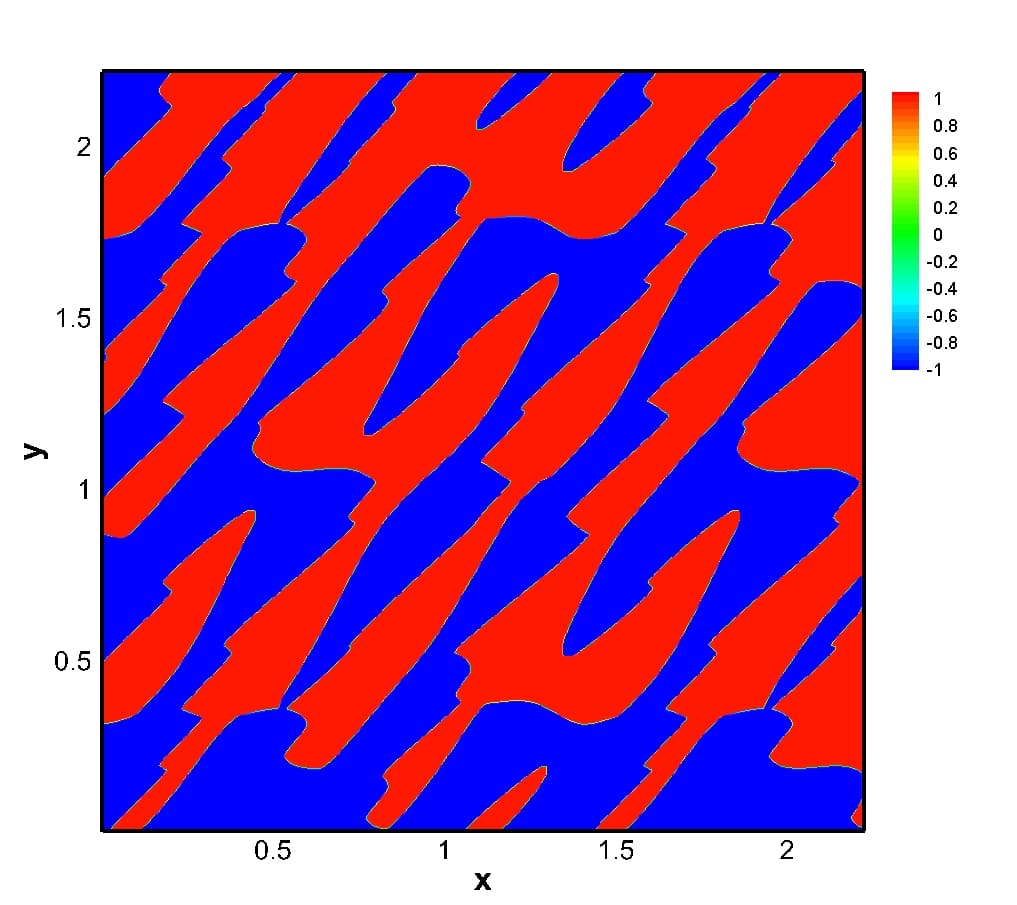}}
  
\centering
  \subcaptionbox{2S1T: $L^1 = 8.5011 \times 10^{-10}$}{\includegraphics[height=2.2in]{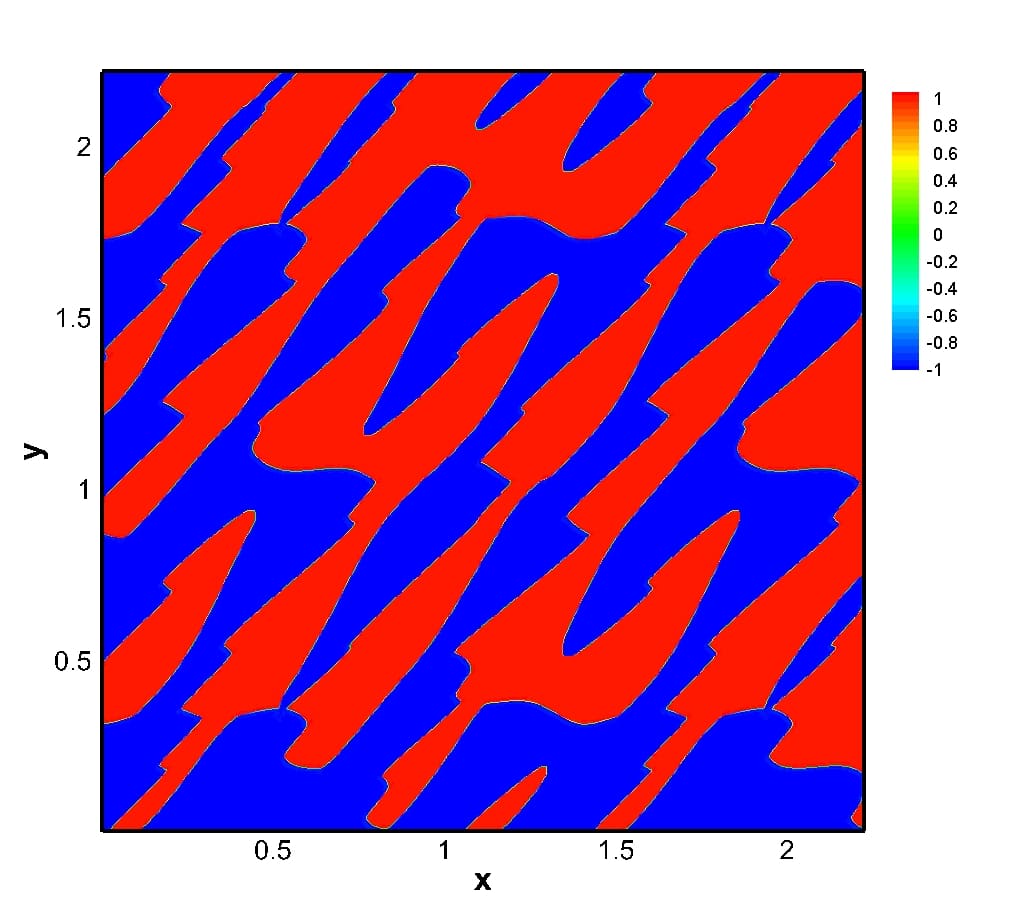}}
  \subcaptionbox{2S3T: $L^1 = 0$}{\includegraphics[height=2.2in]{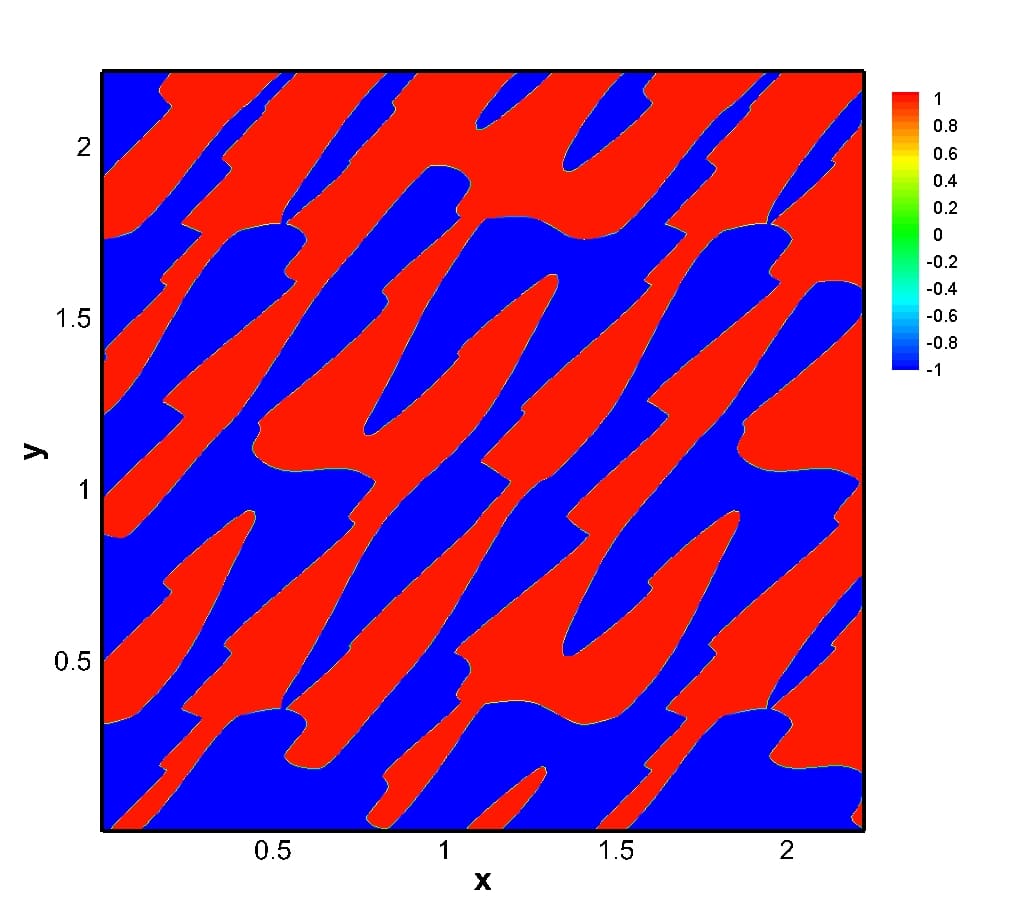}}
   
\caption{2-D contours of the velocity $v$ with an $[800\times800]$ grid, $\kappa = 2$ at $\tau = -10^{-8}$ with different orders of spatial and temporal discretization.}
\label{fig:co2D-SCH}
\end{figure}

The cosmological Burgers equation is solved with a second-order spatial and third-order temporal discretization (2S3T). Observe that a third-order SSP Runge-Kutta scheme for the temporal discretization is selected. In addition, we compute this example with several schemes: 
first-order space and third-order time (1S3T), 1S1T, and 2S1T. Figure~\ref{fig:co2D-SCH} shows the solutions for the above-mentioned schemes. 
The $L^1$ norm for these schemes based on the best scheme (2S3T) is also calculated. Similar to the expanding tests, the $L^1$ norm is very small. Higher-order temporal schemes are more accurate than higher-order spatial ones. Low-order schemes can be used to be able to reduce computational cost.


\begin{figure}[htbp]
\centering
  \subcaptionbox{1S1T: $L^1 = 4.2872 \times 10^{-7}$}{\includegraphics[height=2.2in]{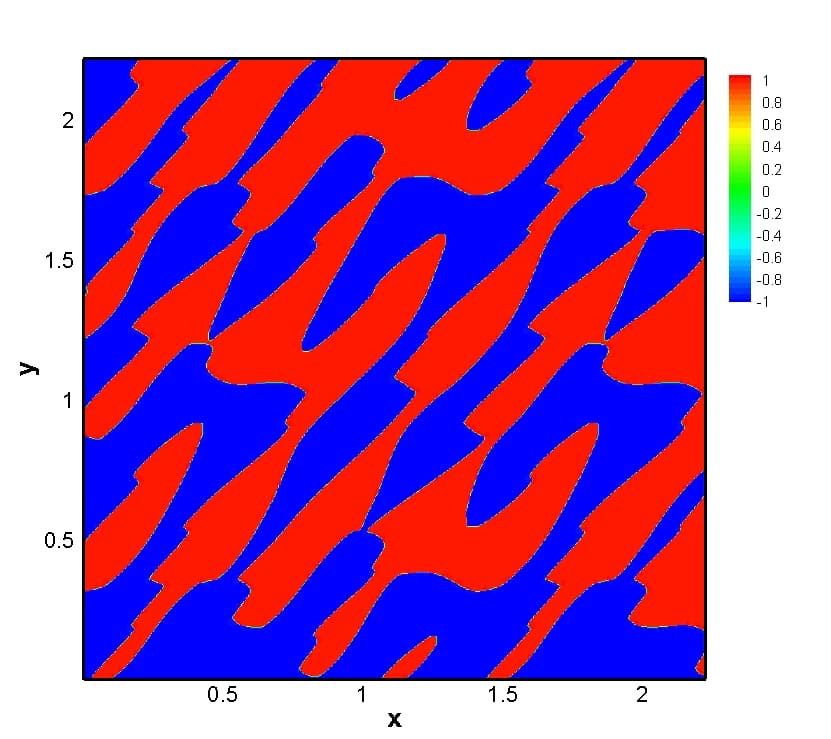}}
  \subcaptionbox{1S3T: $L^1 = 3.8952 \times 10^{-11}$}{\includegraphics[height=2.2in]{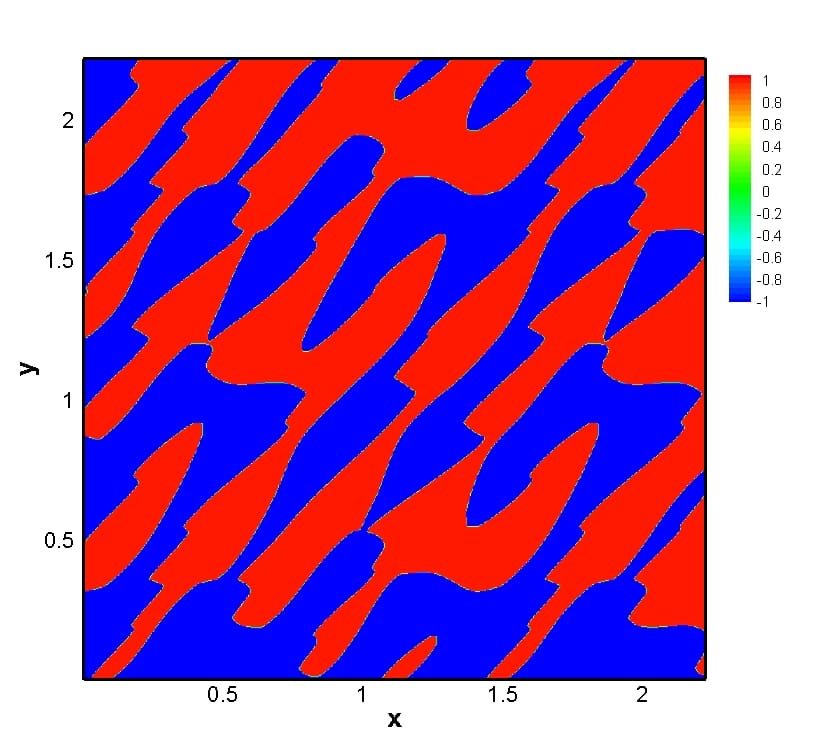}}
  
\centering
  \subcaptionbox{2S1T: $L^1 = 2.2605 \times 10^{-7}$}{\includegraphics[height=2.2in]{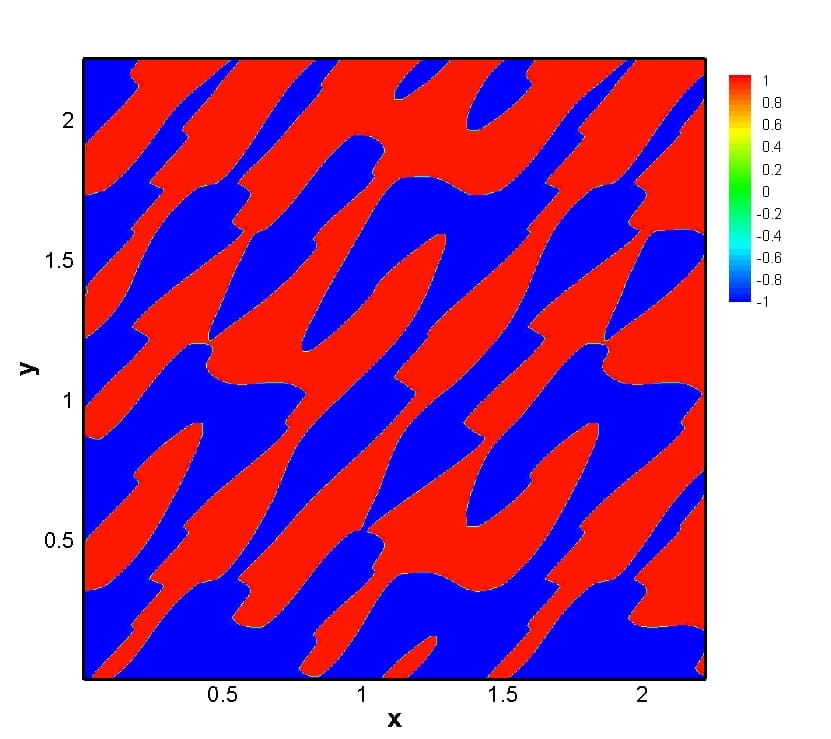}}
  \subcaptionbox{2S3T: $L^1 = 0$}{\includegraphics[height=2.2in]{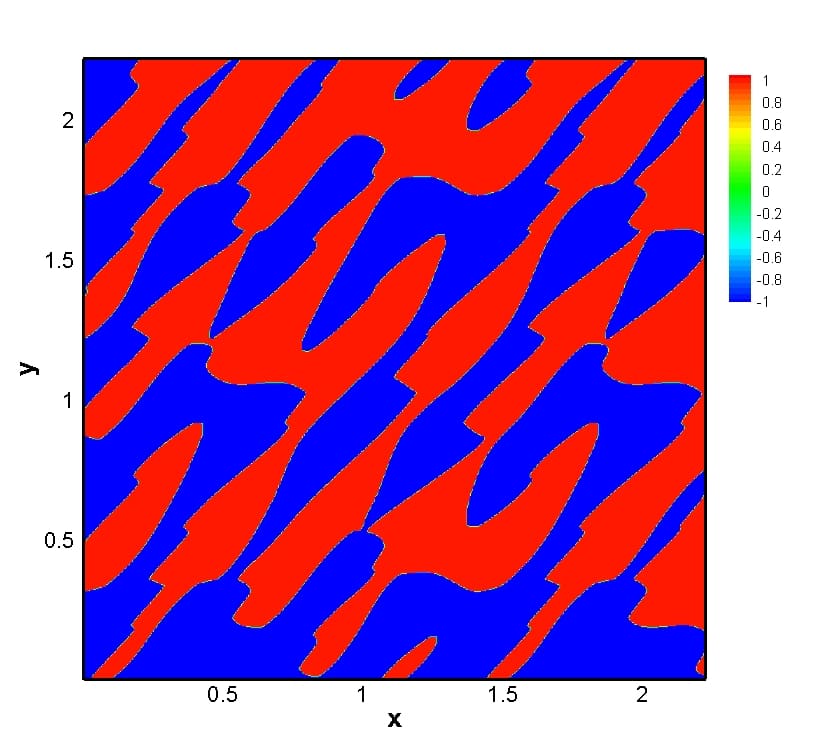}}
   
\caption{2-D contours of the velocity $v$ with an $[800\times800]$ grid, $\kappa = 4$ at $\tau = -10^{-8}$ with different orders of spatial and temporal discretization.}
\label{fig:cok4}
\end{figure}

\begin{figure}[htbp]
\centering
  \subcaptionbox{$\tau = -10^{-1}$}{\includegraphics[height=2.2in]{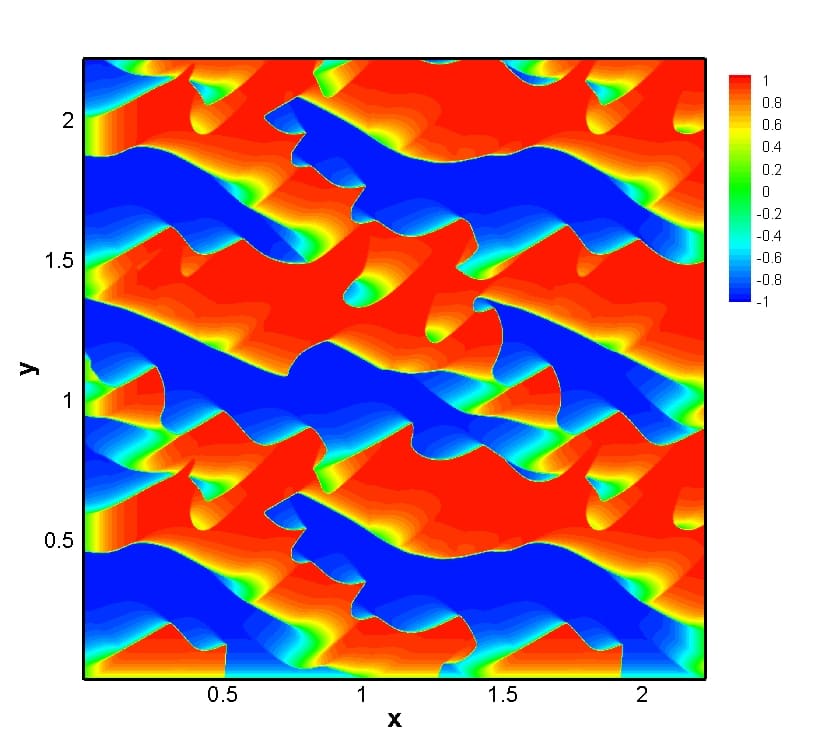}}
  \subcaptionbox{$\tau = -10^{-8}$}{\includegraphics[height=2.2in]{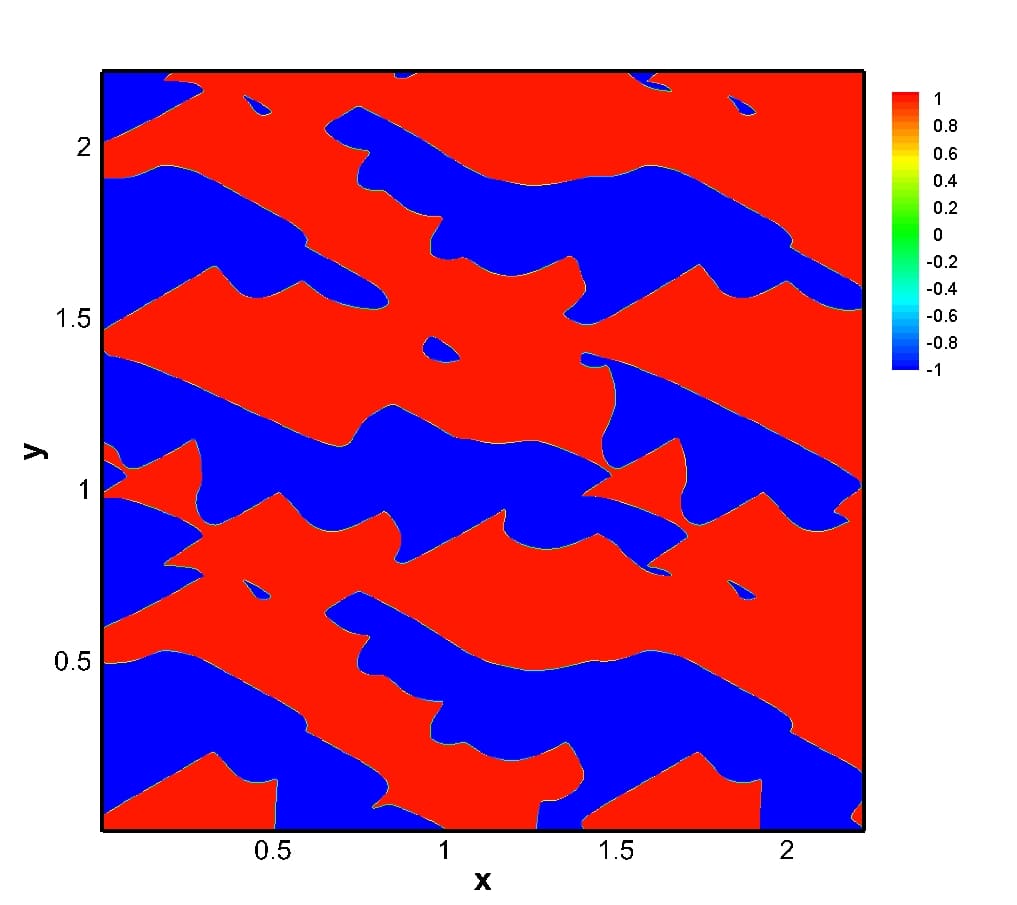}}
  
\centering
  \subcaptionbox{$\tau = -10^{-1}$}{\includegraphics[height=2.2in]{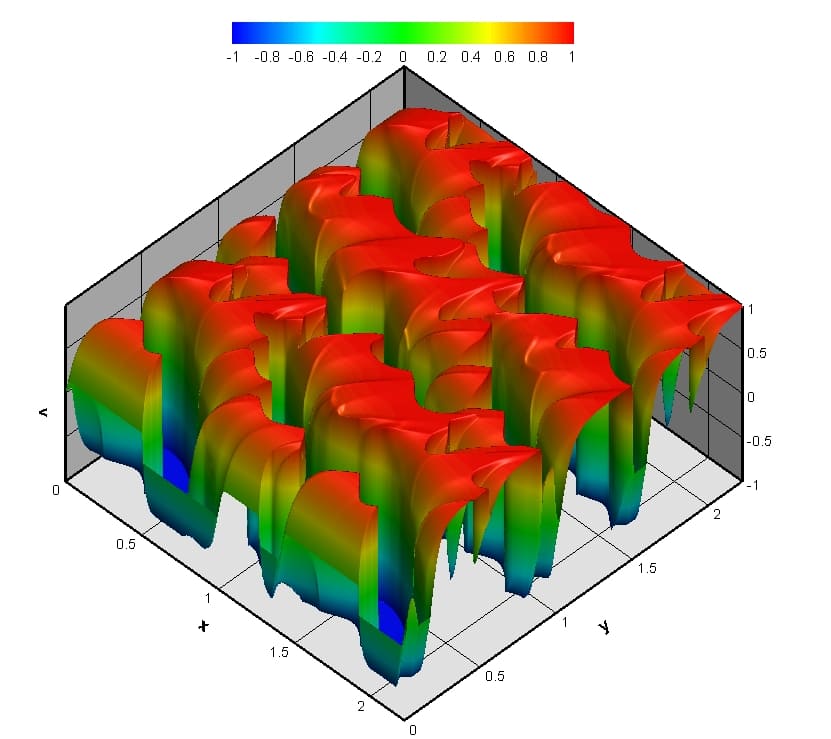}}
  \subcaptionbox{$\tau = -10^{-8}$}{\includegraphics[height=2.2in]{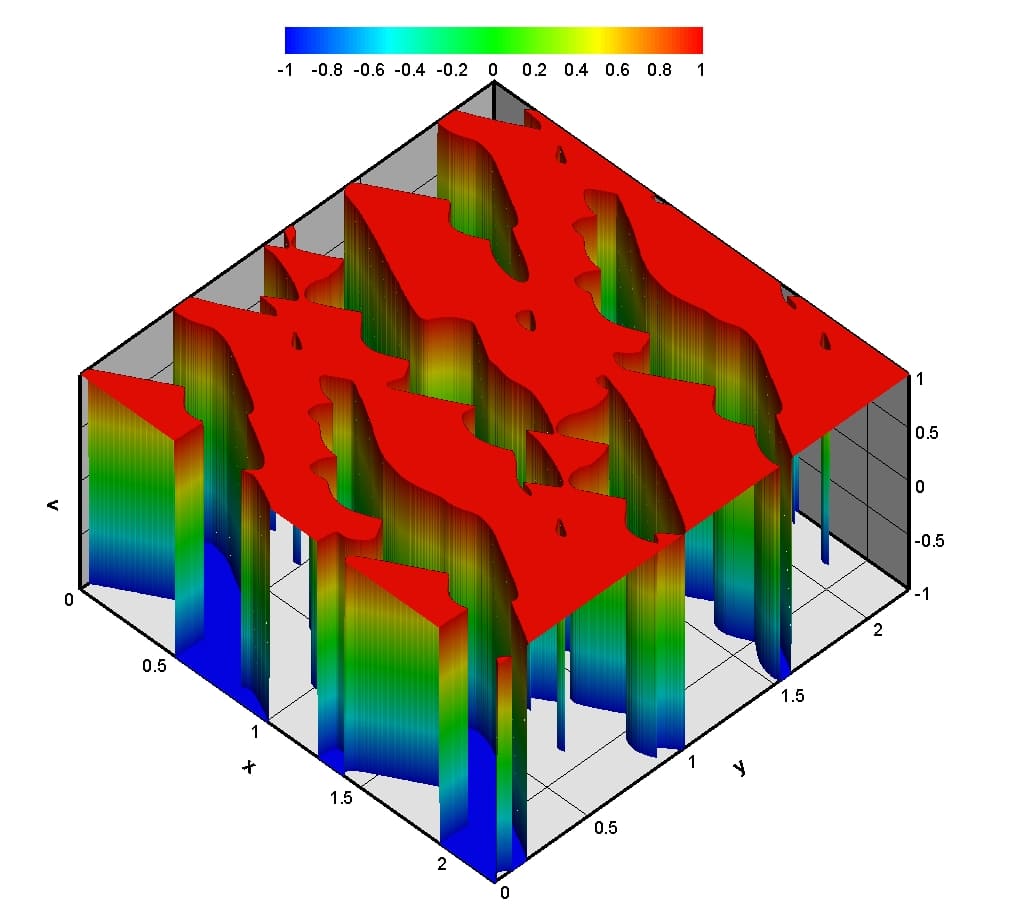}}
   
\caption{2D and 3D contours of the solution with flux  $g(v) = v^3 /2$ at two different times.}
\label{fig:co-g=u31}
\end{figure}

\begin{figure}[htbp]
\centering
  \subcaptionbox{$\tau = -10^{-1}$}{\includegraphics[height=2.2in]{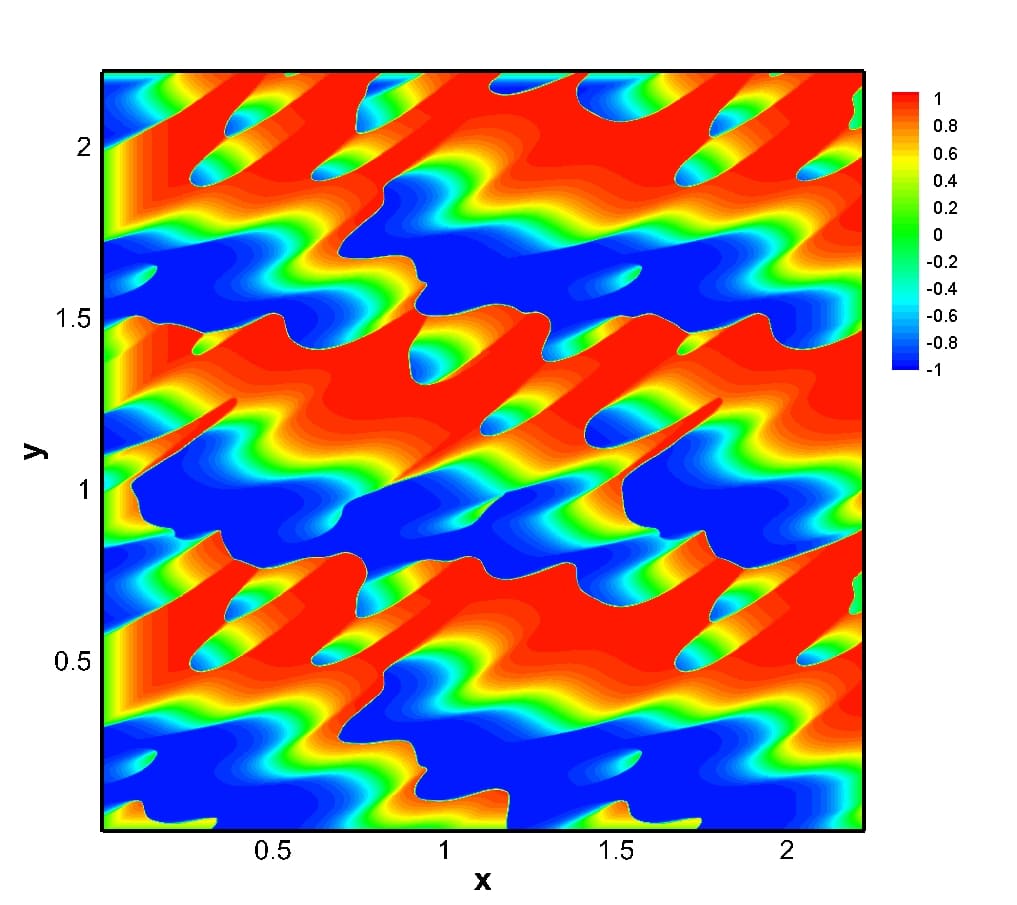}}
  \subcaptionbox{$\tau = -10^{-8}$}{\includegraphics[height=2.2in]{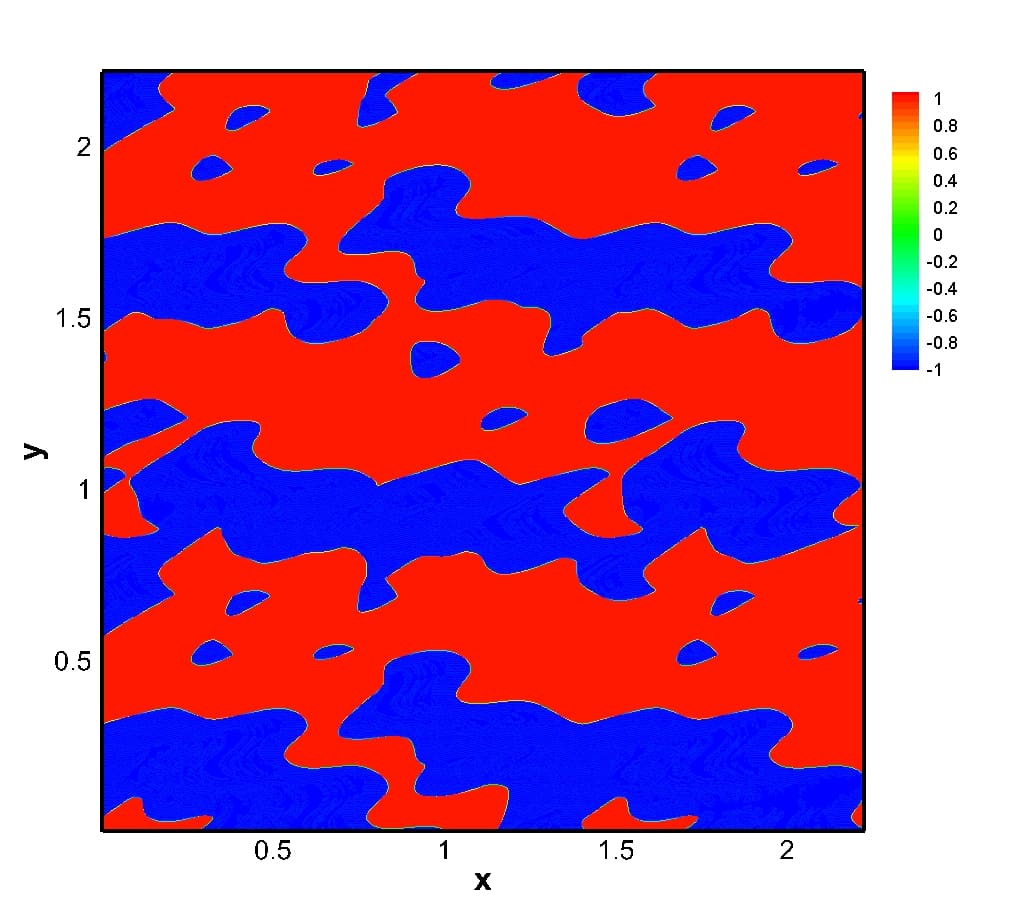}}
  
\centering
  \subcaptionbox{$\tau = -10^{-1}$}{\includegraphics[height=2.2in]{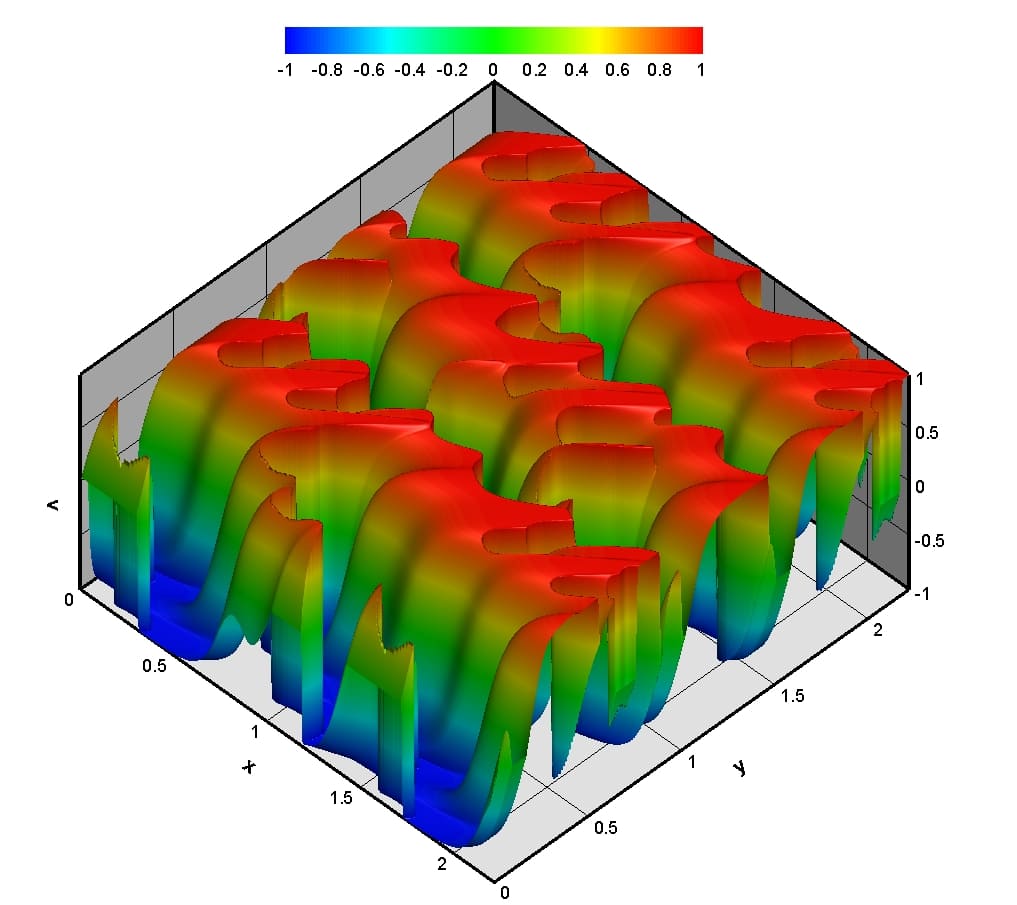}}
  \subcaptionbox{$\tau = -10^{-8}$}{\includegraphics[height=2.2in]{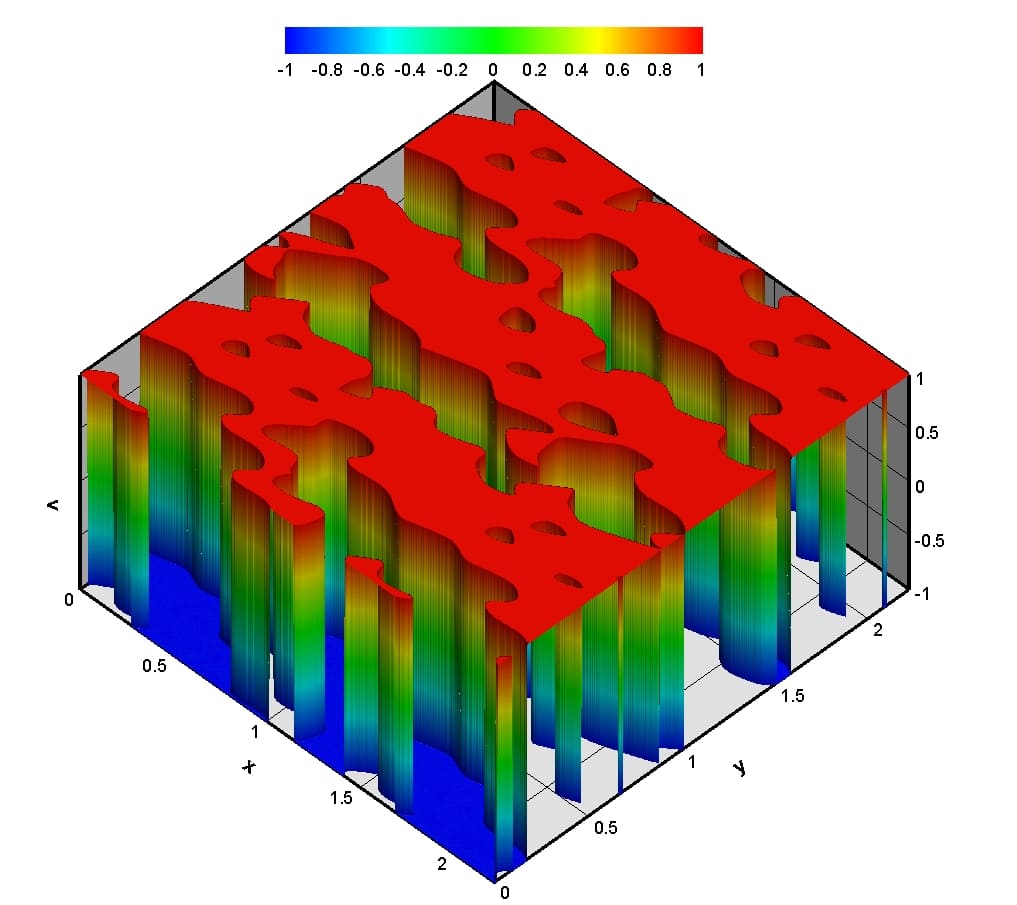}}
   
\caption{2D and 3D contours of the solution with flux $g(v) = (1-\beta) v^2/2 + \beta v^3/3$  and $\beta=1/2$ at two different times.}
\label{fig:co-g=u32}
\end{figure}

Similarly as in the previous tests, the contracting test is solved with $\kappa = 4$. The results show that velocity approaches to $\pm 1$. Figure~\ref{fig:cok4} provides the velocity with four different schemes. The $L^1$ norm for these solutions is very small and the order of temporal discretization is more effective than the order of spatial one. 
Moreover, the cosmological Burgers equation with fluxes in \eqref{eq:gv3} and \eqref{eq:gv32} for a contracting background are solved. $v(\tau, x, y)$ advances to $\pm1$ (Figures~\ref{fig:co-g=u31} and~\ref{fig:co-g=u32}) and shocks are created towards the $x$ direction. 


\subsection{Comparison between the $(1+1)$ and $(2+1)$ models} 
\label{Convergence}

\begin{figure}[htbp]
\centering
  \subcaptionbox{}{\includegraphics[height=2.2in]{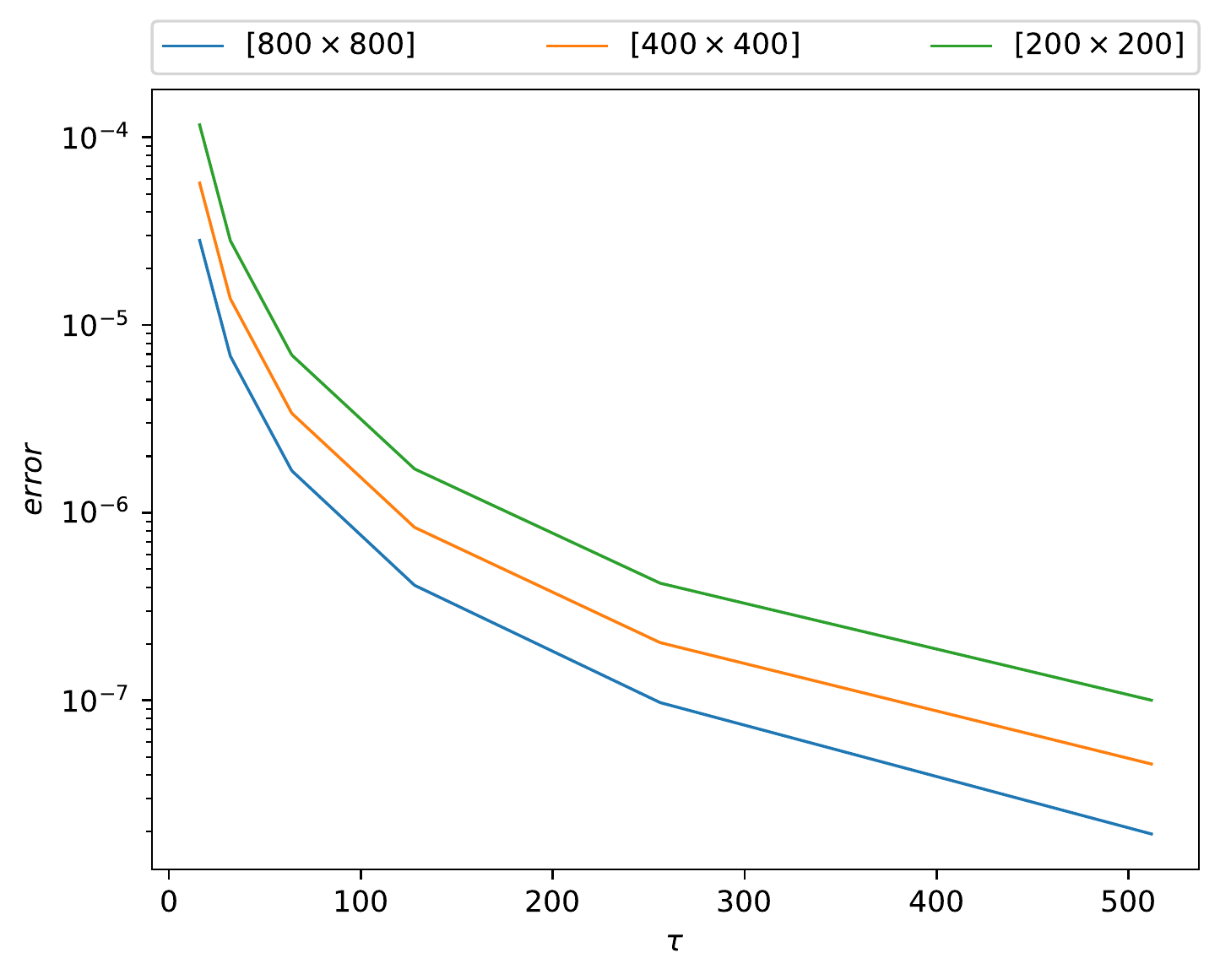}}
  \subcaptionbox{}{\includegraphics[height=2.2in]{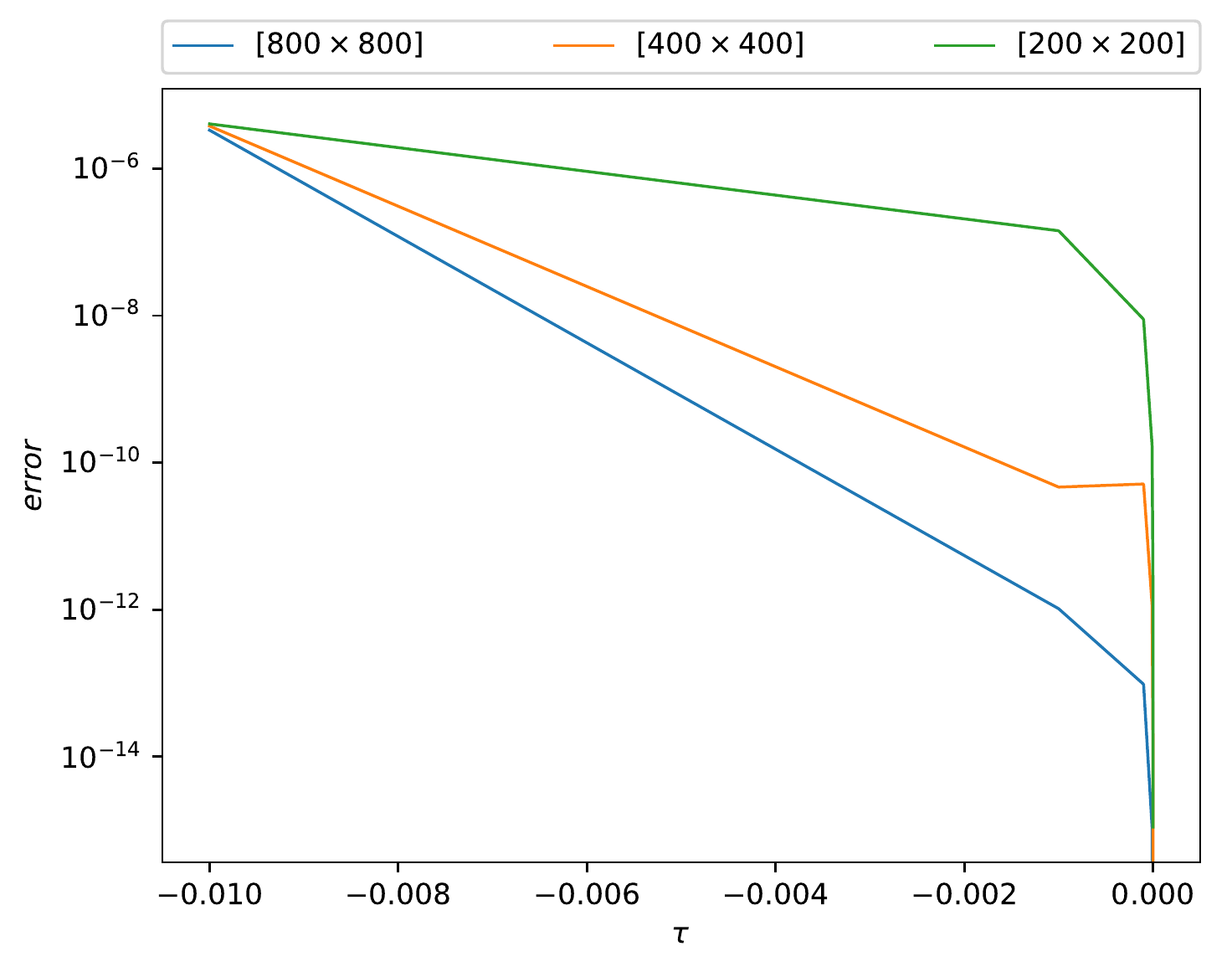}}
\caption{Convergence of the solution $v$ at different times $\tau$ with $\kappa = 2$.
 (a) Expanding background. (b) Contracting background.}
\label{fig:converg}
\end{figure}

\begin{figure}[htbp]
\centering
 \subcaptionbox{$\kappa = 2$}{\includegraphics[height=2.2in]{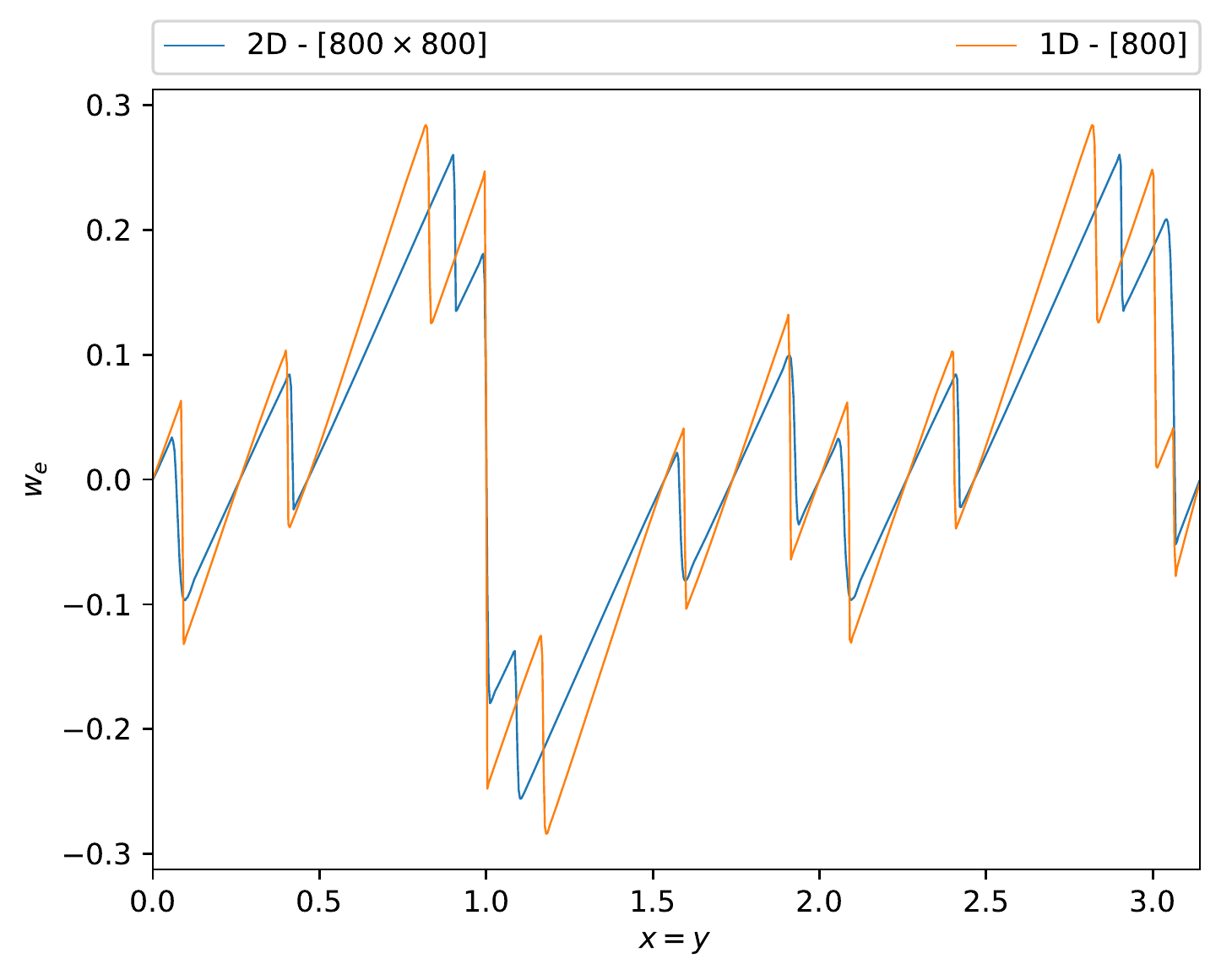}}
 \subcaptionbox{$\kappa = 4$}{\includegraphics[height=2.2in]{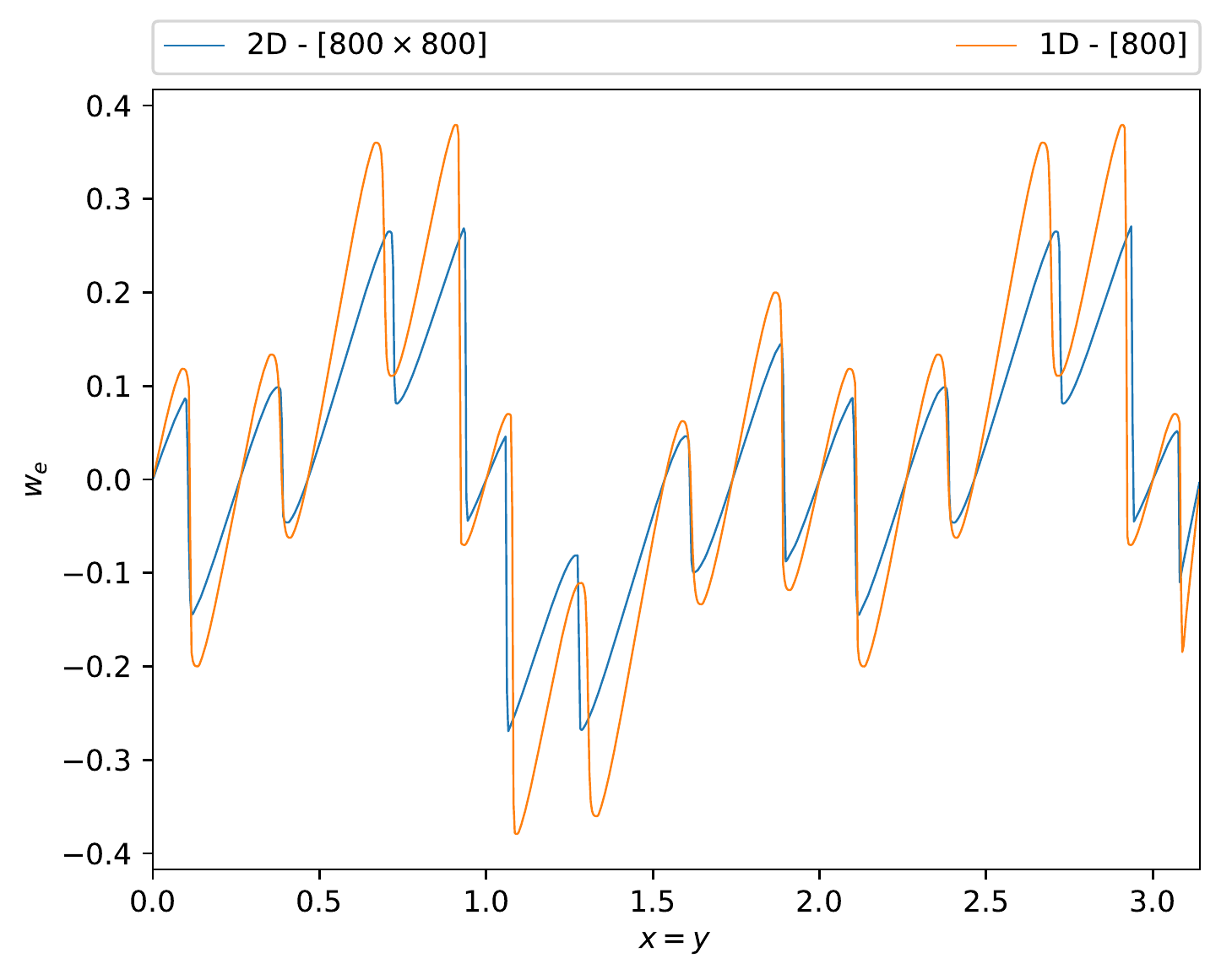}}

\caption{Solution based on the scheme (2T4S) on an expanding background at $\tau = 1024$ in $(2+1)$ dimensions on $x=y$ and in $(1+1)$ dimensions.}
\label{fig:comparisonEx}
\end{figure}

\begin{figure}[htbp]
\centering
 \subcaptionbox{$\kappa = 2$}{\includegraphics[width=2.8in]{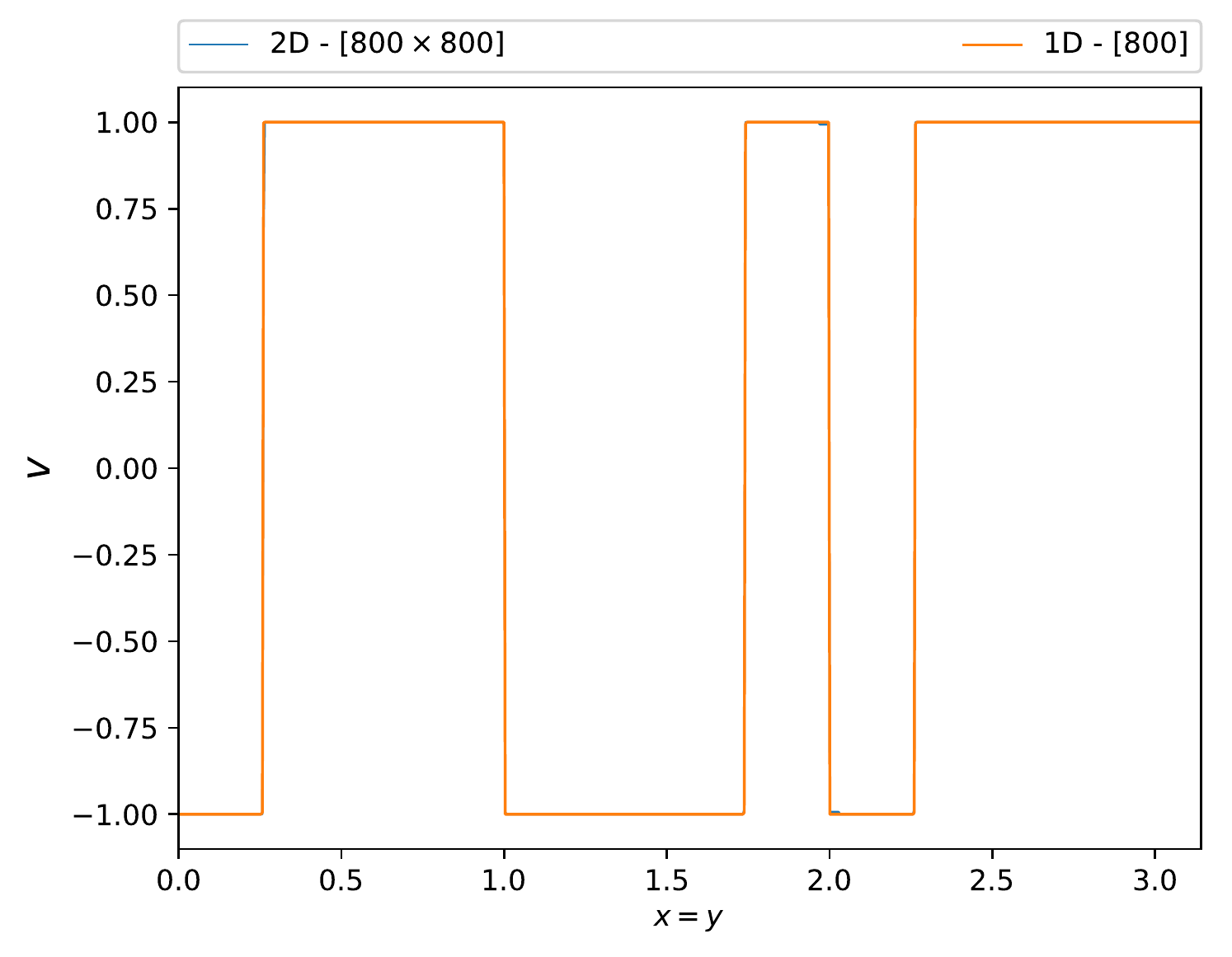}}
 \subcaptionbox{$\kappa = 4$}{\includegraphics[width=2.8in]{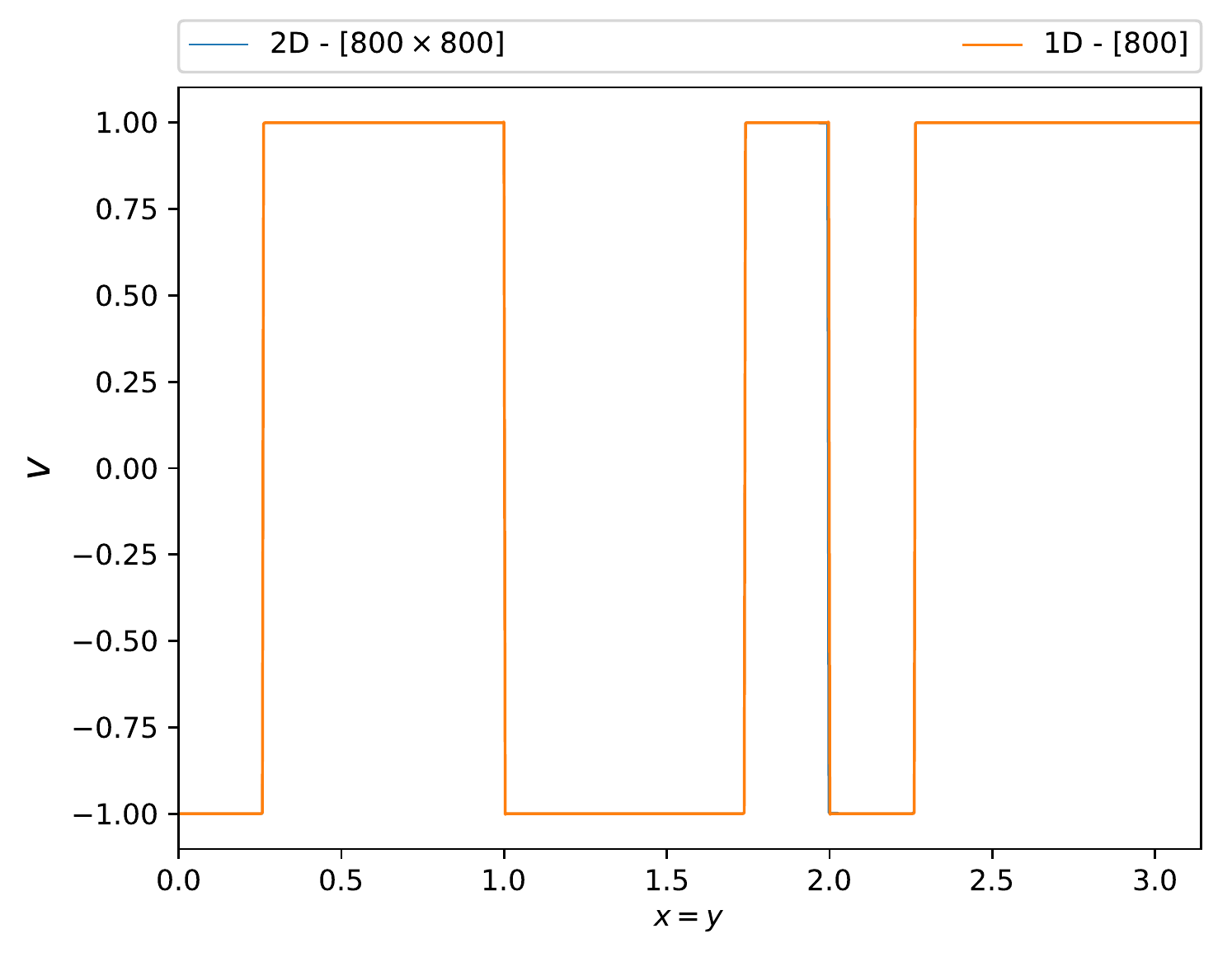}}

\caption{Solutions based on the  scheme (2S3T) on a contracting background in $(2+1)$ dimensions on $x=y$ and $(1+1)$ dimensions, at the time $\tau = -10^{-4}$.}
\label{fig:comparisonCo}
\end{figure}

To study the convergence of the solutions, we compute the $L^1$ norm with different grid refinements at $\tau=16$, $32$, $64$, $128$, $256$, $512$, $1024$ for the expanding background tests, and 
$\tau=-10^{-2}$, $- 10^{-3}$, $- 10^{-4}$, $- 10^{-5}$, $- 10^{-6}$, $- 10^{-7}$, 
and $- 10^{-8}$ for the contracting background. 
Figure~\ref{fig:converg} shows that the error decreases as grid refines in both tests. 
The solution $v(\tau, x, y)$ where $x=y$ in $(2+1)$ dimensions ($[800 \times 800]$) are compared with the corresponding solutions  $v(\tau, x)$ in $(1+1)$ dimensions with $\kappa = 2$ and $4$. The rescaled solutions for the expanding case presented in Figure~\ref{fig:comparisonEx} show some differences between the $(1+1)$-- and $(2+1)$---solutions due to the fact that the initial condition is not quite symmetric with respect to the diagonal of the domain $x=y$. However, both solutions follow a similar evolution.
Finally, the numerical solutions in the contracting case are presented in Figure~\ref{fig:comparisonCo}, and 
we observe that they follow a similar trend. 


\vskip.3cm

\noindent{\bf Acknowledgments.} The authors were supported by an Innovative Training Network (ITN) under the grant 642768 ``ModCompShock'' managed by the third author (PLF). During the preparation of this paper, the second author (MAG) was a six-month visiting fellow at Sorbonne University. The third author (PLF) was also partially supported by the Centre National de la Recherche Scientifique (CNRS). 


\end{document}